\documentclass[reqno, a4paper, 11pt]{amsart}
\usepackage{setspace}
\usepackage{fullpage}
\usepackage{amsmath}
\usepackage{amsaddr}
\usepackage{amsthm}
\usepackage{amsfonts}
\usepackage{amssymb}
\usepackage{graphicx,caption,epstopdf}
\usepackage{mathrsfs}
\usepackage{xcolor}
\usepackage{fullpage}
\usepackage[mathlines]{lineno}

\usepackage{etoolbox}

\newcommand*\linenomathpatchAMS[1]{
  \expandafter\pretocmd\csname #1\endcsname {\linenomathAMS}{}{}
  \expandafter\pretocmd\csname #1*\endcsname{\linenomathAMS}{}{}
  \expandafter\apptocmd\csname end#1\endcsname {\endlinenomath}{}{}
  \expandafter\apptocmd\csname end#1*\endcsname{\endlinenomath}{}{}
}

\expandafter\ifx\linenomath\linenomathWithnumbers
  \let\linenomathAMS\linenomathWithnumbers

  \patchcmd\linenomathAMS{\advance\postdisplaypenalty\linenopenalty}{}{}{}
\else
  \let\linenomathAMS\linenomathNonumbers
\fi

\let\oldtocsection=\tocsection
 
\let\oldtocsubsection=\tocsubsection
 
\let\oldtocsubsubsection=\tocsubsubsection
 
\renewcommand{\tocsection}[2]{\hspace{0em}\oldtocsection{#1}{#2}}
\renewcommand{\tocsubsection}[2]{\hspace{1em}\oldtocsubsection{#1}{#2}}
\renewcommand{\tocsubsubsection}[2]{\hspace{2em}\oldtocsubsubsection{#1}{#2}}

\linenomathpatchAMS{gather}
\linenomathpatchAMS{multline}
\linenomathpatchAMS{align}
\linenomathpatchAMS{alignat}
\linenomathpatchAMS{flalign}

\usepackage{comment}
\usepackage[english,algoruled,lined,noresetcount,norelsize]{algorithm2e}
\usepackage{amsmath,amsthm,amssymb,enumerate,mathscinet}
\usepackage{fullpage}
\usepackage{graphicx}

\usepackage{mathrsfs}  
\usepackage[backref,colorlinks]{hyperref} 

\usepackage[abbrev,msc-links]{amsrefs}   
\usepackage[utf8]{inputenc}
\usepackage{tikz}

\usepackage{enumerate}
\usepackage{cleveref}     

\usepackage{mathtools}

\usepackage{lineno}
\theoremstyle{plain}
\newtheorem{thm}{Theorem}[section]
\newtheorem{fact}[thm]{Fact}
\newtheorem{prop}[thm]{Proposition}
\newtheorem{clm}[thm]{Claim}
\newtheorem{cor}[thm]{Corollary}
\newtheorem{lem}[thm]{Lemma}

\newtheorem*{clm*}{Claim}

\theoremstyle{definition}

\newtheorem{dfn}[thm]{Definition}

\newtheorem{conj}[thm]{Conjecture}
\newtheorem{obs}[thm]{Observation}
\usepackage{mathrsfs}

\numberwithin{equation}{section}

\newenvironment{claimproof}[1]{\par\noindent\underline{Proof of Claim:}\space#1}{\hfill $\blacksquare$}

\newenvironment{altclaimproof}[1]{\par\noindent\underline{Proof of Claim \ref{clm:fixingdiamondtrees}:}\space#1}{\hfill $\blacksquare$}

\def\COMMENT#1{}
\let\COMMENT=\footnote

\usepackage{bm}

\usepackage{amsmath}

\usepackage{enumitem}

\DeclareMathOperator{\codeg}{codeg}

\let\polishlcross=\l
\def\l{\ifmmode\ell\else\polishlcross\fi}

\let\eps=\varepsilon
\let\theta=\vartheta
\let\phi=\varphi

\def\EE{\mathbb E}
\def\PP{\mathbb P}

\def\NN{\mathbb N}

\def\RR{\mathbb R}
\def\AA{\mathbb A}
 
\def\cB{{\mathcal B}}
\def\cC{{\mathcal C}}
\def\cD{{\mathcal D}}

\def\cR{{\mathcal R}}
\def\cS{{\mathcal S}}
\def\cT{{\mathcal T}}

\def\cW{{\mathcal W}}

\DeclareMathAlphabet\bc{OMS}{cmsy}{b}{n}

\DeclareMathAlphabet\bfc{OMS}{cmsy}{b}{n}
\DeclareMathAlphabet{\pzc}{OT1}{pzc}{m}{it}

\begin{document}

\title{Clique factors in pseudorandom graphs}
\date{\today}

\author{Patrick Morris}
\address{Universitat Polit\`ecnica de Catalunya (UPC), Barcelona, Spain}
\email{pmorrismaths@gmail.com}
\thanks{The  author was supported by the Deutsche Forschungsgemeinschaft (DFG, German Research
Foundation) under Germany's Excellence Strategy - The Berlin Mathematics
Research Center MATH+ (EXC-2046/1, project ID: 390685689) and by a Walter Benjamin fellowship of the DFG - project number 504502205. }

\date{\today}

\begin{abstract}
An $n$-vertex graph is said to to be \emph{$(p,\beta)$-bijumbled} if for any vertex sets $A,B\subseteq V(G)$, we have
\[e(A,B)=p|A||B|\pm \beta \sqrt{|A||B|}.\]
 We prove that for any $3\leq r\in \NN$ and $c>0$ there exists an  $\eps>0$ such that any $n$-vertex $(p,\beta)$-bijumbled graph with $n\in r \mathbb{N}$, $p>0$,  $\delta(G)\geq cpn$ and  $\beta \leq \eps p^{r-1}n$, contains a~$K_r$-factor. This implies a corresponding result for the stronger pseudorandom notion  of $(n,d,\lambda)$-graphs.

For the case of triangle factors, that is when  $r=3$, this result  resolves a conjecture of Krivelevich, Sudakov and Szab\'o from 2004 and it is tight due to a pseudorandom triangle-free construction of Alon. In fact,  in this case even more is  true: as a corollary to this result and a result of Han, Kohayakawa, Person and the author, we can conclude that the same condition of $\beta=o(p^2n)$ actually guarantees  that a $(p,\beta)$-bijumbled graph $G$ contains every graph on $n$ vertices with maximum degree at most 2.

\end{abstract}

\maketitle
 
\newpage 
 
\tableofcontents

\newpage

\section{Introduction} \label{sec:intro}

We say a graph $G$ contains a \emph{$K_r$-factor} if there is a collection of vertex disjoint copies of $K_r$ that completely cover the vertex set of $G$. When $r=3$, we often refer to a $K_3$-factor as a \emph{triangle factor}.   As a natural generalisation of a perfect matching in a graph, $K_r$-factors  are a fundamental object in  
graph theory with a wealth of results studying various aspects and variants, in particular exploring  probabilistic e.g. \cite{balogh2012corradi,han2019tilings,JohanssonKahnVu_FactorsInRandomGraphs,kim2003perfect,kri97}, extremal  e.g. \cite{allen2015density,balogh2016triangle,nenadov2018ramsey,treglown2016degree},  and algorithmic e.g.  \cite{caprara2002packing,kann1994maximum,kierstead2010fast,kirkpatrick1983complexity},  considerations. 
However, unlike perfect matchings, it is not easy to verify whether a graph $G$ contains a $K_r$-factor or not. Certainly it is necessary that the number of vertices of  $G$ must be divisible by $r$ but given this,  it was proved by Schaeffer~\cite{karp1975computational} (in the case r = 3) and by Kirkpatrick and Hell~\cite{kirkpatrick1983complexity} (in general)  that determining if a graph on $n\in r \mathbb{N}$ vertices contains a $K_r$-factor is an NP-complete problem. Given that we cannot hope for a nice characterisation of graphs which contain $K_r$-factors, there has been a focus on providing sufficient conditions which are computationally easy to verify. One classical such theorem is due to Hajnal and Szemer\'edi \cite{hs} who showed that a $K_r$-factor is guaranteed if the host graph is sufficiently dense. The case of triangle factors was previously shown by Corr\'adi and Hajnal~\cite{corradi1963maximal}.

\begin{thm} \label{thm:HajnalSzem}
If $ 3\leq r\in \NN$ and $G$ is a graph on $n\in r \mathbb{N}$ vertices with minimum degree $\delta(G)\geq (1-1/r)n$, then $G$ contains a $K_r$-factor. 
\end{thm}
This theorem is tight, as can be seen, for example, by taking $G$ to be a complete graph with a clique of size $n/r+1$ removed to leave an independent set of vertices, say $I$. One then has that $\delta(G) = (1-1/r)n-1$ and $G$ does not have a $K_r$-factor. Indeed, any copy of $K_r$ in a family of vertex disjoint $K_r$s can use at most one vertex of $I$ but a $K_r$-factor should contain $n/r<|I|$ copies of $K_r$. All examples verifying the tightness of Theorem \ref{thm:HajnalSzem} share some features with the graph given here. For example they contain much larger independent sets than almost all graphs of this density. Therefore, one might hope to capture more graphs having a $K_r$-factor by adding a condition that precludes the atypical behaviour of the extremal examples. 

\vspace{2mm}
This naturally leads us to the notion of \emph{pseudorandom graphs}, which are, roughly speaking, graphs which imitate random graphs of the same density. The study of pseudorandom graphs, initiated in the 1980s by Thomason \cite{Th87a,Th87b}, has become a central and vibrant field at the intersection of Combinatorics and Theoretical Computer Science. We refer to the excellent survey of Krivelevich and Sudakov \cite{KS06} for an introduction to the topic.  One way of imposing pseudorandomness is through the spectral notion of the  eigenvalue gap. This then leads to the study of $(n,d, \lambda)$-graphs $G$ which are $d$-regular $n$-vertex graphs with \emph{second eigenvalue} $\lambda$. By second eigenvalue, what is actually meant is the  second largest eigenvalue in absolute value as follows.  Given an $n$-vertex $d$-regular graph $G$, we can look at the eigenvalues of the adjacency  matrix $A$ of $G$ which, as $A$ is a symmetric $0/1$-matrix, are real and can be ordered  as $\lambda_1\geq \lambda_2 \geq \ldots \geq \lambda_n$. The second eigenvalue is then defined to be $\lambda:=\max\{|\lambda_2|,|\lambda_n|\}$. 
It turns out that this parameter $\lambda$ controls the pseudorandomness of the graph $G$,  with smaller values of $\lambda$ giving graphs that have stronger pseudorandom properties.  More concretely, the relation is given by the following property of $(n,d,\lambda)$-graphs, see e.g.~\cite[Theorem 2.1]{KS06},  which is known as the \emph{Expander Mixing Lemma} and shows that $\lambda$ controls the edge distribution between vertex sets. For any vertex subsets $A,B$ of an $(n,d,\lambda)$-graph $G$, one has that \begin{equation} \label{eq:ELM} \left|e(A,B)-\frac{d}{n}|A||B|\right|\leq \lambda \sqrt{|A||B|}, \end{equation}
where $e(A,B):=|\{uv\in E(G):u\in A, b\in B\}|$ denotes the number\footnote{Note that edges that lie in $A\cap B$ are counted twice.} of edges in $G$ with one endpoint in  $A$ and the other in $B$. Note that $\frac{d}{n}$ is the density of the graph $G$, and hence one would expect to see $\frac{d}{n}|A||B|$ edges between the vertex sets $A$ and $B$ in a random graph $G$. The pseudorandom parameter $\lambda$ then controls the discrepancy from this paradigm.

    It  follows from simple linear algebra, see e.g. \cite{KS06}, that for an $(n,d,\lambda)$-graph, one has that $\lambda\leq d$ always and moreover, as long as $d$ is not too close to $n$, say $d\leq 2n/3$, one has that $\lambda=\Omega(\sqrt{d})$. Thus, we think of $(n,d,\lambda)$-graphs with $\lambda=\Theta(\sqrt{d})$ as being \emph{optimally pseudorandom}. For example, it is known that random regular graphs are optimally pseudorandom $(n,d,\lambda)$-graphs with high probability\footnote{Here, and throughout, we say that a property holds \emph{with high probability} if the probability that it holds tends to~$1$ as the number of vertices $n$ tends to infinity.}~\cite{broder1998optimal,tikhomirov2019spectral}. 

\vspace{2mm}

A prominent theme in the study of pseudorandom graphs has been to give conditions on the parameters, $n$, $d$ and $\lambda$ that guarantee certain properties of an $(n,d,\lambda)$-graph. For example, it follows easily from \eqref{eq:ELM} that any $(n,d,\lambda)$-graph $G$ with $\lambda< d^2/n$ contains a triangle as there is an edge in the neighbourhood of every vertex. In particular, any optimally pseudorandom graph with $d=\omega(n^{2/3})$ must contain a triangle.  Moreover, this condition is tight due to a  triangle-free construction of an $(n,d,\lambda)$-graph  due to Alon~\cite{Alon94} with $d=\Theta(n^{2/3})$ and $\lambda=\Theta(n^{1/3})$. Alon's construction is optimally pseudorandom and Krivelevich, Sudakov and Szab\'o~\cite{KSS04} generalised it to the whole possible range of densities. That is, for any $d=d(n)$ such that $\Omega(n^{2/3})=d \leq n$, they gave  a sequence of infinitely many $n$ and triangle-free $(n',d,\lambda)$-graphs with $n'=\Theta(n)$ and $\lambda=\Theta(d^2/n)$.  In general, finding optimal conditions for subgraph appearance in $(n,d,\lambda)$-graphs seems very hard. Indeed the only  tight  conditions that are known  are those for fixed size odd cycles~\cite{AK98,KS06}. With respect to spanning structures, it is only perfect matchings that have been well understood~\cite{brouwer2005eigenvalues,cioaba2005perfect,KS06}.  Whilst such questions are interesting in their own right, they also have implications in other areas of mathematics. As an example, we mention the beautiful connection given by Alon and Bourgain~\cite{alon2014additive} (see also~\cite{ABHKP16}) who used the existence of certain subgraphs in pseudorandom graphs to prove the existence of additive patterns in large multiplicative subgroups of finite fields.

The purpose of this paper is to answer what has become one of the central problems in this area, by giving a tight condition for an $(n,d,\lambda)$-graph to contain a triangle factor. 

\begin{thm} \label{thm:maintriangle}
There exists $\eps>0$ such that any $(n,d,\lambda)$-graph with $n\in 3 \mathbb{N}$, $d>0$ and $\lambda\leq \eps d^2/n$, contains a triangle factor. 
\end{thm}
Theorem \ref{thm:maintriangle} was conjectured by Krivelevich, Sudakov and Szab\'o \cite{KSS04} in 2004. Focusing solely on optimally pseudorandom graphs, that is, setting $\lambda=\Theta(\sqrt{d})$, Theorem \ref{thm:main} gives that any optimally pseudorandom graph with $d=\omega(n^{2/3})$ contains a triangle factor. Comparing this to Theorem \ref{thm:HajnalSzem}, we see that imposing pseudorandomness, which is easy to compute via the second eigenvalue, allows us to capture much sparser graphs which are guaranteed to contain a triangle factor.

\vspace{2mm}

Theorem~\ref{thm:maintriangle} (and the more general Theorem~\ref{thm:main} below) conclude a body of work towards the conjecture of Krivelevich, Sudakov and Szab\'o and the proof of the theorem, discussed in Section~\ref{sec:proofmain}, builds upon the many beautiful ideas of various authors, which have arisen in this study. The first step towards the conjecture was given by Krivelevich, Sudakov and Szab\'o \cite{KSS04} themselves, who showed that $\lambda \leq \eps d^3/(n^2\log n)$ for some sufficiently small $\eps$ is enough to guarantee a triangle factor. This was improved to $\lambda \leq \eps d^{5/2}/n^{3/2}$  by Allen, B\"ottcher, H\`an, Kohayakawa and Person \cite{ABHKP17}  who also proved that the same condition guarantees the appearance of the square of a Hamilton cycle, a supergraph of a triangle factor. Recently, Nenadov~\cite{Nen18} got very close to the conjecture, showing that $\lambda \leq \eps d^2/(n\log n)$ guarantees a triangle factor. Concentrating solely on optimally pseudorandom graphs, these results  imply that having degree $d=\omega(n^{4/5}(\log n)^{2/5})$, $\omega(n^{3/4})$ and $\omega((n \log n)^{2/3})$ respectively, guarantees the existence of a triangle factor. 

In a different direction, one can fix the condition that $\lambda \leq \eps d^2/n$ for some small $\eps>0$ and prove the existence of other structures giving evidence for a triangle factor. Again, this was initiated by Krivelevich, Sudakov and Szab\'o~\cite{KSS04} who proved that with this condition, one can guarantee the existence of a \emph{fractional triangle factor}. That is, they showed that there is some function $w$ which assigns a weight $w(T)\in[0,1]$ to each triangle $T$ in a pseudorandom graph $G$, such that for every vertex $v\in V(G)$, one has that the sum $\sum_{v\in T}w(T)$ of the  weights of triangles containing $v$ is precisely equal to $1$.  Imposing $\{0,1\}$-weights recovers the notion of a triangle factor and a fractional triangle factor is thus a natural relaxation. Another interesting result of Sudakov, Szab\'o and Vu~\cite{sudakov2005generalization} showed that when we have $\lambda\leq \eps d^2/n$, we have many triangles and these are well distributed in the $(n,d,\lambda)$-graph $G$. Indeed they proved a Tur\'an-type result showing that any triangle-free subgraph of such a graph $G$ must contain at most half the edges of $G$. 
A more recent result due to Han, Kohayakawa and Person \cite{HKP18a,HKP18b} shows that $\lambda \leq \eps d^2/n$ guarantees the existence of a \emph{near triangle factor}; that there are vertex disjoint triangles covering all but $n^{647/648}$ vertices of such an $(n,d,\lambda)$-graph.  

\vspace{2mm}

We will deduce Theorem~\ref{thm:maintriangle} from a more general theorem (Theorem~\ref{thm:main} below) which deals with $K_r$-factors for all $r\geq 3$ and works with a larger class of pseudorandom graphs where we do not restrict solely to regular graphs. Indeed, we will work with  following notion of  \emph{bijumbledness}, whose usage dates back to the original works of Thomason~\cite{Th87a,Th87b},  and whose definition captures the key property of edge distribution, given for $(n,d,\lambda)$-graphs by \eqref{eq:ELM}.

\begin{dfn} \label{def:bijumbled}
Let $n\in\NN$, $p=p(n)\in[0,1]$ and $\beta=\beta(n,p)>0$. An $n$-vertex graph $G=(V,E)$ is \emph{$(p,\beta)$-bijumbled} if for every pair vertex subsets $A,B\subseteq V$, one has that \begin{equation} \label{eq:bijumbled} \left|e(A,B)-p|A||B|\right|\leq \beta \sqrt{|A||B|}. \end{equation}

\end{dfn}
Note that, due to \eqref{eq:ELM}, $(n,d,\lambda)$-graphs are $(d/n,\lambda)$-bijumbled.  As with $(n,d,\lambda)$-graphs, we are interested in finding conditions on the parameters $n$, $p$ and $\beta$, that guarantee the existence of certain subgraphs in  $n$ vertex $(p,\beta)$-bijumbled graphs.  Our main theorem gives conditions for the existence of $K_r$-factors for all $r\geq 3$ in this setting.

\begin{thm} \label{thm:main}
For every $3 \leq r\in \mathbb{N}$ and $c>0$ there exists an $\eps>0$ such that any $n$-vertex $(p,\beta)$-bijumbled graph with $n\in r \mathbb{N}$, $p>0$, $\delta(G)\geq cpn$ and  $\beta \leq \eps p^{r-1}n$, contains a $K_r$-factor. 
\end{thm}
 We remark that the condition that $\delta(G)\geq cpn$ is natural. Indeed Definition~\ref{def:bijumbled} implies that almost all vertices will have degree at least $cpn$ and some lower bound on minimum degree is necessary to avoid isolated vertices.  Theorem~\ref{thm:maintriangle} follows directly from Theorem~\ref{thm:main} and much of the context and past results discussed above have analogous statements when $r\geq 4$ with many authors also working in the more general setting of $(p,\beta)$-bijumbled graphs. In particular, for all $r\geq 3$, a condition of $\beta=o(p^{r-1}n)$ guarantees a copy of $K_r$ and before Theorem~\ref{thm:main} the best condition known for ensuring a $K_r$-factor was $\beta=o(p^{r-1}n/\log n)$ due to Nenadov~\cite{Nen18}. Another result due to Han, Kohayakawa, Person and the author \cite{HKMP18} appeared at roughly the same time as that of Nenadov and gave a condition of $\beta=o(p^{r}n)$ for a $K_r$-factor, which for $r\geq 4$ gives a stronger result than the previously best known condition of Allen, B\"ottcher, H\`an, Kohayakawa and Person \cite{ABHKP17}. Although this condition is weaker than Nenadov's only when the bijumbled graph is very dense, it turns out that the proof methods of both results
 will be useful in proving Theorem~\ref{thm:main}.

There is one key difference in the picture for the case when $r=3$ and when $r\geq 4$:  the tightness of the condition $\beta=o(p^{r-1}n)$ for \emph{both} the clique and the clique factor when $r\geq 4$ is unknown. We defer a more in depth discussion of this to our concluding remarks (Section \ref{sec:conclude}) and conclude this introduction by again focusing on the most interesting case of triangle factors  where we know that Theorem~\ref{thm:main} and Theorem~\ref{thm:maintriangle} are tight  
 due to the construction of Alon (and its generalisation to the whole range of densities by Krivelevich, Sudakov and Szab\'o) discussed above. Indeed, one of the reasons that the Krivelevich-Sudakov-Szab\'o conjecture  (Theorem~\ref{thm:maintriangle})  has attracted so much attention is that it marks a distinct difference between the behaviour of random graphs and that of (optimally) pseudorandom graphs. In random graphs, we know that triangles appear at density roughly $p=n^{-1}$, whilst for triangle factors the threshold is considerably denser, namely $p=n^{-2/3}(\log n)^{1/3}$ \cite{JohanssonKahnVu_FactorsInRandomGraphs} (see also recent results~\cite{heckel2018random,kahn2019asymptotics,kahn2020hitting,riordan2018random}  that  imply that  this threshold is sharp). On the other hand, 
there exists triangle-free, 
optimally pseudorandom graphs with density  roughly $n^{-1/3}$, 
but Theorem~\ref{thm:main} asserts that any pseudorandom graph whose density is a constant factor  larger than this is guaranteed to have not only a triangle  but  a triangle factor.  Furthermore, it follows from Theorem~\ref{thm:main} and (the proof of) a result of Han, Kohayakawa, Person and the author~\cite{HKMP18b} that even more is true. 

\begin{cor} \label{cor:2-factor}
For every $c>0$ there exists an $\eps>0$ such that any $n$-vertex $(p,\beta)$-bijumbled graph with $\delta(G)\geq cpn$, $p>0$ and  $\beta \leq \eps p^{2}n$ is~\emph{$2$-universal}. That is, given any graph $F$ on at most $n$ vertices, with maximum degree $2$, $G$ contains a copy of $F$. In particular,  any $(n,d,\lambda)$-graph $G$ with 
$\lambda\leq \eps d^2/n$ is $2$-universal.\end{cor}

Our proof of Theorem~\ref{thm:main} incorporates  discrete algorithmic techniques, probabilistic methods,  fractional relaxations and linear programming duality, 
and the method of absorption. 
In the next section we discuss the proof in detail and reduce the problem to proving two intermediate propositions and a lemma. These will then be proven in what follows after developing the necessary theory. 

\vspace{2mm}

{\bf{Remark.}} An accompanying conference version \cite{confversion} of this work  deals solely with the setting of Theorem~\ref{thm:maintriangle}. More technical parts of the proof are omitted there and we hope that it serves as a gentle introduction to the present paper.

\vspace{2mm}

{\bf{Acknowledgements.}} Much of this work was done during a visit of the author to IMPA, Rio de Janeiro, Brazil. I am very grateful to both Pedro Ara\'ujo and  Robert Morris  with whom I discussed many aspects of this project,  for their hospitality, support and enthusiasm. I am also very grateful to my coauthors from previous projects Jie Han, Yoshiharu Kohayakawa and Yury Person, for introducing me to this area and many of the techniques and approaches used to tackle these sorts of problems. Finally I am grateful to the anonymous referee, to Tibor Szab\'o and again to Robert Morris for providing many helpful suggestions  aiding the presentation and readability  of the paper.

\section{Proof of main theorem} \label{sec:proofmain}
The proof of Theorem \ref{thm:main} rests on the shoulders of the previous results~\cite{ABHKP17,HKP18a,HKP18b,HKMP18,HKMP18b,KSS04,Nen18} working towards the conjecture of Krivelevich, Sudakov and Szab\'o. Indeed it is fair to say that the solution of the conjecture would not have been possible without the insights and ideas of the many authors who tackled this problem. In this section, we discuss these as well as our novel ideas and lay out the key concepts and scheme of the proof. In doing so, we will reduce the theorem to several intermediate results, whose proofs will be the subject of the rest of the paper.

Our proof, like some of its predecessors \cite{ABHKP17,HKMP18,Nen18},  works by the method of absorption. It turns out that finding many vertex disjoint copies of $K_r$ in a $(p,\beta)$-bijumbled graph $G$ as in Theorem~\ref{thm:main}, is easy. This follows from a simple consequence of Definition~\ref{def:bijumbled} 
which guarantees that any small linear sized set of vertices contains a copy of $K_r$, see e.g. Corollary \ref{cor:transversalcliques}  \eqref{cor:r-2small} 
for a precise statement. Therefore we can greedily choose copies of $K_r$ to be in our $K_r$-factor and continue this process until we are left with some small leftover set of vertices $L$, where small here means, of size at most $\eps r n$, say. However, at this point we get stuck; we have no way of guaranteeing the existence of a $K_r$ in $L$ and so we do not know how to get a larger set of vertex disjoint copies of $K_r$. The idea of absorption is to put aside an \emph{absorbing set} of vertices which can \emph{absorb} the leftover vertices $L$ into a $K_r$-factor. That is, before running this greedy process to build a $K_r$-factor, we find some special set of vertices $X\subset V(G)$ which has the property that for \emph{any} small set of vertices $L\subset V(G)\setminus X$, there is a $K_r$-factor in $G[X\cup L]$ (provided the trivial divisibility constraint that $r|(|X|+|L|)$). If we can find such an $X$ in $G$, then we can put it to one side and run the greedy argument to cover almost all the vertices which do not lie in $X$, with vertex disjoint copies of $K_r$. We can then use the absorbing property to \emph{absorb} the leftover vertices $L$ and get a full $K_r$-factor.

\begin{figure}[h]
    \centering
    \captionsetup{width=0.94\linewidth}
  \includegraphics[scale=0.84]{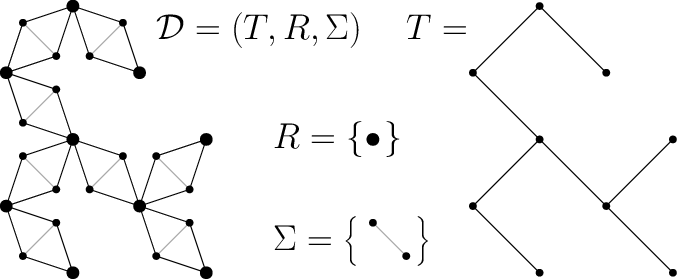}
    \caption{   \label{fig:DiamondTree}   An example of a $K_3$-diamond tree $\cD=(T,R,\Sigma)$  of order $9$ shown on the left. The removable vertices $R$ are the larger vertices of $\cD$ and the interior cliques $\Sigma$ are the edges given in grey. The auxiliary tree $T$ is depicted on the right.} 
  \end{figure}

This leaves the challenge of defining some structure in $G$ which has this absorbing property and finding such a structure (on some vertex set $X$) in $G$. The building blocks of our absorbing structure will be  subgraphs that we call \emph{$K_r$-diamond trees}. In words, a $K_r$-diamond tree $\cD=(T,R,\Sigma)$ is the graph obtained by taking a tree $T$ and replacing each edge $e\in E(T)$ by a copy of $K_{r+1}^-$ whose degree $r-1$ vertices are the vertices of $e$ and whose degree $r$ vertices are new and distinct from previous choices, see Figure~\ref{fig:DiamondTree} for an example. The following definition formalises this notion. 

\vspace{4mm}

\begin{dfn} \label{def:diamondtree}
A \emph{$K_r$-diamond tree} $\cD$ of \emph{order} $m$ in a graph $G$ is a tuple $\cD=(T,R,\Sigma)$ where $T$ is an (auxiliary) tree of order $m$ (i.e. with $m$ vertices), $R\subset V(G)$ is a subset of $m$ vertices of $G$ and\footnote{Here and throughout, we use the notation $K_{r-1}(G)$ to denote the family of $(r-1)$-cliques in $G$.} $\Sigma\subset K_{r-1}(G)$ is a set of $m-1$ copies of $K_{r-1}$ in $G$ such that the following holds.  There are bijective maps $\rho:V(T)\to R$ and $\sigma:E(T)\to \Sigma$ such that: 
\begin{itemize}
    \item The copies of $K_{r-1}$ in $\Sigma$ are pairwise vertex disjoint in $G$ and they are also disjoint from $R$ i.e. $V(S)\cap V(S')=\emptyset$  and $V(S)\cap R=\emptyset$ for all $S,S'\in \Sigma$;
    \item For all $e=uv\in E(T)$, we have that $V(\sigma(e))\subseteq N^G(\rho(u))\cap N^G(\rho(v))$. That is,  the $r-1$-clique $\sigma(e)\in K_{r-1}(G)$ can be extended to a copy of $K_r$ in $G$ by adding the vertex $\rho(u)$ and likewise with $\rho(v)$.  
\end{itemize}
We refer to $R$ as the set of \emph{removable vertices} of $\cD$ and to $\Sigma$ as the set of \emph{interior cliques} of $\cD$.  We define the vertices of $\cD$ to be all the removable vertices and the vertices in interior cliques. That is,~$V(\cD):=\left(\cup_{S\in \Sigma}V(S)\right) \cup R$. Finally we define the \emph{leaves} of the diamond tree  to be the vertices which are images of leaves in $T$ under $\rho$.
\end{dfn}

 Note that a $K_r$-diamond tree of order $m$ has exactly $(m-1)r+1$ vertices.  Krivelevich~\cite{kri97} used $K_3$-diamond trees in an absorption argument for triangle factors in random graphs which is often cited as one of the first appearances of the absorption method. Nenadov~\cite{Nen18} also used this idea in his result that got within a $\log$-factor of Theorem~\ref{thm:main}. The utility of these subgraphs in absorption arguments comes from the following key observation which shows that they can contribute to a $K_r$-factor in many ways.

\begin{figure}[h]
    \centering
    \captionsetup{width=.94\linewidth}
  \includegraphics[scale=0.84]{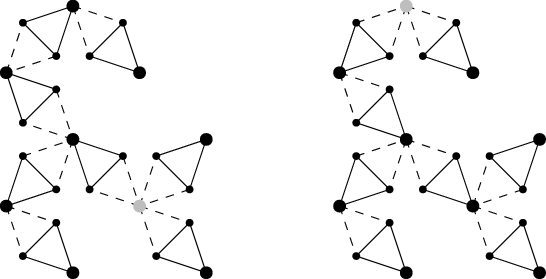}
    \caption{   \label{fig:RemovableVertices} Some examples of the $K_3$-factors found after removing a removable vertex from the $K_3$-diamond tree in Figure~\ref{fig:DiamondTree} (see Observation~\ref{obs:diamondtree}).} 
  \end{figure}

\begin{obs} \label{obs:diamondtree}
Given a $K_r$-diamond tree $\cD=(T,R,\Sigma)$ in $G$, we have that for \emph{any} removable vertex~$v\in R$, there is a $K_r$-factor of $G[V(\cD)\setminus \{v\}]$. Indeed, consider $u=\rho^{-1}(v)$ in $V(T)$ and the map $\phi:E(T)\to V(T)\setminus\{u\}$ which maps each edge $e$ of $T$ to the vertex in $e$ which has the larger distance from $u$ in $T$. Then $\phi$ is a bijection and taking the copies of $K_r$ on $\sigma(e)\cup \rho(\phi(e)) $ for each edge $e\in E(T)$ gives the required $K_r$-factor. See Figure~\ref{fig:RemovableVertices} for some examples. 
\end{obs}

Observation \ref{obs:diamondtree} works for any underlying auxiliary tree $T$. It turns out that in the $(p,\beta)$-bijumbled graphs $G$ we are interested in, one can find $K_r$-diamond trees of any order up to linear size. Indeed, one can use the argument of  Krivelevich \cite{kri97}  to construct  these or a different argument due to Nenadov \cite{Nen18}. The method of Nenadov gives diamond trees whose  auxiliary tree is a path whilst the argument of Krivelevich gives no control over the underlying auxiliary tree which defines the diamond tree found. As a key part of our argument, we will need to prove the existence of diamond trees which have extra structure, as we discuss shortly.

\vspace{2mm}

In order to utilise the absorbing power of diamond trees, we need to group them together in collections. The following definition of an \emph{orchard} captures how we do this.

\begin{dfn} \label{def:orchard}
We say a collection $ \bc{O}=\{\cD_1,\ldots,\cD_k\}$ of pairwise vertex disjoint $K_r$-diamond trees in a graph $G$ is a \emph{$(k,m)_r$-orchard} if there are $k$ diamond trees in the collection and each has order at least $m$ and at most $2m$.  We refer to $k$ as the \emph{size} of the orchard and $m$ as the \emph{order} of the orchard\footnote{Note that we abuse notation slightly here. Indeed we refer to \emph{the} order of an orchard although this may not be uniquely defined by the orchard. We take the convention that when we refer to the order of an orchard, we simply fix one of the possible orders arbitrarily, noting that these possible orders differ by  a factor of at most 2.}. We denote by $V(\bc{O})$, the vertices featuring in diamond trees in $\bc{O}$, that is $V(\bc{O})=\bigcup_{i\in[k]}V(\cD_i)$. Finally, if $\bc{O}'\subseteq \bc{O}$ is a subset of diamond trees in an orchard $\bc{O}$, we call $\bc{O}'$ a \emph{suborchard} of $\bc{O}$.
\end{dfn}

The term orchard here is supposed to be instructive, indicating that this is a `neat' collection of diamond trees that all have a similar order and are completely disjoint from one another. As noted in Observation \ref{obs:diamondtree}, a $K_r$-diamond tree can contribute to a $K_r$-factor in many ways. By grouping together many vertex disjoint $K_r$-diamond trees into a $(k,m)_r$-orchard such that $km =\Omega(n)$, we get a structure with a  strong absorbing property, as the following lemma shows. We say a~$(K,M)_r$-orchard $\bc{O}$ \emph{absorbs} a $(k,m)_r$-orchard $\bc{R}$ if there is a $((r-1) k,M)_r$-suborchard $\bc{O}'\subset \bc{O}$, such that there is a $K_r$-factor in $G[V(\bc{R})\cup V(\bc{O}')]$.

\begin{lem} \label{lem:absorptionbetweenlayers}
For any $3\leq r\in \NN$ and $0<\zeta,\eta<1$ there exists an $\eps>0$ such that the following holds for any $n$ vertex $(p,\beta)$-bijumbled graph $G$ with $\beta\leq \eps p^{r-1}n$. 
Let $\bc{O}$ be a $(K,M)_r$-orchard  in $G$ such that $KM\geq\zeta n$. Then there exists a set $B \subset V(G)$ such that $|B|\leq \eta p^{2r-4}n$ and  $\bc{O}$ absorbs any $(k,m)_r$-orchard  $\bc{R}$ in $G$ with
\begin{equation} \label{eq:absorptionconstraints}
  V(\bc{R})\cap (B \cup V(\bc{O}))=\emptyset, \qquad   k\leq K/(8r) \qquad \mbox{ and } \qquad
  kM\leq mK. \end{equation}
\end{lem}

Morally, Lemma~\ref{lem:absorptionbetweenlayers} says that \emph{large orchards absorb small orchards}. Here, by large we refer to both the size and the order of the orchards. Indeed the second condition in~\eqref{eq:absorptionconstraints} shows that the larger orchard has to have a larger size than the smaller orchard. This is the critical condition when we want absorption between orchards of similar order. The third condition shows that the ratio between the orders of the orchards is constrained by the ratio of the sizes. That is, the larger $\bc{O}$ is compared to $\bc{R}$ with respect to their sizes, the smaller $\bc{R}$ can be than $\bc{O}$ with respect to their orders. This will be the critical condition when we want absorption between orchards of (polynomially) different orders. The first condition in~\eqref{eq:absorptionconstraints} simply states that in order for $\bc{O}$ to absorb $\bc{R}$, we need that $\bc{R}$ avoids some small set of bad vertices $B$. This will be easy to implement in applications.

Lemma~\ref{lem:absorptionbetweenlayers} will be proven in Section~\ref{sec:absorbing orchards}. It provides us with an absorption property between two distinct orchards. We will also need an absorption property within orchards themselves, showing that we can find a large suborchard which hosts a $K_r$-factor in $G$. Given Observation \ref{obs:diamondtree}, in order to find  $K_r$-factors on suborchards it suffices to find copies of $K_r$ which traverse sets of removable vertices. We therefore make the following definition. 

\begin{dfn} \label{def:trianglehypergraph}
Given a $(k,m)_r$-orchard $\bc{O}=\{\cD_1,\ldots,\cD_k\}$ in a graph $G$, the \emph{$K_r$-hypergraph} generated by $\bc{O}$, denoted $H=H(\bc{O})$,  is the $r$-uniform hypergraph with vertex set $V(H)=\bc{O}$ and with $\{\cD_{i_1},\ldots,\cD_{i_r}\}$ for distinct $i_1,\ldots,i_r\in [k]$ forming a hyperedge in $H$ if and only if there is a copy of $K_r$ traversing\footnote{Here and throughout, when we say a copy of $K_r$ \emph{traverses} $r$ disjoint sets of vertices if it contains one vertex from each set.} the sets $R_{i_1},\ldots,R_{i_r}$ in $G$, where $R_{i_j}$ is the set of removable vertices of $\cD_{i_j}$ for all~$j$. 
\end{dfn}

Appealing to Observation \ref{obs:diamondtree} then gives the following,  as finding copy of $K_r$  traversing $r$ sets of removable vertices  removes exactly one vertex from each set.

\begin{obs} \label{obs:matchings in triangle hypergraph}
If $\bc{O}$ is an orchard of $K_r$-diamond trees in a graph $G$ and $H(\bc{O})$ contains a perfect matching, then~$G[V(\bc{O})]$ contains a $K_r$-factor.
\end{obs}

We will be particularly interested in orchards which contain near $K_r$-factors in a robust way. This gives us the notion of  a shrinkable orchard. 

\begin{dfn} \label{def:shrinkable}
Given $0<\gamma<1$, we say a $(k,m)_r$-orchard $\bc{O}$ in a graph $G$ is \emph{$\gamma$-shrinkable} if there exists a suborchard $\bc{Q}\subset \bc{O}$ of size at least $\gamma k$ such for \emph{any} suborchard $\bc{Q}'\subseteq \bc{Q}$, we have that there is a matching in $H:=H(\bc{O}\setminus \bc{Q}')$ covering all but $k^{1-\gamma}$ of the vertices of $H$.
\end{dfn}

Our first key proposition gives the existence of shrinkable orchards. It will be discussed in Section~\ref{sec:shrinkable orchards} and proven in  Sections \ref{sec:small order shrinkable} and \ref{sec:large order shrinkable}.

\begin{prop} \label{prop:shrinkable orchards}
For any $3\leq r\in \NN$ and $0<\alpha,\gamma<1/2^{12r}$ there exists an $\eps>0$ such that the following holds for any $n$ vertex $(p,\beta)$-bijumbled graph $G$ with $\beta\leq \eps p^{r-1}n$ 
and any vertex subset $U\subseteq V(G)$ with $|U|\geq n/2$. For any $m\in \mathbb{N}$ with $1\leq m \leq n^{7/8}$  there exists a $\gamma$-shrinkable $(k,m)_r$-orchard $\bc{O}$ in $G[U]$ with $k\in \NN$ such that $\alpha n\leq km \leq  2\alpha n$. 
\end{prop}

Given Lemma~\ref{lem:absorptionbetweenlayers} and Proposition~\ref{prop:shrinkable orchards}, we know that we can find orchards which contain large $K_r$-factors and that large orchards can absorb smaller orchards. This suggests the following approach for giving an absorbing structure which can absorb leftover vertices in our $(p,\beta)$-bijumbled graph $G$ (we keep the discussion at a high level  here to highlight the key idea, the details of this scheme are elaborated in discussions within the proof of Theorem~\ref{thm:main}). Find a sequence of vertex disjoint shrinkable orchards, each on a linear number of vertices. Each orchard in the sequence will have a larger order than that of the previous orchard and the first orchard in the sequence will be composed of linearly many $K_r$-diamond trees of constant size. 
We can then run a \emph{cascading absorption} through the sequence of orchards. That is, given some small leftover set of vertices $L$ (which is itself a $(|L|,1)$-orchard), we use the first orchard in the sequence to absorb $L$. We then use that the first orchard is shrinkable and so we can cover most of what remains of the first orchard with vertex disjoint copies of $K_r$. There will be some $K_r$-diamond trees of the first orchard left at the end of this and for these we appeal to Lemma~\ref{lem:absorptionbetweenlayers} to absorb this small suborchard using the second orchard. Then again, the second orchard is shrinkable and so  the remainder of the second orchard can be almost fully covered with vertex disjoint $K_r$s, leaving some small leftover suborchard uncovered. We then  repeat to absorb this leftover with the third orchard and continue in this fashion. In this way we cascade the absorption through the orchards and each time we do this, we increase the order of the orchard which we need to absorb. 

This approach is promising but we need to cut this process off at some point and find a \emph{full} $K_r$-factor on the vertices which have not already been covered by vertex disjoint copies of $K_r$. The next proposition states that once the orchard has a large enough order, we can find a structure that can \emph{fully absorb} any leftover.

\begin{prop} \label{prop:finalabsorption}
For any $3\leq r\in \NN$ and $0<\alpha,\eta<1/2^{3r}$ there exists an $\eps>0$ such that the following holds for any $n$ vertex $(p,\beta)$-bijumbled graph $G$ with $\beta\leq \eps p^{r-1}n$ 
and any vertex subset $W\subseteq V(G)$ with $|W|\geq n/2$.

 There exist vertex subsets $A,B\subset V$ such that $A\subset W$,  $|A|\leq \alpha n$, $|B|\leq \eta p^{2r-4} n$ and for  any $(k,m)_r$-orchard $ \bc{R}$  
whose vertices lie in $V(G)\setminus (A\cup B)$, we have that if $|A|+|V(\bc{R})|\in r \mathbb{N}$, $k\leq \alpha^2 n^{1/8}$ and $m\geq n^{7/8}$ then $G[A\cup V(\bc{R})]$ has a $K_r$-factor.  
\end{prop}

In our absorption scheme sketched above,  $\bc{R}$ will be the leftover of the last orchard (the one with the largest order) after having cascaded the absorption through the sequence of orchards. Proposition~\ref{prop:finalabsorption} states that the vertex set $A$ can fully absorb this $\bc{R}$. Hence when constructing our absorbing structure, we will first find $A$ and then construct our sequence of shrinkable orchards so that they avoid $A$ and also the small set of bad vertices $B$ given by Proposition~\ref{prop:finalabsorption}. Finally, we remark that the necessity  of $W$ in Proposition~\ref{prop:finalabsorption} comes from the fact that before we find our absorbing structure, we will put aside some small set $Y$ which will be used later in the proof to help with bad vertices, and so need to find $A$ in $W=V(G)\setminus Y$. 

We are now in a position to prove Theorem \ref{thm:main}, using only Lemma \ref{lem:absorptionbetweenlayers}, Propositions \ref{prop:shrinkable orchards} and \ref{prop:finalabsorption}, some simple properties of $(p,\beta)$-bijumbled graphs and Chernoff's Theorem (Theorem \ref{thm:chernoff}), a well-known result which gives concentration of binomial random variables. 

\begin{proof}[Proof of Theorem \ref{thm:main}]
For convenient reference throughout the proof, let us fix our constants \begin{equation} \label{eq:final proof constants}\gamma:=\frac{c}{2^{24r}}, \quad \lambda:=\gamma^2, \quad   \alpha:=\lambda^2, \quad \zeta=\alpha^2 \quad  \eta:=\zeta^2 \quad  \mbox{and} \quad  t:=\frac{7}{8\lambda}.\end{equation} We further fix $\eps>0$ much smaller than all these constants and sufficiently small enough to apply Lemma \ref{lem:absorptionbetweenlayers} and Propositions \ref{prop:shrinkable orchards} and \ref{prop:finalabsorption} with these parameters.  
We also use some simple consequences of Definition~\ref{def:bijumbled} which imply that, by choosing $\eps>0$ sufficiently small, we guarantee that any vertex subset of size $\zeta pn$ contains a copy of $K_{r-1}$ whilst any vertex set of size $\zeta n$ contains a copy of $K_r$, see e.g. Corollary~\ref{cor:transversalcliques}. Finally we note that if $\delta(G)\geq (1-1/r)n$ then it follows from Theorem \ref{thm:HajnalSzem} that $G$ has a $K_r$-factor  
and so we can assume that $\delta(G)< (1-1/r)n$. For such $n$ vertex $(p,\beta)$-bijumbled graphs  $G$ with $\beta \leq \eps p^{r-1}n$, a well-known fact (see Fact \ref{fact:largen}) implies that by choosing $\eps>0$ to be sufficiently small, we can assume that $n$ is sufficiently large in what follows as otherwise no $n$ vertex $(p,\beta)$-bijumbled graphs with $\beta\leq \eps  p^{r-1}n$ exist and the theorem is vacuously true.  Moreover, another well-known fact (see Fact~\ref{fact:dense}) implies that $(p,\beta)$-bijumbled graphs cannot be too sparse. In particular, with our condition on $\beta$,  by choosing $\eps>0$ to be sufficiently small, we can also assume that $p\geq n^{-1/3}$ in what follows.

\vspace{2mm}

Before finding our $K_r$-factor in $G$ we need to do some preparation. We begin by setting aside a randomly chosen subset of vertices $Y\subset V(G)$. We let each vertex be in $Y$ with probability $\alpha$. It follows from Chernoff's Theorem (see Theorem \ref{thm:chernoff}) and a union bound that with high probability as $n$ tends to infinity, we have that $|Y|\leq 2 \alpha n$ and $\deg_Y(v)\geq c\alpha pn/2$ for all $v\in V(G)$. Indeed, this follows because $\EE[\deg_Y(v)]\geq c\alpha pn=\Omega(n^{2/3})$ for each $v\in V(G)$. Therefore, as $n$ is large, we can fix such an instance of $Y$. We will use the vertices of $Y$ to find copies of $K_r$ containing `bad' vertices later in the argument. 

Next, we apply Proposition \ref{prop:finalabsorption} (with $W=V(G)\setminus Y$) to obtain  vertex sets $A\subset V(G)\setminus Y$  and $B$ such that  $|A|\leq \alpha n$, $|B|\leq \eta p^{2r-4} n$ and we have the following key absorption property. For  any $(k,m)_r$-orchard $ \bc{R}$ 
whose vertices lie in $V(G)\setminus (A\cup B)$, we have that if $|A|+|V(\bc{R})|\in r \mathbb{N}$, $k\leq \zeta n^{1/8}$ and $m\geq n^{7/8}$ then $G[A\cup V(\bc{R})]$ has a $K_r$-factor. That is, $A$ can absorb orchards whose order is sufficiently large.

\vspace{2mm}

As sketched above, the idea is now to provide constantly many (namely, $t+1$) vertex-disjoint shrinkable orchards $\bc{O}_0,\bc{O}_1,\ldots \bc{O}_t$,  each on a linear number of vertices and whose vertices are disjoint from $A$. The order of these orchards will increase slightly (namely, by a factor of $n^\lambda$) at each step in the sequence. Due to our definition of $t$ \eqref{eq:final proof constants}, we can have that $\bc{O}_0$ has $\Omega(n)$ diamond trees of constant order while $\bc{O}_t$ has diamond trees of order $\Omega(n^{7/8})$. The point is that we will be able to repeatedly apply Lemma~\ref{lem:absorptionbetweenlayers} and the fact that each orchard is shrinkable to create a cascading absorption through the shrinkable orchards. Indeed $\bc{O}_0$ will be able to absorb leftover vertices and each $\bc{O}_i$ will be able to absorb any leftover $K_r$-diamond trees in $\bc{O}_{i-1}$, after using that $\bc{O}_{i-1}$ is shrinkable to cover almost all of the vertices of $\bc{O}_{i-1}$ with disjoint copies of $K_r$. Once this absorption reaches $\bc{O}_t$, we will be able to use $A$ to absorb the leftover $K_r$-diamond trees in $\bc{O}_t$ and complete a $K_r$-factor. In fact, when absorbing between orchards we do not use \emph{all} of $\bc{O}_i$ to absorb leftover diamond trees in $\bc{O}_{i-1}$ but rather a suborchard $\bc{Q}_i\subset \bc{O}_i$ which contains a $\gamma$ proportion of the $K_r$-diamond trees in $\bc{O}_i$. Indeed this $\bc{Q}_i$ is provided by the fact that $\bc{O}_i$ is shrinkable (see Definition~\ref{def:shrinkable}) and guarantees that removing diamond trees from $\bc{Q}_i$ will not prevent us from covering almost all of what remains of $\bc{O}_i$ with vertex disjoint copies of $K_r$.

\vspace{2mm}

In detail, we collect what we require in the following claim. 

\begin{clm} \label{clm:final proof}
There exists vertex disjoint orchards $\bc{O}_0,\bc{O}_1,\ldots,\bc{O}_t$ in $G$ such  that following properties hold. 
\begin{enumerate}[label=(\roman*)]
\item  \label{clm:disjointness} For all $0\leq i \leq t$, we have that $V(\bc{O}_i)\cap (A \cup B\cup Y)=\emptyset$.
\item  \label{clm:sizes} For each $0 \leq i \leq t$, fixing $m_i:=n^{i\lambda}$, we have that $\bc{O}_i$ is a $(k_i,m_i)_r$-orchard for some~$k_i$ such that $\alpha n \leq k_im_i\leq 2 \alpha n$. 
\item \label{clm:shrinkable} Each $\bc{O}_i$ is $\gamma$-shrinkable with respect to some suborchard $\bc{Q}_i\subset \bc{O}_i$ such that \[k_i^*:=|\bc{Q}_i|\geq\gamma k_i.\]
\item \label{clm:absorbing} For $1\leq i \leq t$, given any suborchard $\bc{P}\subset \bc{O}_{i-1}$ such that $|\bc{P}|\leq k_{i-1}^{1-\gamma}$, we have that $\bc{Q}_i$ absorbs $\bc{P}$. 
\end{enumerate}
\end{clm}

Before verifying the claim, let us see how we can derive the theorem using the claim. So suppose we have found such orchards $\bc{O}_0,\ldots, \bc{O}_t$ and fix
\[X:=A \bigcup \left(\bigcup_{i=0}^tV(\bc{O}_i)\right).\] 
Furthermore, note that as $k^*_0m_0= k_0^*\geq \gamma\alpha n\geq \zeta n$, by Lemma~\ref{lem:absorptionbetweenlayers}, there exists some set $B_0 \subset V(G)$ such that  $|B_0|\leq \eta p^{2r-4}n$  and $\bc{Q}_0$ absorbs any $(k,1)$-orchard\footnote{Note that a $(k,m)_r$-orchard with $m=1$ is simply a set of vertices. Each $K_r$-diamond tree in the orchard has order $1$ and so is a single isolated vertex.} $\bc{R}$ such that 
\begin{equation} \label{eq:lowest level}k\leq \zeta n\leq \frac{\gamma \alpha n}{8r}\leq \frac{k_0^*}{8r}
\end{equation} 
and $V(\bc{R}) \cap (B_0\cup V(\bc{Q}_0))=\emptyset$.  Indeed the condition on $k$ comes from \eqref{eq:absorptionconstraints}, using that $m_0=1$ and our lower bound on $k_0^*$. Fix $Z:=B_0 \setminus X$ and note that $z:=|Z|\leq \eta p^{2r-4}n$ as $Z$ is a subset of $B_0$. Note also that $X\cap Y=\emptyset$ due to part \ref{clm:disjointness} of Claim~\ref{clm:final proof} and how we defined $A$.

\begin{figure}[h]
    \centering
  \includegraphics[scale=0.84]{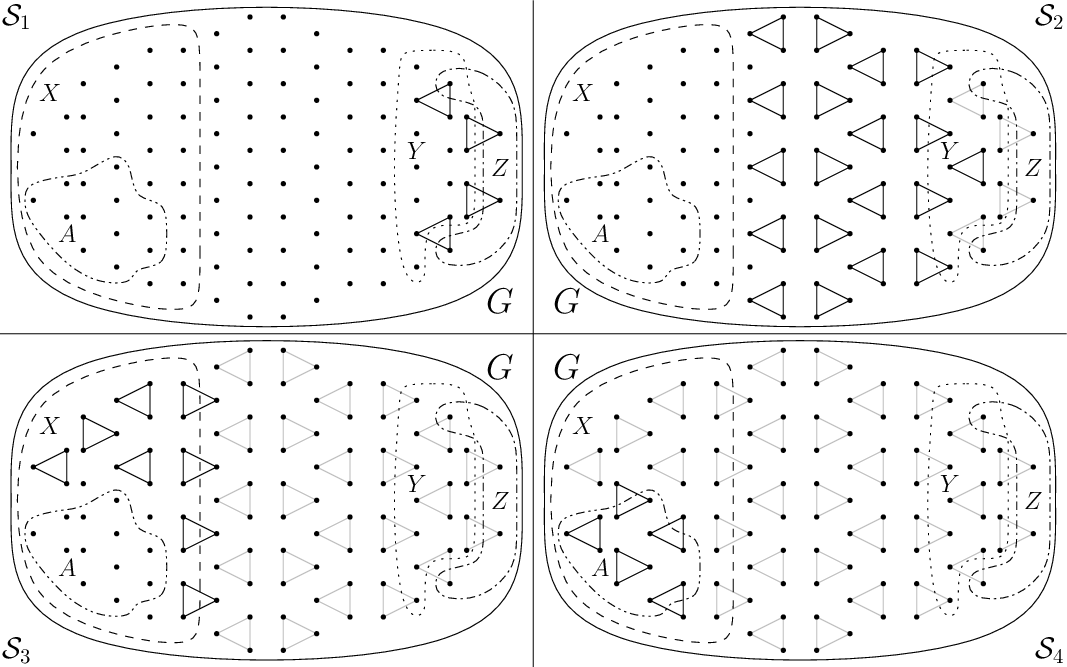}
    \caption{   \label{fig:phases} A schematic to demonstrate the triangles found (and the vertex sets they cover) by our four-phase   algorithm that finds a $K_3$-factor in $G$. } 
  \end{figure}

\vspace{2mm}

We are now ready to find our $K_r$-factor $\cS$ which we  do algorithmically in four phases. See Figure~\ref{fig:phases} for a visual guide to the cliques found in each phase.   So let us initiate with $\cS_1=\emptyset$. In the first phase we find copies of $K_r$ containing the vertices in $Z$, using some vertices in $Y$. So let us order the  vertices of $Z$ arbitrarily as $Z:=\{b_1,\ldots,b_z\}$ and fix $Y_1:=Y \setminus Z$.  Now for $1\leq j \leq z$, we find an $r$-clique $S_j$ containing~$b_j$ and $r-1$ vertices of $ Y_j$. We add $S_j$ to $\cS_1$, fix $Y_{j+1}:=Y_j \setminus V(S_j)$ and move to step $j+1$. To see that we can always find such a clique, note that  for each $j\in[z]$ we have that 
\[\deg_{Y_j}(b_j)\geq \deg_{Y}(b_j)-|Z|-r(j-1) \geq \frac{c\alpha pn}{2}- r \eta p^{2r-4}n \geq \zeta pn, \]
recalling our key property of $Y$ and using the definitions of our constants \eqref{eq:final proof constants}. A simple consequence of \eqref{eq:bijumbled} (see e.g. Corollary \ref{cor:transversalcliques} \eqref{cor:forbiddingsubgraph} \ref{cor:r-1cliqueinpn})  
implies that there is a copy of $K_{r-1}$ in $N(b_j)\cap Y_j$ and so this forms an $r$-clique $S_j$ with~$b_j$. In this way, we see that we succeed at every step $j$ and at the end of the first phase we have a set of vertex disjoint $r$-cliques  $\cS_1$ in $G$ of size $z$, such that every vertex in $Z$ is contained in a clique in $\cS_1$.

In the second phase we find the majority of the $K_r$-factor which we do greedily. We initiate with $\cS_2=\emptyset$ and $W= V(G) \setminus (X \cup V(\cS_1))$. Now whilst $|W|\geq \zeta n$, we can find an $r$-clique  $S$ in $W$. Again, this is a simple consequence of \eqref{eq:bijumbled}, see  Corollary \ref{cor:transversalcliques} \eqref{cor:r-2small}.   
We add $S$ to $\cS_2$ and delete its vertices from $W$. Therefore at the end of the second phase, we are left with some vertex set $L\subset V(G) \setminus X$ such that $|L|\leq \zeta n$ and $\cS_1 \cup \cS_2$ form a $K_r$-factor in $G[V(G)\setminus (X\cup L)]$. 

\vspace{2mm}

\begin{figure}[h]
    \centering
  \includegraphics[scale=0.84]{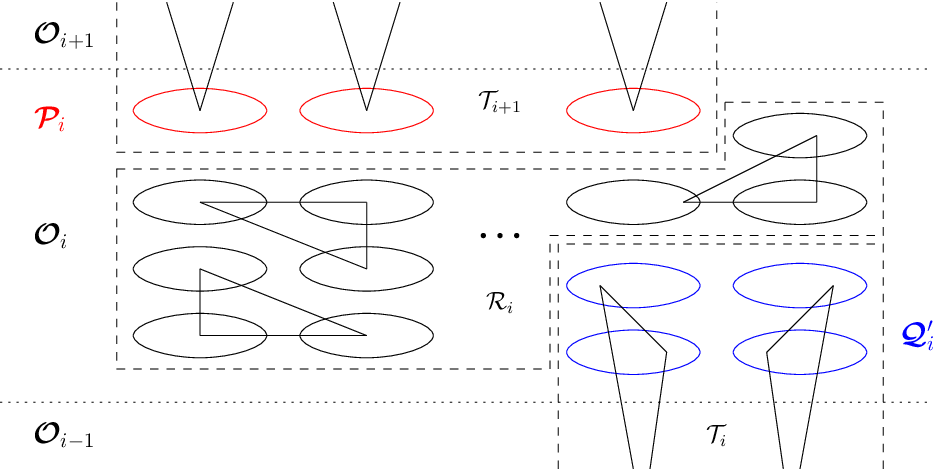}
    \caption{   \label{fig:phase3} A closer look at phase 3 of the algorithm in the case $r=3$. } 
  \end{figure}

In our third phase, we will find vertex disjoint $r$-cliques $\cS_3$  which cover $L$ and use almost all the vertices of $X \setminus A$. 
We begin by fixing $\ell:=|L|$ and noting that $L$ is a $(\ell,1)$-orchard which we relabel as $\bc{P}_{-1}$. Now we run the following procedure for $0 \leq i \leq t$  (see Figure~\ref{fig:phase3}). We first absorb $\bc{P}_{i-1}$ using $\bc{Q}_i$. That is, we find a suborchard $\bc{Q}'_i\subset \bc{Q}_i$ such that there is a  $K_r$-factor $\cT_i$ in $G[V(\bc{P}_{i-1}) \cup V(\bc{Q}'_i)]$.  We add the $r$-cliques in $\cT_i$ to $\cS_3$. Then, using that $\bc{O}_i$ is $\gamma$-shrinkable (Claim~\ref{clm:final proof} \ref{clm:shrinkable}), we can define some $\bc{P}_i \subset \bc{O}_i \setminus \bc{Q}'_i$ such that $|\bc{P}_i|\leq k_i^{1-\gamma}$ and there is a perfect matching in the $K_r$-hypergraph $H(\bc{O}_i \setminus (\bc{Q}'_i \cup \bc{P}_i))$. By Observation \ref{obs:matchings in triangle hypergraph}, this perfect matching gives a $K_r$-factor $\cR_i$ in $G[V(\bc{O}_i \setminus (\bc{Q}'_i \cup \bc{P}_i))]$. We add $\cR_i$ to $\cS_3$ and move to step $i+1$ or finish if $i=t$.  Note that in order to find $\cT_i$ and $\bc{Q}'_i$ in each step $i\geq 1$, we appeal to Claim \ref{clm:final proof} \ref{clm:absorbing}, whilst when $i=0$, the existence of $\cT_1$ and $\bc{Q}'_1$ is guaranteed by the fact that $L$ is a $(\ell,1)$-orchard with $\ell\leq \zeta n$ as in \eqref{eq:lowest level} and $L$ is disjoint from $Z$ and hence $B_0$.

\vspace{2mm}

Let $\bc{R}:= \bc{P}_t\subset \bc{O}_t$. We have that $\cS_1 \cup \cS_2 \cup \cS_3$ is a $K_r$-factor in $G[V(G) \setminus (A \cup V(\bc{R}))]$. Hence as $r|n$, we must have that $r|(|A|+|V(\bc{R})|)$. Moreover, $\bc{R}$ is a $(k,m)_r$-orchard with $k\leq k_t^{1-\gamma}\leq(2\alpha n^{1/8})^{1-\gamma}<\alpha^2 n^{1/8}$  and $m=n^{7/8}$. Finally, note that $V(\bc{R}) \cap B=\emptyset$ due to property  \ref{clm:disjointness} of Claim \ref{clm:final proof}. Therefore, by the key property of the absorbing vertex set $A$ in Proposition \ref{prop:finalabsorption}, we have that there is a $K_r$-factor $\cS_4$ in $G[A\cup V(\bc{R})]$. It follows that $\cS:=\cS_1 \cup \cS_2 \cup \cS_3 \cup \cS_4$ is a  $K_r$-factor in $G$, completing the proof.

\vspace{2mm} 

It remains to establish Claim~\ref{clm:final proof} and find the shrinkable orchards as stated. 
We will do this  algorithmically in decreasing order. The reason for this is that in order for \ref{clm:absorbing} to hold we will appeal to Lemma~\ref{lem:absorptionbetweenlayers} and therefore there will be some set of bad vertices $B_i$ which we want $\bc{O}_{i-1}$ to avoid. In fact, we will ensure that $\bc{O}_{i-1}$ avoids $B_j$ for all $i\le j \le t$. This is not necessary but eases our definitions (as we do not have to reintroduce vertices into the pool $U_i$ of available vertices); the important condition in what follows is that $\bc{O}_{i-1}$ avoids $B_i$ for all $i$.

We start by fixing   $U_{t+1}:= V(G)\setminus (A\cup B\cup Y)$. Now for $t\geq i \geq 0$ in descending order, we apply Proposition \ref{prop:shrinkable orchards} to find a $\gamma$-shrinkable $(k_i,m_i)_r$-orchard $\bc{O}_i$  such that $\alpha n\leq k_im_i\leq 2 \alpha n$ and $V(\bc{O}_i)\subset U_{i+1}$.   
We then define $U_i$ as follows.  As $\bc{O}_i$ is $\gamma$-shrinkable, it defines some suborchard $\bc{Q}_i\subset \bc{O}_i$ as in condition \ref{clm:shrinkable} of the claim. 
Now as  $k^*_im_i\geq \gamma\alpha n\geq  \zeta n$, it   follows from Lemma \ref{lem:absorptionbetweenlayers} that there exists some $B_i \subset V(G)$ with $|B_i|\leq \eta p^{r-1} n$ such that if $k$ and $m$ satisfy $k\leq k^*_i/8r$ and $km_i\leq  mk^*_i$ and $\bc{R}$ is a $(k,m)_r$-orchard with $V(\bc{R})\subset V(G)\setminus (B_i\cup V(\bc{Q}_i))$ then $\bc{Q}_i$ absorbs $\bc{R}$.   We fix $U_i:=U_{i+1}\setminus (V(\bc{O}_i) \cup B_i)$ and move onto the next index $i-1$.

Let us first check that the process succeeds in finding the shrinkable orchards $\bc{O}_t,\ldots, \bc{O}_0$ at each step. Note that we start with $|U_{t+1}|\geq n-3\alpha n- \eta p^{2r-4} n \geq n-4 \alpha n$. Moreover at each step $i$, we remove  at most $\eta p^{2r-4}n \leq \alpha n$ vertices which lie in $B_i$ and at most $4r \alpha n$ vertices from $U_{i+1}$ which lie in the orchard~$\bc{O}_i$. Indeed the orchard is composed of $k_i$ vertex disjoint $K_r$-diamond trees of order at most $2m_i$,  the number of vertices in each diamond tree is less than $r$ times its order and $k_im_i\leq 2 \alpha n$.  Hence for all $t\geq i\geq 0$, we have that 
\[|U_i|\geq n-(t+2)\cdot 5r\alpha n  \geq  n/2,\] 
 using that $t\leq 1/\lambda$, $\alpha =\lambda^2$ and the definition of $\lambda$  (see \eqref{eq:final proof constants})  here. Hence Proposition~\ref{prop:shrinkable orchards} gives the existence of $\bc{O}_i$ at each step and verifies part \ref{clm:shrinkable} of the claim. Note that the conditions  \ref{clm:disjointness} and \ref{clm:sizes} also hold simply from how we defined the $\bc{O}_i$ and the fact that we found them in the sets $U_i$, each of which is a subset of $U_{t+1}$.

Thus it remains to verify the absorption property between orchards, namely \ref{clm:absorbing}. For each $1 \leq i \leq t$, we chose $\bc{O}_{i-1}$ to have vertices in $U_i$ and hence $V(\bc{O}_{i-1}) \cap (B_i \cup V(\bc{Q}_i))=\emptyset$.  Therefore  we have by Lemma~\ref{lem:absorptionbetweenlayers} that $\bc{Q}_i$ absorbs any suborchard $\bc{P}\subset \bc{O}_{i-1}$,  with $|\bc{P}|\leq k_{i-1}^{1-\gamma}$ if $k_{i-1}^{1-\gamma}\leq k_i^*/8r$ and $k_{i-1}^{1-\gamma}m_i\leq k_i^*m_{i-1}$. 

Now as $m_{i}=n^\lambda m_{i-1}$ and $n^{-\lambda} \leq 1/8r$ for sufficiently large $n$, it suffices to show that $k_{i-1}^{1-\gamma}\leq k_i^*n^{-\lambda}$. To see this,  note that due to the fact that  
$ \alpha n \leq k_{i-1}m_{i-1}, k_im_i\leq 2\alpha n$, 
we have  \[k_{i-1} \leq \frac{2 \alpha n}{m_{i-1}}=\frac{2\alpha n^{1+\lambda}}{m_i}\leq 2k_in^{\lambda} \leq \frac{2k^*_{i}n^{\lambda}}{\gamma}\] and using this as a lower bound for $k_i^*$, it suffices to show that 
\[k_{i-1}^\gamma\geq \frac{2n^{2\lambda}}{\gamma}.\]
This is certainly true as $k_{i-1}\geq k_t \geq \alpha n^{1/8}>n^{4\lambda/\gamma}$, 
 recalling that $4\lambda/\gamma=4\gamma$ from \eqref{eq:final proof constants}. This shows that \ref{clm:absorbing} holds for all $i$ and concludes the proof of the claim and hence the proof. 
\end{proof}

We remark that this proof scheme builds on that of Nenadov~\cite{Nen18} (which  in turn is influenced by that of Krivelevich~\cite{kri97}) who proved that $\beta \leq \eps p^2n/\log n$ suffices for a triangle factor in an $n$ vertex $(p,\beta)$-bijumbled graph. Indeed, Nenadov also uses a result akin to Lemma~\ref{lem:absorptionbetweenlayers} albeit between orchards whose orders only differ by a constant factor. His absorbing structure then contains a sequence of  $\Theta(\log n)$ orchards whose order increases by a constant factor along the sequence. Therefore the last orchard in the sequence contains constantly many diamond trees of large order (of order $\Theta(n/\log n)$). These can be fully absorbed because any three large sets host a transversal triangle and so transversal triangles between removable sets can be greedily found, completing a triangle factor in the last step. Similarly, the $(k,m)_3$-orchards used in his argument are not imposed to be shrinkable but can be seen to host a triangle factor on all but $o(k)$ of the diamond trees by again applying a greedy approach of finding transversal triangles.  The necessity of the $\log n$ in the condition of Nenadov is thus due to needing $\Theta(\log n)$ orchards in the absorbing structure and thus requiring slightly stronger properties of the $(p,\beta)$-bijumbled graph, for example the existence of triangles on sets of $\Omega(n/\log n)$ vertices.

\vspace{2mm} 

The key challenge in this paper is then to prove Propositions~\ref{prop:shrinkable orchards} and \ref{prop:finalabsorption}. Both results rely heavily on  a technique we develop to provide the existence of $K_r$-diamond trees in which we have some control over the set of removable vertices. This control is rather weak; we cannot guarantee  that any fixed vertices appear as removable vertices but we can give some flexibility over the choice of removable vertices. See Proposition~\ref{prop: choosing removable vertices} for the technical statement of what we prove. 

In order to prove Proposition~\ref{prop:shrinkable orchards}, we build on the approach of Han, Kohayakawa and Person~\cite{HKP18a,HKP18b}. Indeed, their result showing the existence of a near   $K_r$-factor (covering all but some~$n^{1-\eps'}$ vertices) in $(n,d,\lambda)$-graphs can be seen as a step towards proving the existence of shrinkable orchards of order~$1$. The approach involves showing the existence of a near-perfect matching in a subhypergraph $H'$ of the $K_r$-hypergraph generated by $V(G)$. In order to do this, one needs to carefully choose $H'$ and this is done by finding many fractional  $K_r$-factors in $G$ which do not put too much weight on (copies of $K_r$ containing) any given edge. Therefore, the methods of Krivelevich, Sudakov and Szab\'o~\cite{KSS04}, who proved the existence of singular fractional  $K_r$-factors, become pertinent. They use the power of linear programming duality to prove that certain expansion properties guarantee the existence of fractional factors. In our setting, it turns out that we need several distinct arguments to prove the existence of shrinkable orchards of different orders. We follow the scheme of using fractional factors (in fact, fractional perfect matchings in $K_r$-hypergraphs) but need to adapt the method for different applications and we rely crucially on probabilistic methods to actually prove the existence of orchards which satisfy the necessary expansion properties. 

\vspace{2mm}

It can be seen that Proposition~\ref{prop:shrinkable orchards} alone (for all orders of orchards) would lead via the same proof scheme to a condition of $\beta\leq \eps p^{r-1}n/(\log \log n)$. In order to close the gap and achieve Theorem~\ref{thm:main}, Proposition \ref{prop:finalabsorption} is necessary. To prove this, we appeal to a different absorption argument whose roots go back to an ingenious argument of Montgomery~\cite{M14a,M19} in his work on spanning trees in random graphs.  The approach, sometimes called the absorption-reservoir method, uses a bipartite graph, which we call a template (see Section~\ref{sec:templates}) as an auxiliary graph to define an absorbing structure. This idea was previously used by Han, Kohayakawa, Person and the author~\cite{HKMP18} to find clique factors in pseudorandom graphs and we used this approach again in our result on $2$-universality~\cite{HKMP18b}. Here we combine this idea with the absorbing power of orchards and prove Proposition~\ref{prop:finalabsorption} with a three-stage algorithm which finds the absorbing structure necessary. 

\vspace{2mm}

The rest of this paper is organised as follows. In the next section, we run through the necessary preliminaries, providing the background theory that we will use. This includes properties of bijumbled graphs, the study of perfect fractional matchings via linear programs, probabilistic methods and the absorption-reservoir method of Montgomery~\cite{M14a,M19}. In Section~\ref{sec:diamondtree} we then study what kinds of diamond trees we can guarantee in our bijumbled graph. The key result here is Proposition~\ref{prop: choosing removable vertices}, which will be crucial at various points in our proof. We then turn to addressing the necessary results for the cascading absorption through the orchards in Section~\ref{sec:orchards}. We prove Lemma~\ref{lem:absorptionbetweenlayers} in Section~\ref{sec:absorbing orchards} and discuss  Proposition~\ref{prop:shrinkable orchards} in Section~\ref{sec:shrinkable orchards}, reducing it to two intermediate propositions which tackle small and large order shrinkable orchards separately. We go on to prove the existence of shrinkable orchards of small order in Section~\ref{sec:small order shrinkable}  and large order in Section~\ref{sec:large order shrinkable}. Finally we prove  Proposition~\ref{prop:finalabsorption} which provides the final absorption in the proof of Theorem~\ref{thm:main}, in Section~\ref{sec:absorbingstructure}.

\section{Preliminaries} \label{sec:Preliminaries}
\subsection{Notation} \label{sec:Notation}
For a graph $G$ and $r\in \NN$, we define $K_r(G)$ to be the set of copies of $K_r$ in $G$. When referring to (a copy of) a clique $S\in K_r(G)$, we will identify the copy with the set of vertices that hosts it. That is, we think of $S\in K_r(G)$  as a set of $r$ vertices which host a clique in $G$ rather than the copy of the clique itself. Given a set of $r$-cliques $\Sigma\subseteq K_r(G)$, we use the notation $V(\Sigma)$ to denote all vertices that feature in cliques in $\Sigma$, i.e. $V(\Sigma):=\cup_{S\in\Sigma}S$.   
We call a set of cliques $\Sigma\subset K_r(G)$ a \emph{matching} of cliques  if it is composed of pairwise vertex disjoint cliques, that is, $S\cap S'=\emptyset$ for any $S\neq S'\in \Sigma$. 
Now given subsets  $S,W\subset V(G)$ of vertices, we let  $N^G_W(S)$ denote the common neighbours of the vertices in $S$ which lie in $W$. That is, $N^G_W(S):=\left(\cap_{v\in S}N^G(v)\right)\cap W$. Likewise, we define $\deg^G_W(S):=|N^G_W(S)|$ to be the cardinality of this neighbourhood.   
If the graph $G$ is clear from context then we drop the superscripts. Also if $S=\{u\}$ is a single vertex, we will drop the set brackets. 
We say that a clique $S\in K_r(G)$ \emph{traverses}  vertex subsets $U_1,\ldots,U_r\subseteq V(G)$ if there exists some ordering of $S$ as  $S=\{u_1,\ldots,u_r\}$ such that $u_i\in U_i$  for all $i\in[r]$. Note that when the $U_i$ are pairwise disjoint this simplifies to requiring that $S$ contains one vertex from each $U_i$. However at times we will deal with not necessarily disjoint sets $U_i$ and so this more delicate definition is needed.  

\vspace{2mm}

If $H$ is an $r$-uniform hypergraph for some $r\in \NN$ and $v,u\in V(H)$, 
$\deg^H(v)$ denotes the number of edges in $H$ containing $v$, and $\codeg^H(u,v)$ denotes the number of edges of $H$ which contain both $u$ and $v$.  If the hypergraph $H$ is clear from context, we drop the superscripts.  If $H$ is an $r$-uniform hypergraph with $r\geq 3$ and $J$ is a   \emph{$2$-uniform} graph  on the same vertex set $V(H)$, then $H_J$ denotes the subhypergraph of $H$ given by all edges of $H$ that  contain some edge of $J$.

\vspace{2mm}

For sets $X, Y\subseteq V$ with $Y$ not necessarily contained in $X$, we use $X\setminus Y$ to denote $X \cap (V\setminus Y)= X \setminus (X\cap Y)$. For graphs $\tilde{G}$ and $G$ on the same vertex set with $\tilde{G}$ a subgraph of $G$, we let $G\setminus \tilde{G}$ denote the graph on $V(G)$ given by the set of edges that feature in $G$ but not in $\tilde{G}$. If $H'$ and $H$ are $r$-uniform hypergraphs with $H'$ a subgraph of $H$, then $H\setminus H'$ is defined similarly.

We use the notation $x=y \pm z$ to denote that $x\leq y+z$ and $x \geq y-z$ and we say a property holds with high probability (whp, for short), if the probability that it holds tends to 1 with some parameter $n$ (usually the number of vertices of a graph). Finally, we drop ceilings and floors unless necessary, so as not to clutter the arguments. 

\subsection{Properties of bijumbled graphs} \label{sec:bijumbledprops}

Here we collect some properties of bijumbled graphs. These range from simple consequences of Definition \ref{def:bijumbled} to more involved statements catered to our purposes. We begin by showing that we can assume that the graphs we consider have an arbitrarily large number of vertices. 

\begin{fact} \label{fact:largen}
Given any $3\leq r\in \NN$ and $n_0\in \NN$, there exists $\eps>0$ such that  any $n$ vertex~$(p,\beta)$-bijumbled 
graph $G$ with $n\in r \NN$, $p>0$,  $\delta(G)< (1-1/r) n$  and~$\beta \leq \eps p^{r-1}n$ must have $n\geq n_0$.
\end{fact}
\begin{proof}
Let $\eps>0$ such that $\eps<1/(2n_0r)$. Suppose for a contradiction that there exists an $n$ vertex $(p,\beta)$-bijumbled graph with $\delta(G)< (1-1/r) n$, $\beta \leq \eps p^{r-1}n$ and $n< n_0$. Then due to the upper bound on the minimum degree of $G$, there exists a vertex $u\in V(G)$ and a set $W\in V(G)\setminus \{u\}$ such that $|W|=n/r$ and $\deg^G_W(u)=0$. However, from Definition~\ref{def:bijumbled}, we have that 
\[e(\{u\},W)\geq p|W|-\eps p^{r-1}n\sqrt{\frac{n}{r}}\geq \frac{pn}{r}\left(1-\eps\sqrt{nr}\right)\geq  \frac{pn}{2r}>0,\]
a contradiction.
\end{proof}

Fact~\ref{fact:largen} shows that by choosing $\eps>0$ sufficiently small, we guarantee that any bijumbled graph $G$ we are interested in either has a large number of vertices or has $\delta(G)\geq (1-1/r)n$, in which case Theorem~\ref{thm:HajnalSzem} implies the existence of a $K_r$-factor and we are done. 
We will use this at various points in our argument and simply state that we choose $\eps>0$ sufficiently small to force $n$ to be sufficiently large.

The following well known fact states that bijumbled graphs cannot to be too sparse.

\begin{fact} \label{fact:dense}
For any $3\leq r\in \NN$ and any $C>0$, there exists an $\eps>0$ such that if $G$ is an $n$ vertex $(p,\beta)$-bijumbled graph with  $p>0$ and 
$\beta \leq \eps p^{r-1}n$, then $p\geq C n^{-1/(2r-3)}\geq Cn^{-1/3}$. 
\end{fact}

\begin{proof}

Let $\eps>0$ be such that $\eps^2\leq 1/(32C^{2r-3})$  
and small enough that we can assume that 
\begin{enumerate}[label=(\roman*)]
    \item \label{case:smalln} $n\geq 9$;
    \item \label{case:pnotconstant} $p\leq 1/16$.  
\end{enumerate}
Indeed, from Fact~\ref{fact:largen}, we can choose $\eps$ so that \ref{case:smalln} holds and $Cn^{-1/(2r-3)}<1/16$ and so we are done if we are not in case \ref{case:pnotconstant}. 
We will also restrict to the case that 
\begin{enumerate}[label=(\roman*)]
\setcounter{enumi}{2}
    \item \label{case:vsmallp}$p\geq 1/(2n)$.
\end{enumerate}
To see that we can do this, suppose for a contradiction that there exists a $(p,\beta)$-bijumbled graph $G=(V,E)$ with $pn<1/2$.  We appeal to  Definition~\ref{def:bijumbled} and upper bound $2e(G)=e(V,V)$ by $pn^2+\eps p^{r-1}n^2<n-1$.  
Hence there must be some vertex $u\in V$ which is isolated in $G$. But then defining $W:=V\setminus \{u\}$, the lower bound of Definition~\ref{def:bijumbled} gives that $e(\{u\},W)\geq p(n-1)-\eps p^{r-1}n\sqrt{n-1}\geq pn(1/2-\eps pn)>0$, a contradiction.

\vspace{2mm}

We now turn to proving the statement in full generality. Our aim is to construct large (disjoint) vertex subsets $U$ and $W$ such that $e(U,W)=0$. We do this in the following greedy fashion. We initiate a process by setting $U=\emptyset$ and $W=V(G)$. Now whilst $|W|\geq 3n/4$, there exists some $u\in W$ with $\deg_W(u)\leq 2p|W|\leq 2pn$. Indeed this follows from Definition \ref{def:bijumbled} as 
\[\sum_{w\in W}\deg_W(w)=e(W,W)\leq p|W|^2+\eps p^{r-1}|W|\leq 2p|W|^2.\] 
We then choose such a $u$, delete it from $W$ and add it to $U$ and also remove $N_W(u)$ from $W$. 

Let $U$ and $W$ be the resulting sets after this process terminates. It is clear that $e(U,W)=0$ as we have removed all the neighbours of each vertex $u\in U$ from $W$ during the process. We also claim that $|W|\geq n/2$ and $U\geq 1/(16p)$. Indeed, the last step removed at most $1+2pn$ vertices from $W$. Due to our assumptions \ref{case:smalln} and \ref{case:pnotconstant}, we have that $1+2pn< n/4$ and so as $W$ had size greater than $3n/4$ before this step, we indeed have that $|W|\geq n/2$ as the process terminates. To see the lower bound on the size of $U$, note that if this was not the case, then \[|V(G)\setminus W|=\left|\cup_{u\in U} \left(\{u\}\cup N^G(u)\right)\right|\leq \sum_{u\in U }\left|\{u\}\cup N^G(u)\right|\leq |U|(1+2pn)\leq \frac{1}{16p}+\frac{n}{8}\leq \frac{n}{4},\]
using assumption \ref{case:vsmallp} in the last inequality. This
implies then that $|W|\geq 3n/4$, a contradiction as the process terminated. 

Thus $|W|\geq n/2$, $|U|\geq 1/(16p)$
and from Definition \ref{def:bijumbled}, we have that \[0=e(U,W)\geq p|U||W|-\eps p^{r-1}n\sqrt{|U||W|},\]
implying that $p^{2r-3}\geq 1/(32\eps^2 n)$. Given our upper bound on $\eps$, this implies that  $p\geq Cn^{-1/(2r-3)}$ as required. 
\end{proof}

Our first lemma shows that few vertices have degree much smaller or much larger than expected to a given set. 
\begin{lem} \label{lem:bad degrees}
For any $3\leq r\in \NN$ and $\eta >0$ there exists an $\eps>0$  such that if $G$ is an $n$-vertex $(p,\beta)$-bijumbled graph with $\beta \leq \eps p^{r-1}n$ 
then for $W\subseteq V(G)$ we have that:
\begin{enumerate}[label=(\roman*)]
\item  \label{lem:badverticessmalldegs} The number of vertices $v\in V(G)$ such that $\deg_W(v)< p|W|/2$, is less than \[\frac{\eta   p^{2r-4}n^2}{|W|}.\]
\item \label{lem:bad vertices high degs} For any $q$  such that $2p\leq q \leq 1$, the number of vertices $v\in V(G)$ such that $\deg_W(v)> q|W|$, is less than \[\frac{\eta   p^{2r-2}n^2}
{q^2|W|}.\]
\end{enumerate}
\end{lem}

\begin{proof}

Fix $\eps>0$ such that such that~$4\eps^2<\eta $. We prove only \ref{lem:bad vertices high degs}, the proof of \ref{lem:badverticessmalldegs} is both similar and simpler. We set $B$ to be the set of `bad' vertices i.e. vertices $v$ such that $\deg_W(v)>q|W|$. Thus we have that 
\[q|B||W|<e(B,W)\leq p|B||W|+ \eps p^{r-1}n \sqrt{|B||W|},\]
using the definition of $B$ and \eqref{eq:bijumbled}. 
Rearranging gives that \[|B|<\frac{\eps^2p^{2r-2}n^2}{\left(q-p\right)^2|W|},\]
and using that $p\leq q/2$ gives the desired conclusion using our choice of $\eps$.
\end{proof}

Next, we state some further consequences of  Definition~\ref{def:bijumbled}, showing that we can find cliques traversing large enough subsets of vertices. The following lemma is very general and will be used at various points in our argument. Due to its generality, there are some technical features. Whilst these are all necessary for certain parts of our argument, we do not need all of these at once. In fact, for easy reference, we list the consequences of Lemma~\ref{lem:generaltransversalcliques} that we will use in Corollary~\ref{cor:transversalcliques}. This may also serve to digest the statement of Lemma~\ref{lem:generaltransversalcliques}, seeing how it is applied in practice.

\begin{lem} \label{lem:generaltransversalcliques}
For any $3\leq r\in \NN$ and $0<\alpha<1/2^{2r}$ there exists an $\eps>0$ such that the following holds for any $n$ vertex $(p,\beta)$-bijumbled graph $G$ with $\beta\leq \eps p^{r-1}n$. Suppose that there are integers $x_i$, $i\in[r+1]$ such that $x_1\geq \ldots \geq x_{r+1}\geq 0$ and for some $r^*\in [r]$, one has that 
\begin{equation} \label{eq:x_iconds}
    x_{i}+x_{i+1}+2i\leq 2r-2,
    \end{equation}
for all  $1\leq i \leq r^*$. Define~$y:=\max\{x_{i+1}+i:i\in[r^*]\}$. Then for any collection of subsets $U_i\subseteq V(G)$ such that $|U_i|\geq \alpha p^{x_i}n$ for all $i\in[r+1]$ and  for any subgraph $\tilde{G}$ of $G$ with maximum degree less than $\alpha^2 p^{y}n$, defining $G':=G\setminus \tilde{G}$, we have that there exists a clique $S\in K_{r^*}(G')$ traversing $U_1,\ldots,U_{r^*}$ such that \[\deg^{G'}_{U_j}(S)\geq \alpha {p}^{r^*}|U_j|\] for $r^*+1\leq j \leq r+1$.

\end{lem}

\begin{proof}
 Fix $\eps>0$ small enough to apply Lemma~\ref{lem:bad degrees}~\ref{lem:badverticessmalldegs} with $\eta:=\alpha^2/(2^{4r}r)$. Further, fix $y$ and $\tilde{G}$ as in the statement, setting $ G':=G\setminus \tilde{G}$. We will prove inductively that for $i=1,\ldots,r^*$, there exists an $i$-clique $S_i\in K_{i}(G')$ traversing $U_1,\ldots, U_i$ such that $\deg^{G'}_{U_j}(S_i)\geq (p/4)^{i}|U_j|$ for all $j$  with $i+1\leq j\leq r+1$. Note that  $S_{r^*}$ is the desired copy of $K_{r^*}$ in the statement, using that~$\alpha\leq 1/4^{r^*}$ here.    
 
 So fix some $i\in[r^*]$. If $i\geq 2$, by induction we deduce the existence of $S_{i-1}$ as claimed and for $i\leq j \leq r+1$, define $W_j\subseteq U_j$ so that $W_j:=N^{G'}_{U_j}(S_i)$. If $i=1$, we simply set $W_j:=U_j$ for all $j$. We thus have that \begin{equation} \label{eq:W_jlowerbound}
     |W_j|\geq \left(\frac{p}{4}\right)^{i-1}|U_j|\geq \alpha 4^{1-i}p^{x_j+i-1}n, 
     \end{equation} for $i\leq j\leq r+1$.  
 Now we appeal to Lemma~\ref{lem:bad degrees}~\ref{lem:badverticessmalldegs} and conclude that for each $j$ with $i+1\leq j\leq r+1$, there is some set $B_j\subset V(G)$ such that $\deg_{W_j}^G(v)\geq p|W_j|/2$ for all $v\in V(G)\setminus B_j$ and
 \begin{equation}
     \label{eq:B_jupperbound}
|B_j|\leq \frac{\eta  p^{2r-4}n^2}{|W_j|}\leq \frac{\eta  4^{i-1} p^{2r-3-i-x_{j}}n}{\alpha }\leq \frac{\alpha p^{2r-3-i-x_{i+1}}n}{4^ir}\leq  \frac{\alpha p^{x_i+i-1}n}{4^ir}\leq \frac{|W_i|}{2r}.  \end{equation} 
 Here, we used \eqref{eq:W_jlowerbound} in the second inequality, the definition of $\eta $ and the fact the $x_j\leq x_{i+1}$  in the third, \eqref{eq:x_iconds} in the fourth and \eqref{eq:W_jlowerbound} once again in the final inequality. We can thus conclude from \eqref{eq:B_jupperbound} that there exists a vertex 
  $w_i\in W_i$ such that $w_i\notin B_j$ 
  for all $i+1\leq j \leq r-1$. 
 We claim that choosing $S_i=S_{i-1}\cup \{w_i\}$ completes the inductive step. Indeed $S_i\in K_{i}(G')$ as $w_i$ was chosen from the common neighbourhood of $S_{i-1}$ in $G'$. 
 Also, fixing some $i+1\leq j\leq r-1$, we have that  $N^G(w_i)$ intersects $W_j=N^{G'}_{U_j}(S_{i-1})$ in at least $p|W_j|/2$ vertices. Furthermore, at most \[\alpha^2p^{y}n\leq \alpha^2 p^{x_{i+1}+i}n\leq\frac{\alpha}{2^{2r}}p^{x_{j}+i}n\leq p|W_j|/4\] edges adjacent to $w_i$ lie in $\tilde{G}$, using the definition of $y$, the upper bound on $\alpha$, the fact that $x_j\leq x_{i+1}$ and \eqref{eq:W_jlowerbound}. Therefore  we can conclude that for all $i+1\leq j\leq r$, we have $\deg^{G'}_{U_j}(S_i)\geq \deg^{G'}_{W_j}(S_i)\geq p|W_j|/4\geq (p/4)^i|U_j|$,  
 as required. This completes the induction and the proof. 
\end{proof}

We now collect some easy consequences of Lemma~\ref{lem:generaltransversalcliques} for reference later in the proof. 

\begin{cor} \label{cor:transversalcliques}
For any $3\leq r\in \NN$ and $0<\alpha<1/2^{2r}$ there exists an $\eps>0$ such that the following holds for any $n$ vertex $(p,\beta)$-bijumbled graph $G$ with $\beta\leq \eps p^{r-1}n$. 
We have that:
\begin{enumerate}
\item \label{cor:forbiddingsubgraph} For any subgraph $\tilde{G}$ of $G$ with maximum degree less than $\alpha^2 p^{r-1}n$ we have the following.
\begin{enumerate} [label=(\roman*)]
 \item \label{cor:r-1cliqueinpn}  For any vertex subsets $U_1, \ldots, U_{r-1}\subseteq V(G)$ such that $|U_i|\geq \alpha pn$ for $i\in[r-1]$,   there exists an $(r-1)$-clique  $S\in K_{r-1}(G\setminus \tilde{G})$, traversing the $U_i$. 
    \item \label{cor:1small} For any  $ U_1, \ldots, U_{r}\subseteq V(G)$ such that  $|U_1|\geq \alpha p^{2r-4}n$ and $|U_i| \geq \alpha n$ for~$2\leq i \leq  r$, there exists an $r$-clique $S\in K_r(G\setminus \tilde{G})$,  traversing the~$U_i$. 
    \end{enumerate}
    
    \vspace{2mm}
    
    \item \label{cor:r-2small} For any  $ U_1, \ldots, U_{r}\subseteq V(G)$ such that  $|U_1|\geq \alpha p^{r-1}n$, $|U_i| \geq \alpha p n$ for~$2\leq i \leq  r-2$ and $|U_{r-1}|,|U_r|\geq \alpha n$, 
    there exists an $r$-clique $S\in K_r(G)$,  traversing the~$U_i$.
    
    \vspace{2mm}
    \item \label{cor:popular K_r-1} For any $W_0, W_1, W_2\subseteq V(G)$ such that $|W_0|,|W_1|,|W_2|\geq \alpha n$, there exists an $S\in K_{r-1}(G[W_0])$ such that $\deg_{W_i}(S)\geq \alpha^2 p^{r-1}n$ for $j=1,2$.
    \end{enumerate}

\end{cor}

\begin{proof}
Fix $\eps>0$ small enough to apply Lemma \ref{lem:generaltransversalcliques}. This is predominantly a case of plugging in the values and checking the conditions of Lemma~\ref{lem:generaltransversalcliques}. For part~\eqref{cor:forbiddingsubgraph}, we let $G'=G\setminus \tilde{G}$. Then for \eqref{cor:forbiddingsubgraph}\ref{cor:r-1cliqueinpn}, we take $r^*=r-2$, $x_i=1$ for $1\leq i\leq r+1$ and $y=r-1$. We thus have that for $i\in[r^*]$,  $x_{i}+x_{i+1}+2i=2+2i\leq 2r-2$ and $x_{i+1}+i=1+i\leq r-1=y$. Therefore taking $U_i$ for $1\leq i\leq r-1$ with $|U_i|\geq \alpha pn$ (and defining $U_{r+1}=U_r=U_{r-1}$), Lemma~\ref{lem:generaltransversalcliques} gives us an $(r-2)$-clique $S'\in K_{r-1}(G')$ traversing $U_1,\ldots, U_{r-2}$ such that $\deg^{G'}_{U_{r-1}}(S')\geq \alpha^2p^{r-1}n>0$  (here Fact \ref{fact:dense}  shows positivity).  
Therefore choosing any vertex $v\in N^{G'}_{U_{r-1}}(S') $ and fixing $S=S'\cup\{v\}$ gives the required clique.

The other cases are similar. For part~\eqref{cor:forbiddingsubgraph} \ref{cor:1small}, we fix $r^*=r-1$, $x_1=2r-4$, $x_i=0$ for $2\leq i\leq r+1$ and $y=r-1$. Again, it easily checked that the conditions on the $x_i$ are all satisfied and so applying Lemma $\ref{lem:generaltransversalcliques}$ (fixing $U_{r+1}=U_r$) gives a $(r-1)$-clique $S'$ in $G'$ traversing $U_1,\ldots,U_{r-1}$ such that $S'$ has a nonempty $G'$-neighbourhood in $U_{r}$. Therefore adding any vertex in this neighbourhood to $S'$ gives the required $r$-clique $S\in K_r(G')$. 

For part \eqref{cor:r-2small}, we fix $r^*=r-1$, $x_1=r-1$, $x_i=1$ for all $i$ such that  $2\leq i \leq r-2$ and $x_{r-1}=x_r=x_{r+1}=0$. We also let~$\tilde{G}$ be the empty graph and so~$G=G'$. Now note that for $r=3$, we have  $x_1=2$ and $x_2=0$ and so $
x_1+x_2+2=4=2r-2$, whilst for $r\geq 4$, we have $x_1+x_2+2=r+2\leq 2r-2$. The conditions  \eqref{eq:x_iconds} for $2\leq i \leq r^*=r-1$ can be similarly checked.  Therefore Lemma~\ref{lem:generaltransversalcliques} gives an $(r-1)$-clique $S'\in K_{r-1}(G)$ traversing $U_1,\ldots,U_{r-1}$ such that $N^G_{U_r}(S')\neq \emptyset$ and so as above, we extend $S'$ to the required $r$-clique $S$.

Finally, for part \eqref{cor:popular K_r-1} we fix $r^*=r-1$, $x_i=0$ for all $1\leq i \leq r$ and define our sets as $U_i=W_0$ for $i\in[r-1]$ and $U_r=W_1$, $U_{r+1}=W_2$. Applying  Lemma~\ref{lem:generaltransversalcliques} then directly gives the required $(r-1)$-clique $S\in K_{r-1}(G[W_0])$ (again here~$\tilde{G}$ is taken to be empty).
\end{proof}

\subsection{Concentration of random variables} \label{sec:concentration}
We will use the following well known concentration bounds, see e.g. \cite[Theorem 2.1, Corollary 2.4 and Theorem 2.8]{janson2011random}

\begin{thm}[Chernoff bounds] \label{thm:chernoff}
Let $X$ be the sum of a set of mutually  independent Bernoulli random variables and let $\lambda=\EE[X]$. Then for any $0<\delta<3/2$, we have that 
\[\PP[X\geq (1+\delta)\lambda]\leq  e^{-\delta^2\lambda/3 } \hspace{2mm} \mbox{ and } \hspace{2mm} \PP[X\leq (1-\delta)\lambda] \leq  e^{-\delta^2\lambda/2 }.\]
Furthermore, if $x\geq 7 \lambda$, then $\PP[X\geq x]\leq e^{-x}$.
\end{thm}

\subsection{Perfect fractional matchings} \label{sec:fractional matchings lin programming}
Given an $r$-uniform hypergraph $H$, a  \emph{fractional matching} in $H$ is a function $f : E(H) \to \RR_{\geq 0}$ such that
 $\sum_{e:v\in e}f(e)\leq 1$ for all $v\in V(H)$. We say the fractional matching is \emph{perfect} if $\sum_{e:v\in e}f(e)= 1$ for all $v\in V(H)$. The \emph{value} of a fractional matching $f$ is $|f| := \sum_{e\in E(H)}f(e).$
 The maximum value $|f| $ over all choices of fractional matching $f$ of $H$, we call the \emph{fractional matching number of $H$}, which we denote by $\nu^*(H)$. 
 
  A \emph{fractional cover} of $H$ is a function $g:V(H) \to \RR_{\geq 0}$ such that for all $e\in E(H)$, one has $\sum_{v\in e}g(v)\geq 1$. The \emph{value} of a fractional cover $g$ is $|g|:=\sum_{v\in V(H)}g(v)$. The \emph{fractional cover number} of $H$, denoted $\tau^*(H)$ is then the minimum value of a fractional cover $g$ of $H$. 
  
 For an $r$-uniform hypergraph $H$, the fractional matching number of $H$ can be encoded as the optimal solution of a linear program. Taking the dual of this linear program gives another linear program which outputs the fractional cover number as an optimal solution. The duality theorem from linear programming thus tells us that $\nu^*(H)=\tau^*(H)$ for any hypergraph $H$.  Using this, as well as the so called `complementary slackness conditions' that follow from the duality theorem, one can derive the following simple consequences, see e.g. \cite[Proposition 2]{Kriv96} or \cite[Proposition 2.4]{HKP18b}.
 
 \begin{prop}\label{prop:fractional properties}
For any $r$-uniform hypergraph $H$ on $N$ vertices, the following hold. 
\begin{enumerate}
    \item \label{item:fractional matching upper bound} $\nu^*(H)\leq N/r$ with equality if and only if there exists a perfect fractional matching in $H$. 
    \item \label{item: fractional matching at least matching}$\nu^*(H)\geq \nu(H)$ where $\nu(H)$ denotes the size of the largest matching in $H$.
    \item \label{item:frac cover subset} If $g:V(H)\to \RR_{\geq 0}$ is a fractional cover and $U\subset V(H)$, then $g':=g|_{U}:U\to \mathbb{R}_{\geq 0}$ is a fractional cover of $H[U]$ and hence $|g'|=\sum_{u\in U}g(u)\geq \tau^*(H[U])=\nu^*(H[U])$.
    \item \label{item:frac nonzero vertices}If $g:V(H) \to \RR_{\geq 0}$ is an optimal fractional cover i.e.  $|g|=\tau^*(H)$, then $\nu^*(H)\geq |W|/r$ where $W:=\{v\in V(H):g(v)>0\}$.
\end{enumerate}
\end{prop}

We now give two lemmas, exploring some simple conditions which guarantee the existence of a perfect fractional matching.  

\begin{lem} \label{lem:vertex fans}
Suppose   $H$ is an $N$-vertex,  $r$-uniform hypergraph such that given any vertex $v\in V(H)$ and any subset $W\subseteq V(H) \setminus \{v\}$ of at least $N/(2r)$ vertices, there exists an edge in $H$ containing $v$ and $r-1$ vertices of $W$.  
Then $H$ has a perfect fractional matching.
\end{lem}
\begin{proof}
Suppose for a contradiction that $H$ does \emph{not} have a perfect fractional matching. Thus, by Proposition \ref{prop:fractional properties} \eqref{item:fractional matching upper bound}, if we take $g:V(H) \to \RR_{\geq 0}$ to be an optimal fractional cover of $H$ so that $|g|=\nu^*(H)=\tau^*(H)$ we have that $|g|< N/r$. Hence, if we order the vertices in decreasing weight order according to $g$, we have by Proposition \ref{prop:fractional properties} \eqref{item:frac nonzero vertices}, that $g(w)=0$, where $w$ is the final vertex in this order. Take $W\subset V(H)\setminus \{w\}$ to be the set of $N/(2r)$ vertices preceding $w$ in the order. Then by the condition of the lemma, there exists an edge using $w$ and $r-1$ vertices of $W$. Since $g(w)=0$, we have that there is some vertex $w'$ in $W$ with $g(w')\geq 1/(r-1)$. Therefore all vertices preceding $W$ in the order (as well as $w'$) have at least this weight and in total 
\[|g|\geq \frac{\left(N-N/r\right)}{r-1}\geq N/r,\]
a contradiction. 
\end{proof}

Given a vertex subset $U\subset V:=V(H)$ in a hypergraph $H$, a \emph{fan} \emph{focused at $U$} in $H$ is a subset $F\subset E(H)$ of edges of $H$ such that $|e\cap U|=1$ for all $e\in F$ and $e\cap e'\cap (V\setminus U)=\emptyset$ for all $e \neq e'\in F$. In words, each edge of a fan intersects $U$ in exactly one vertex and outside of $U$, the edges in a fan are pairwise disjoint. The \emph{size} of a fan is simply the number of edges in the fan.  If $U=\{u\}$ is a single vertex, we simply refer to a fan focused at $u$.

Lemma \ref{lem:vertex fans} shows that if $H$ has the property that the link of every vertex $v$ has no large independent sets,  then it must have a perfect fractional matching. In fact, it is not necessary that we need such an expansion property to hold locally at every vertex and can instead focus on subsets of vertices, if we have an added condition that every vertex has a large enough fan focused  at it.  
This is the content of the following lemma. 

\begin{lem} \label{lem:two stage fans}
Suppose $H$ is an $N$-vertex,  r-uniform hypergraph and there exists $ r \leq M \leq N/(2r)$ such that the following hold:
\begin{enumerate}[label=(\roman*)]
    \item \label{cond:small vertex fan} 
    For all $v\in V(H)$ there is a fan focused at $v$ in $H$ of size $M$.
   \item \label{cond:expansion}  For every subset $W_0\subset V(G)$ with $|W_0|=M$ and every subset $W_1 \subset V(G) \setminus W_0$ with $|W_1|\geq N/(2r)$, there exists an edge of $H$ with one vertex in $W_0$ and the other $r-1$ vertices in $W_1$.
   \end{enumerate}
Then $H$ has a perfect fractional matching.
\end{lem}
\begin{proof}
 We start by noticing that \ref{cond:expansion} leads to the following two consequences.
 \begin{enumerate}[label=(\alph*)]
      \item \label{cond:big fans} 
      For all $U \subset V(H)$ with $|U|=(r-1)M$, fixing $V':=V(H)\setminus U$  we have that for all $U'\subset V' $ such that $|U'|=M$, there is a fan of size $N/r-M$ focused at $U'$ in $H[V']$.
    \item \label{cond:large matching in two stage} Every subset of at least $N/r$ vertices of $H$ induce an edge in $H$.
 \end{enumerate}
 Indeed, for $U'$ as in \ref{cond:big fans} we can build the fan $F_{U'}$ focused at $U'$ greedily. Whilst $|F_{U'}|\leq N/r-M$ we have that $W:=V(G)\setminus (V(F_{U'}) \cup U' \cup U)$ has size at least \[N-(N/r-M)(r-1)-M-(r-1)M \geq N-(2r-1)N/(2r) \geq N/(2r),\] using that $M\leq N/(2r)$ here.   Hence we can find an edge using one vertex of $U'$ and $r-1$ vertices of $W$ which extends the fan $F_{U'}$. The condition \ref{cond:large matching in two stage} also follows easily as taking $W'$ to be a set with $N/r$ vertices, we have that for any $W'' \subset W'$ with $|W''|=M$, there is an edge containing  a vertex in $W''$ and $r-1$ vertices of $W'\setminus W''$ from \ref{cond:expansion}.

Now we turn to the main proof.  We  fix $g:V(H)\rightarrow \RR_{\geq 0}$ to be an optimal fractional cover and suppose for a contradiction that $|g|< N/r$. We deduce the existence of a vertex $w\in V(H)$ with $g(w)=0$ and a fan $F_w$ focused at $w$ of size $M$. Taking $U_1:= \bigcup \{e\setminus \{w\}:e\in F_w\}$, we have that $|U_1|=(r-1)M$ and $\sum_{u\in U_1}g(u)\geq M$.

Now consider $V':=V(H)\setminus U_1$. If $\nu^*(H[V'])\geq N/r-M$ then we can conclude that $\sum_{v\in V'}g(v) \geq N/r-M$ from Proposition \ref{prop:fractional properties} \eqref{item:frac cover subset} which implies that $|g|\geq N/r$, a contradiction. Hence \begin{equation} \label{eq:lower bound on frac cover} \nu^*(H[V'])< \frac{N}{r}-M=\frac{N'-M}{r},\end{equation}
where $N':=|V'|=N-(r-1)M$. We fix  $g': V' \rightarrow \RR_{\geq 0}$ to be some optimal fractional cover of $H[V']$ with $|g'|=\nu^*(H[V'])$.  By Proposition \ref{prop:fractional properties} \eqref{item:frac nonzero vertices}, we therefore have that there is some set $U_2\subset V'$ with $|U_2|=M$ and $g'(u')=0$ for all $u'\in U_2$.  By \ref{cond:big fans} there exists a fan $F_{U_2}$ of size $N/r-M$ focused at $U_2$ in $H[V']$. Taking $Z:= \bigcup \{e:e\in F_{U_2}\}\setminus U_2$, we have that $|Z|= (r-1)(N/r-M)$ and similarly to  before, using that for each edge $e\in F_{U_2}$ we have $\sum_{v\in e} g'(v) \geq 1$ and the fact that $g'(u')=0$ for all $u'\in U_2$, we can conclude that $\sum_{z\in Z}g'(z) \geq |F_{U_2}|=N/r-M$. 

Finally, we look at $V'':=V'\setminus Z$. We have that $N'':=|V''|=N'-(r-1)(N/r)+(r-1)M$  and using \ref{cond:large matching in two stage} and Proposition \ref{prop:fractional properties} \eqref{item: fractional matching at least matching}, we have that \[ \nu^*(H[V''])\geq \frac{N''-(N/r)}{r}=\frac{N'+(r-1)M}{r}-\frac{N}{r}.\]
Hence, by Proposition \ref{prop:fractional properties} \eqref{item:frac cover subset}, we can conclude that $\sum_{v''\in V''}g'(v'')\geq (N'+(r-1)M)/r-N/r$. Combining this with the lower bound on the sum of $g'$ values on $Z$ implies that $|g'|=\nu^*(H[V'])\geq (N'-M)/r$, contradicting \eqref{eq:lower bound on frac cover}. 
\end{proof}

\subsection{Almost perfect matchings in hypergraphs}
\label{sec:almost perfect matchings}
It is well-known that hypergraphs that have roughly regular vertex degrees and small codegrees contain large matchings. This is often referred to as Pippenger's Theorem but there are in fact a family of similar results, all following from the ``semi-random" or ``nibble" method, see e.g.~\cite[Section 4.7]{alon2004probabilistic}. Here we use the following explicit version which follows directly  
from a result of Kostochka and R\"odl \cite{kostochka1998partial}.
\begin{thm} \label{thm:almost perfect matchings}
For any integers $r\geq 3$ and $K\geq 4$ there exists  $\Delta_0>0$ such that for all $\Delta\geq\Delta_0$ the following holds. If $H$ is a $r$-uniform hypergraph on $N$ vertices such that:
\begin{enumerate}
    \item \label{item:almosti} for all vertices $v\in V(H)$, we have $\deg(v)=\Delta\left(1 \pm K\sqrt{\frac{\log \Delta}{\Delta}}\right)$ and 
    \item \label{item:almostii} for all $u\neq v\in V(H)$, we have $\codeg(u,v)\leq  \Delta^{1/(2r-1)}$, 
\end{enumerate}
then $H$ has a matching covering all but at most $  \Delta^{-1/r}N$ vertices. 
\end{thm}

Indeed,~\cite[Theorem 4]{kostochka1998partial} states that for all $r\geq 3$, $K_0\ge 8$ and reals $0<\delta,\gamma<1$, there exists a~$D_0$ such that if ${H}$ is an $r$-uniform hypergraph on $N$ vertices with  \[D-K_0\sqrt{D \log D}\le \deg^{ H}(v)\le D, \] 
for all $v\in V( H)$, where $D\ge D_0$, and  $\codeg^{ H}(u,v)\le C< D^{1-\gamma}$ for all pairs of vertices~$u\neq v$, then $ H$ has a matching covering all but at most $O\left(N\left(\tfrac{C}{D}\right)^{(1-\delta)/(r-1)}\right)$ vertices. In order to derive~\Cref{thm:almost perfect matchings} from this we fix~$K_0=2K$,~$\delta=\tfrac{1}{4r}$ and~$\gamma=\tfrac{2r-2}{2r-1}$. Letting $D_0$ be the resulting constant given by~\cite[Theorem 4]{kostochka1998partial}, we fix $\Delta_0\ge D_0$ to be some large constant. Hence, our conditions~\eqref{item:almosti} and~\eqref{item:almostii} of Theorem~\ref{thm:almost perfect matchings}  guarantee that $ H$ satisfies the conditions of~\cite[Theorem 4]{kostochka1998partial} with~$D=\Delta +K\sqrt{\Delta\log \Delta}$ and~$C=\Delta^{1-\gamma}$. Now note that \[\frac{C}{D}=(1+o(1))\Delta^{-\gamma}=(1+o(1))\Delta^{-(2r-2)/(2r-1)}=o( \Delta^{-(4r-4)/(4r-1)}).\] Combined with the fact that  $\tfrac{1-\delta}{r-1}=\tfrac{4r-1}{4r(r-1)}$ and $\Delta\ge \Delta_0$ is sufficiently large, it follows from~\cite[Theorem 4]{kostochka1998partial}  that the number of vertices uncovered by a largest matching is always less than~$\Delta^{-1/r}N$, as required.

Clearly, in order to prove that a hypergraph $H$ has a large matching, it suffices to establish the conditions of Theorem \ref{thm:almost perfect matchings} for a spanning subgraph $H' \subset H$. An idea introduced by Alon, Frankl, Huang, R\"odl, Ruci\'nski and Sudakov \cite{AFHRRS12}  is to find such an $H'$ as a random subhypergraph of $H$ and guarantee that the conditions of Theorem \ref{thm:almost perfect matchings} hold for $H'$ by using perfect fractional matchings to  dictate the probability with which we take each edge into $H'$. This idea was then used in the context of finding almost $K_r$-factors in pseudorandom graphs by Han, Kohayakawa and Person~\cite{HKP18a,HKP18b}. We will also adopt this idea and so give the following theorem. 

\begin{thm} \label{thm:almost perfect matchings from fractional matchings}
For all $3\leq r\in \NN$ and $0<\eta<1/2$, there exists an $N_0$ such that the following holds for all $N\geq N_0$. Suppose $H$ is an $N$-vertex,  $r$-uniform hypergraph such that there exist $t:=2N^\eta$ perfect fractional matchings  $f_1,\ldots,f_t:E(H) \rightarrow \RR_{\geq 0}$ in $H$  with the property that 
\begin{equation} \label{eq:upper bound on sum of frac values} \sum_{i=1}^t \sum_{e\in E(H):\{u,v\}\subset e}  f_i(e) \leq 2,\end{equation}
for all pairs of vertices $u\neq v\in V(H)$. Then $H$ has a matching covering all but at most $N^{1-\eta/r}$ vertices.
\end{thm}
\begin{proof}
We take a random subgraph $H'\subseteq H$ be keeping every edge $e\in E(H)$ independently with probability
$p_e= \sum_{i=1}^tf_i(e)/2$ noting that $ p_e \in [0,1]$ for all $e\in E(H)$ due to \eqref{eq:upper bound on sum of frac values}. We fix $\Delta:=t/2=N^\eta$ and $K:=4/\eta$  and claim that $H'$ satisfies the conditions of Theorem \ref{thm:almost perfect matchings} whp as  $N$ tends to infinity. 

To check that $H'$ satisfies the conditions of Theorem \ref{thm:almost perfect matchings}, note that for each $v\in V$ we have \[\EE\big[\deg^{H'}(v)\big]=\sum_{e:v\in e} p_e=\sum_{e:v\in e} \sum_{i=1}^t f_i(e)/2=\frac{1}{2}\sum_{i=1}^t \left(\sum_{e:v\in e}f_i(e)\right)=t/2=\Delta,\]
using that each $f_i$ is a perfect fractional matching. Applying Theorem \ref{thm:chernoff} then gives that 
\begin{align}
     \nonumber \PP\left[\deg^{H'}(v)\neq \Delta\left(1 \pm K\sqrt{\frac{\log \Delta}{\Delta}}\right)\right] &\leq 2 \exp\left(-\frac{K^2 \log \Delta}{3}\right)
     \\  \label{eq:small prob bad degrees} &\leq 2 \exp \left(- \frac{K^2\eta \log N}{3} \right)
     \\ \nonumber &\leq \frac{1}{N^2}, \end{align}
for $N$ sufficiently large. 
Similarly, for $u\neq v\in V(H)$, we have that $\EE\left[\codeg^{H'}(u,v)\right]= \sum_{e:\{u,v\}\subset e}p_e\leq 1$ by \eqref{eq:upper bound on sum of frac values} and applying Theorem \ref{thm:chernoff} gives that 
\begin{equation} \label{eq:small prob bad codegrees}\PP\left[\codeg^{H'}(u,v) \geq \Delta^{1/(2r-1)}\right] \leq \exp \left(- \Delta^{1/(2r-1)}\right) \leq \frac{1}{N^3}, \end{equation}
for large $N$. 
Hence taking a union bound over all vertices and pairs of vertices and upper bounding the failure probabilities 
with \eqref{eq:small prob bad degrees} and \eqref{eq:small prob bad codegrees} gives that $H'$ satisfies the conditions of Theorem \ref{thm:almost perfect matchings} whp. Therefore for $N$ (and hence $\Delta$) sufficiently large, we can fix such an instance of $H'$ and apply Theorem \ref{thm:almost perfect matchings} which gives the large matching in $H'$ and hence in $H$, concluding the proof. \end{proof}

It will be useful for us to work with the following corollary to Theorem~\ref{thm:almost perfect matchings from fractional matchings} which gives us a sufficient condition for us to be able to generate the perfect fractional matchings needed in Theorem~\ref{thm:almost perfect matchings from fractional matchings} via a greedy process. Recall that for a $2$-uniform graph $J$ on $V(H)$, $H_J$ denotes the subhypergraph of $H$ given by all edges of $H$ which contain some edge of $J$.

\begin{thm} \label{thm:almost perfect matchings via greedy pfms}
For all $3\leq r\in \NN$ and $0<\gamma<1/2r^2$,  
there exists an $N_0$ such that the following holds for all $N\geq N_0$. Suppose $H$ is an $N$-vertex,  $r$-uniform hypergraph such that given any  \emph{graph} $J$ on $V(H)$ of maximum degree at most $N^{r^2\gamma}$,  we have that  $H\setminus H_J$ contains a perfect fractional matching. 
Then $H$ has a matching covering all but at most $N^{1-\gamma}$ vertices. 
\end{thm}
\begin{proof}
We will prove this by appealing to Theorem \ref{thm:almost perfect matchings from fractional matchings} with $\eta:=r\gamma$ and so we set out to find $t:= 2 N^{\eta}$ perfect fractional matchings $f_1,\ldots,f_t$ such that \eqref{eq:upper bound on sum of frac values} holds. We do this algorithmically, finding the $f_i$ one at a time. We  begin by defining $J_1$ to be the empty ($2$-uniform) graph on $V(H)$ and for $1 \leq i\leq t$ we do the following. We find a perfect fractional  matching $f_i$ in $H\setminus H_{J_{i}}$ 
and add this to our family of perfect fractional matchings. We then define a graph $G_i$ with vertex set $V(H)$ and a pair of vertices $\rho \in \binom{V(H)}{2}$ forming an edge in $G_i$ if 
\[\sum_{\rho \subset e \in E(H)}f_i(e)\geq \frac{N^{-\eta}}{2}.\]
Finally we define 
$J_{i+1}:= J_{i} \cup G_i$ and move to step $i+1$. 

We claim that this algorithm does not stall and we complete our collection of $t$ perfect fractional matchings. In order to check this, we need to verify that 
we can find a perfect fractional matching in $H\setminus H_{J_{j}}$ for each $j\in[t]$. 
This follows because at each step $i$, we have that  for any $v\in V(H)$, 
\[\sum_{u\in V(H) \setminus \{v\}} \left(\sum_{\{u,v\}\subset e \in E(H)}f_i(e)\right)= (r-1) \sum_{v\in e\in E(H)}f_i(e)=r-1,\]
as $f_i$ is a perfect fractional matching. Hence
the number of pairs $\rho \in \binom{V(H)}{2}$  which contain $v$  and form  an edge of $G_i$ is at most $2(r-1)N^{\eta}$. 
As this holds for all choices of  $v\in V(H)$ 
we have that $G_i$ has maximum degree less than $2(r-1) N^{\eta}$. Thus for each $j\in[t]$,  $J_j:=\bigcup_{i=1}^{j-1}G_i$ has maximum degree less than 
\[2(r-1) N^{\eta} \cdot (j-1)\leq 2(r-1)N^{\eta} t= 4(r-1) N^{2\eta}\leq N^{r\eta}= N^{r^2\gamma},\]
for $N$ sufficiently large. 
So $H\setminus H_{J_i}$ does indeed  host a perfect fractional matching by assumption. 

Finally  we need to check condition \eqref{eq:upper bound on sum of frac values} for each pair of vertices $\rho=\{u,v\}\in \binom{V(H)}{2}$. Note that for any pair $\rho\in \binom{V(H)}{2}$ of vertices of $H$ we have that \[\sum_{i=1}^t \left(\sum_{\rho \subset e \in E(H)}f_i(e)\right) \leq \frac{tN^{-\eta}}{2} \leq 1\]
if $\rho$ does not feature as an edge in  any of the $G_i$. On the other hand, we have that if $\rho=\{u,v\}\in E(G_j)$ for some $j\in [t]$, then note that because we forbid the edges of $H$ containing $\{u,v\}$ 
from being used again we have that
\[\sum_{i=j+1}^t \left(\sum_{\rho\subset e \in E(H)}f_i(e)\right)=0 .\]
Also we have that $\rho \notin E(G_i)$ for $i < j$ as otherwise there could be no weight on (edges containing) $\rho $ in $f_j$. Hence 
\[\sum_{i=1}^{j-1} \left(\sum_{\rho\subset e \in E(H)}f_i(e)\right)\leq \frac{(j-1)N^{-\eta}}{2} \leq 1,\]
and using that \[ \sum_{\rho \subset e \in E(H)}f_j(e) \leq \sum_{u \in e \in E(H)}f_j(e)=1 ,\]
gives that \eqref{eq:upper bound on sum of frac values} holds for all $\rho \in \binom{V(H)}{2}$ as required. So by Theorem \ref{thm:almost perfect matchings from fractional matchings}, we have that $H$ contains a matching covering all but at most $N^{1-\eta/r}=N^{1-\gamma}$ vertices, concluding the proof. 
\end{proof}

\subsection{Templates} \label{sec:templates}

In this section we concentrate on a  powerful new  approach introduced by  Montgomery~\cite{M14a,M19}, in his work on spanning trees in random graphs. The general idea is to use the following key notion as an auxiliary graph to define absorbing structures in the host graph of interest.

\begin{figure}[h]
    \centering
  \includegraphics[scale=0.84]{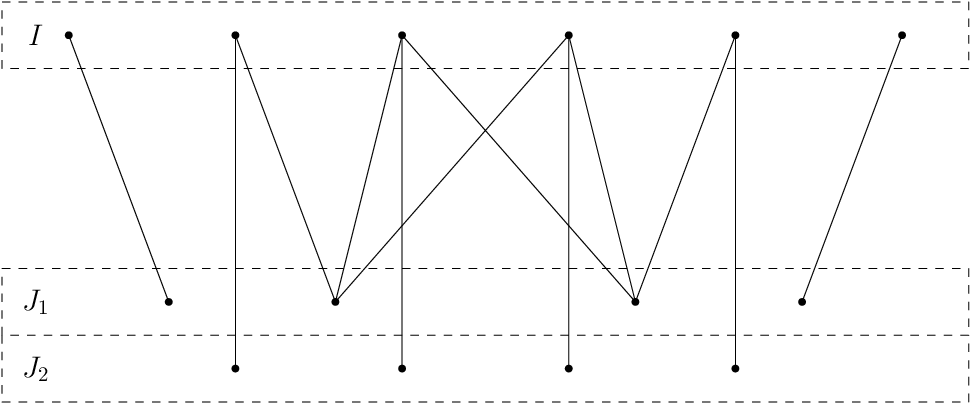}
    \caption{   \label{fig:template} A template $\pzc{T}$ of flexibility $2$. One can check that the key property is indeed  satisfied.} 
  \end{figure}

\begin{dfn} \label{def:template}
A  \emph{template} $\pzc{T}$ with \emph{flexibility} $t\in \mathbb{N}$ is a
bipartite graph on $7t$ vertices with vertex classes $I$ and
$J_1\cup J_2$, such that $|I|=3t$, $|J_1|=|J_2|=2t$, and for any
$\overline{J}\subset J_2$, with $|\overline{J}|=t$, the induced graph
$\pzc{T}[V(\pzc{T})\setminus\overline{J}]$ has a perfect matching. We call $J_2$ the
\textit{flexible} set of vertices for the template.
\end{dfn}

See Figure~\ref{fig:template} for an example of a template. The definition implies that a template is \emph{robust} with respect to having a perfect matching. It is not hard to come up with examples of templates, indeed a complete bipartite graph certainly satisfies the condition. 
The utility of the notion  for defining absorbing structures that are possible to find in the desired host graphs, comes  with the fact that  \emph{sparse templates} exist.  Indeed, Montgomery~\cite{M14a,M19} proved the following using a probabilistic argument. 
\begin{thm} \label{thm:montyexist}
For all sufficiently large~$t$, there exists a template of flexibility~$t$ and maximum degree~$40$.
\end{thm}  
Han, Kohayakawa, Person and the author~\cite{HKMP18b} then showed how to derandomise the argument for the existence of templates and find templates with bounded maximum degree efficiently in polynomial time. We will use the method of template absorption in proving Proposition \ref{prop:finalabsorption}.

\section{Diamond trees} \label{sec:diamondtree}

Recall the definition of diamond trees from Section~\ref{sec:proofmain}, namely Definition~\ref{def:diamondtree}.  
In this section we prove the existence of  diamond trees in our bijumbled graphs.  The main aim is to prove the following  proposition which  gives us some flexibility over which vertices  feature as removable vertices of our diamond tree. This will turn out to be very valuable at various points in our proof. 

\begin{prop} \label{prop: choosing removable vertices}

For any $3\leq r\in \NN$ and $0<\alpha<1/2^{2r}$ there exists an $\eps>0$ such that the following holds for any $n$ vertex $(p,\beta)$-bijumbled graph $G$ with $\beta\leq \eps p^{r-1}n$. For any $2\leq z\leq \alpha n$ and any pair of disjoint vertex subsets $U, W\subset V(G)$  such that $|U|, |W|\geq 4\alpha r n$, 
there   exist disjoint vertex subsets $X, Y  \subset U$ such that the following hold:
\begin{enumerate} 
    \item \label{property:size right} $|X|+|Y|=z$;
  \item \label{property:small X}  $|X|\leq \max\{1, 2z/d_*\}$ with $d_*=\alpha^2p^{r-1}n$;
    \item \label{property:flexibility} for \emph{any} subset $Y'\subset Y$, there exists a $K_r$-diamond tree $\cD=(T,R,\Sigma)$ such that $R=X\cup Y'$  and $\Sigma\subset K_{r-1}(G[W])$ is a matching of $(r-1)$-cliques in $W$.
\end{enumerate}

\end{prop}

Let us pause to digest the proposition. Firstly, note that by choosing $Y'=Y$ in \eqref{property:flexibility} and varying $z$, we can guarantee the existence of $K_r$-diamond trees of any order up to linear in our bijumbled graph $G$. However, the proposition is much more powerful than just this. The  vertex set $Y$ and property \eqref{property:flexibility} allow us \emph{flexibility} in which vertices appear in the removable set of vertices of the diamond tree we take from the proposition. We can start with $z$ much larger than the desired order of the diamond tree we want and then remove unwanted vertices from $Y$ to end up with some $Y'$ that we include in the removable vertices of the diamond tree. The point is that by starting with a larger $z$ (and hence larger $|Y|$), we can deduce stronger properties about the vertices in $Y$, allowing us to then ensure properties of the set of removable vertices $R$ that we would otherwise have no hope in guaranteeing. 
There is a catch, as we are forced to include the set $X$ in any diamond tree we produce, but note that due to property \eqref{property:small X}, the size of $X$ is negligible compared to the size of $Y$. Indeed due to Fact~\ref{fact:dense}, we have that $d_*$ is polynomial in $n$ (of order at least $\Omega(n^{(r-2)/(2r-3)})$, to be precise). Thus we can choose $Y'$ to be much smaller than $Y$ and still have the vertices in $Y'$ contribute a significant subset (at least half, say) of the removable vertices of the diamond tree we obtain.
We delay applications of Proposition~\ref{prop: choosing removable vertices} to later in the proof but refer the reader to Lemmas~\ref{lem:diamond trees with low degrees},~\ref{lem:diamond trees with many transversal triangles} and \ref{lem:diamond tree that intersects many sets} for a flavour of the consequences of the proposition. 

The rest of this section is concerned with proving Proposition~\ref{prop: choosing removable vertices}.
The idea behind the  proof  is simple, we look to find a large (order $z$) $K_r$-diamond tree in $G$ with the property that many of the removable vertices  are leaves (the set $Y$). This allows us to pick and choose which leaves (the set $Y'$) we include in our desired diamond tree, as we can simply remove the other leaves and their corresponding interior cliques\footnote{That is, for each unwanted leaf $v\in R$ in the diamond tree $\cD=(T,R,\Sigma)$, we remove the interior clique in $\Sigma$ which  corresponds to the edge adjacent to the (preimage of) $v$ in the defining tree $T$, as well as the leaf $v$ itself.}, see Figure~\ref{fig:scattered} for an example. 
In order to find diamond trees with many leaves, we introduce the notion of a scattered diamond tree and deduce the existence of such diamond trees in a suitably pseudorandom graph.

\subsection{Scattered diamond trees} \label{sec:scattered}

One way to ensure a large set of leaves in a tree is to impose a minimum degree on all non-leaf vertices.  This leads to the following definition. 

\begin{dfn} \label{def:scattered}
We say a tree $T$ (of order at least $2$) is \emph{$d$-scattered} if every vertex in $V(T)$ which is not a leaf in $T$ has degree at least $d$. As a convention we will also say that a tree of order 1 (a single vertex) is $d$-scattered for all $d$. We say a diamond tree is $\cD=(T,R,\Sigma)$ is $d$-scattered if its underlying auxiliary tree $T$ is $d$-scattered. 
\end{dfn}

\begin{figure}[h]
    \centering
  \includegraphics[scale=1]{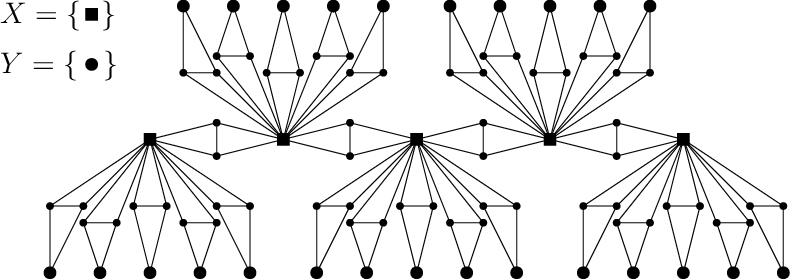}
    \caption{   \label{fig:scattered} A $6$-scattered $K_3$-diamond tree. } 
  \end{figure}

See Figure~\ref{fig:scattered} for an example of a scattered $K_3$-diamond tree.The following simple lemma shows that most of the vertices in a scattered tree (and hence most of the removable vertices in a scattered diamond tree) are leaves. 

\begin{lem} \label{lem:mostly leaves}
Let $d\geq 2$ and suppose that $T$ is a $d$-scattered tree of order $m\geq 3$. Then defining $X\subset V(T)$ to be the vertices\footnote{That is, $X$ is the set of vertices of $T$ which are not leaves.} which have degree greater than 1 in $T$, we have that $|X| \leq \frac{m-2}{d-1}$. 
\end{lem}
\begin{proof}
By the definition of $d$-scattered trees, we have that every vertex in $X$ has degree at least $d$. We define $x:=|X|$. Note that $T[X]$ is a connected subtree of $T$. Indeed, the interior vertices of a path between any two vertices of $T$ must lie in $X$ (as they have degree at least 2). Hence $T[X]$ has exactly $x-1$ edges and we can estimate the number of edges in $T$ as follows:
\[ e(T)=m-1= \sum_{v\in X} \deg(v) - e(T[X]) \geq xd-(x-1).\]
Rearranging, one obtains that $x\leq \frac{m-2}{d-1}$, as required. 
\end{proof}

We will show that we can find large scattered diamond trees in our bijumbled graph. To begin with, we focus on diamond trees for which the auxiliary tree is a star, which we call \emph{diamond stars}. 
The next lemma shows that we can find large diamond stars in a suitably pseudorandom graph.

\begin{lem} \label{lem:diamondstars}
For any $3\leq r\in \NN$ and $0<\alpha<1/2^{2r}$ there exists an $\eps>0$ such that the following holds for any $n$ vertex $(p,\beta)$-bijumbled graph $G$ with $\beta\leq \eps p^{r-1}n$, fixing $d_*:=\alpha^2p^{r-1}n$. 
Let $U_0,U_1,U_2\subseteq V(G)$ be disjoint vertex subsets such that $|U_i|\geq \alpha n$ for $i=0,1$ and $|U_2|\geq \alpha r n$. Then there exists a $K_r$-diamond tree $\cD^*=(T^*,R^*,\Sigma^*)$ in $G$ such that $T^*$ is a star of order $1+d_*$ centred at $x$, say, with\footnote{Here $\rho^*:V(T^*)\rightarrow R^*$ is the associated bijection in the definition of $\cD^*$.} $\rho^*(x)\in U_0$, $R^*\setminus \{\rho^*(x)\}\subset U_1$ and $\Sigma^*\subset K_{r-1}(G[U_2])$  a matching of $(r-1)$-cliques  in $U_2$.  
\end{lem}
\begin{proof}
Fix $\eps>0$ small enough to apply Corollary~\ref{cor:transversalcliques} \eqref{cor:popular K_r-1}. Shrink $U_0$ (if necessary) to be a set of exactly $\alpha n$ vertices.
We claim that there is a matching $M\subset K_{r-1}(G[U_2])$ of $(r-1)$-cliques such that $|M|=\alpha n$ and each clique $S\in M$ has $\deg_{U_i}(S)\geq d_*,$ 
for $i=0,1$. Indeed, we can find $M$ greedily by applying Corollary~\ref{cor:transversalcliques} \eqref{cor:popular K_r-1} (with $W_i=U_{2-i}$ for $i=0,1,2$) repeatedly, adding a $(r-1)$-clique $S$ to $M$ and removing its vertices from $U_2$ after each application. While $|M| \leq \alpha n$, we have that $|U_2|\geq \alpha n$ and so we are indeed in a position to apply Corollary~\ref{cor:transversalcliques} \eqref{cor:popular K_r-1} throughout the process. 

Now once we have found $M$, for each $S\in M$ and for $i=0,1$, let $N_i(S):=N_{U_i}(S)$, that is, the set of  vertices in $U_i$ which form a $K_r$ with $S$. By construction we have that $|N_0(S)|\geq d_*$  for each  $S$ in $M$ 
and so 
\[|\{(v,S)\in U_0\times M:v\in N_0(S)\}|\geq |M|d_*=\alpha nd_*.\]
Hence, as $|U_0|$ has size $\alpha n$ (as we imposed at the start of the proof), by averaging, there exists a vertex~$v_0\in U_0$ and a subset~$\Sigma^*$ of~$d_*$  cliques in~$M$ such that~$v_0$ is in~$N_0(S)$ for all~$S\in \Sigma^*$. We can now construct our diamond star greedily, with $v_0$ as the image of the large degree vertex. 
Sequentially, for each clique $S$ in $\Sigma^*$, choose a vertex $u$ in $N_1(S)$ which has not been previously chosen and add the copy of $K_{r+1}^-$ on $S$, $v_0$ and $u$ to the diamond star (adding $u$ to $R^*$). As $N_1(S)\geq d_*$ for all $S\in \Sigma^*\subseteq M$, there is always an option for $u$ and so this process succeeds in building the required diamond star.
\end{proof}

Our next lemma follows the scheme of Krivelevich~\cite{kri97} to construct  large diamond trees. We adapt his proof to guarantee that the diamond tree obtained is scattered.  

\begin{lem} \label{lem:scattereddiamondtrees}
For any $3\leq r\in \NN$ and $0<\alpha<1/2^{2r}$ there exists an $\eps>0$ such that the following holds for any $n$ vertex $(p,\beta)$-bijumbled graph $G$ with $\beta\leq \eps p^{r-1}n$, fixing $d_*:=\alpha^2p^{r-1}n$. 
For any $2\leq z\leq \alpha n$ and any pair of disjoint vertex subsets $U, W\subset V(G)$  such that $|U|, |W|\geq 4\alpha r n$, 
there  exists a $d_*$-scattered  $K_r$-diamond tree $\cD_{sc}=(T_{sc},R_{sc},\Sigma_{sc})$ of order $m$ such that $z\leq m \leq z+d_*$,  $R_{sc}\subset U$ and $\Sigma_{sc}\subset K_{r-1}(G[W])$ is a matching of  $(r-1)$-cliques in $G[W]$. 
\end{lem}
\begin{proof}
Our proof is algorithmic and works by building a diamond tree forest, that is, a set of pairwise  vertex disjoint diamond trees. At each step of the algorithm, we will add to one of the trees in our forest, boosting the degree of a vertex in the underlying auxiliary tree by $d_*$, using Lemma \ref{lem:diamondstars}. By discarding trees when the sum of the orders of the trees gets too large, we will show that one of the trees in our forest will eventually obtain the desired order after finitely many steps of the algorithm. The details follow.

\vspace{2mm}

Initiate the process by fixing $U_0\subset U$ to be an arbitrary subset of $\alpha n$ vertices, $W_0=\emptyset \subset W$ to be empty and $\cD_1,\ldots,\cD_{\ell}$ with $\ell=\alpha n$ to be the diamond trees which are defined to be the single vertices in $U_0$. That is, for $i \in [\ell]$, the  $K_r$-diamond tree $\cD_i=(T_i,R_i,\Sigma_i)$ corresponds to an auxiliary tree $T_i$ which is just a single vertex and thus $R_i$ is also a single vertex and $\Sigma_i$ is empty. In general, at each step of the process we will have a family $\cD_1,\ldots, \cD_\ell$ (for some $\ell \in \mathbb{N}$) of vertex disjoint $K_r$-diamond trees such that for each~$i$, the diamond tree $\cD_i=(T_i,R_i,\Sigma_i)$ is $d_*$-scattered, has $R_i \subset U_0$ and $\Sigma_i \subset G[W_0]$.
Furthermore, we will have that $U_0=\bigcup_{i\in [\ell]}R_i$ and $W_0=\bigcup_{i\in \ell}V(\Sigma_i)\subset W$ and maintain throughout that $\alpha n \leq |U_0| \leq 2 \alpha n$ and~$|W_0|\leq 2(r-1) \alpha n$

Now at each step, given such a set $U_0$ and family $\cD_1,\ldots, \cD_\ell$, we apply Lemma \ref{lem:diamondstars} with $U_1=U \setminus U_0$ and $U_2=W \setminus W_0$, noting that the conditions on the size of $U$ and $W$ in the statement of the lemma and the imposed conditions on the size of $U_0$ and $W_0$ throughout the process indeed allow Lemma \ref{lem:diamondstars} to be applied.
Thus, we find a $K_r$-diamond star $\cD^*=(T^*,R^*,\Sigma^*)$ of order $d_*+1$  with centre $v_0\in U_0$, $R^* \setminus \{v_0\}\subset U \setminus U_0$ and $\Sigma^*\subset K_{r-1}(G[U_2])$ a matching of  $(r-1)$-cliques. 
As $U_0$ is the union of the removable vertices of the family of diamond trees, we have that there is some $i_0\in[\ell]$ such that $v_0 \in R_{i_0}$. We then update $\cD_{i_0}$ by adjoining the diamond star to the tree at $v_0$, we add all the vertices of $R^*$ to $U_0$ and all the vertices of the $(r-1)$-cliques in $\Sigma^*$, to $W_0$. Now if there is a $K_r$-diamond tree among the (new) family $\cD_1, \ldots, \cD_\ell$ which has order at least $z$, we take such a diamond tree as $\cD_{sc}$ and finish the process. If not, then we look at the size of $U_0$. If $|U_0|< 2\alpha n$, we continue to the next step. If $|U_0|\geq 2\alpha n$, then we sequentially discard arbitrary  $K_r$-diamond trees $\cD_j=(T_j,R_j,\Sigma_j)$ from the family. That is,  we choose a $\cD_j$ in the family,  delete $R_j$ from $U_0$ and delete the vertices that belong to $(r-1)$-cliques in~$\Sigma_j$ from $W_0$. We continue discarding diamond trees until $|U_0|\leq 2\alpha n$. Note that as $|R_j| \leq z \leq \alpha n$ for all $j$, the updated $U_0$ at the end of this discarding process will have size at least $\alpha n$ as required. We then move to the next step.

All the diamond trees in our family are $d_*$-scattered throughout the process and also  $W_0$, as the set of vertices featuring in interior cliques of a family of $K_r$-diamond trees whose orders add up to less than $2\alpha n$, has size less than $2(r-1)\alpha n$ throughout. It is also clear that as the order of any diamond tree in our collection grows by at most $d_*$ in each step, the order of the diamond tree which is found by the algorithm will be at most $z+d_*$. It only remains to check that the algorithm terminates but this  is guaranteed because the number of diamond trees is decreasing throughout the process. Indeed, we never add new diamond trees to the family and every $\alpha n/d_*$ steps we have to discard at least one diamond tree from the family. If the algorithm does not terminate after finding an appropriate $\cD_{sc}$, then eventually we will be left with just one diamond tree $\cD_1$ in the family, but at this point the order of $\cD_1$ would be at least $\alpha n\geq z$, contradicting that the algorithm is still running. 
\end{proof}

Using Lemmas \ref{lem:diamondstars} and  \ref{lem:scattereddiamondtrees} we can now deduce Proposition \ref{prop: choosing removable vertices}.

\begin{proof}[Proof of Proposition \ref{prop: choosing removable vertices}]
Fix $\eps>0$ small enough to apply Lemmas \ref{lem:diamondstars} and  \ref{lem:scattereddiamondtrees} and small enough to force $n$ to be sufficiently large in what follows. 
Let us first deal with the case when $z\leq d_*:=\alpha^2p^{r-1}n$. Here, we arbitrarily partition $U$ into $U_0$ and $U_1$ of size at least $\alpha n$, fix $U_2=W$ and  apply Lemma~\ref{lem:diamondstars} to get a $K_r$-diamond star $\cD^*=(T^*,R^*,\Sigma^*)$ of order $1+d_*$ with $R^*\subset U$ and $\Sigma^*\subset K_{r-1}(G[W])$ a matching  of  $(r-1)$-cliques in $W$. Let $x\in R^*$ be the only non-leaf vertex in $R^*$ and define $X=\{x\}$. Further, let $Y\subset R^*\setminus X$ be an arbitrary subset of $z-1$ vertices. Now taking $\rho^*:V(T^*)\rightarrow R^*$ and  $\sigma^*:E(T^*)\rightarrow \Sigma^*$  to be the  defining bijective maps for  $\cD^*$,  note that for any $Y'\subset Y$, the set of vertices $\{\rho^{*-1}(v):v\in Y'\cup X\}\subset V(T^*)$ span a sub-tree (or rather a sub-star) of $T^*$, say $T$. Therefore, taking $\cD=(T,X\cup Y',\Sigma)$ where $\Sigma:=\{\sigma^*(e):e\in E(T)\}$ defines a  $K_r$-diamond tree with removable vertices $Y'\cup X$. Therefore 
\eqref{property:size right}, \eqref{property:small X} and \eqref{property:flexibility} of the proposition are all satisfied.

When $d_*<z\leq \alpha n$, the proof is similar. We apply Lemma \ref{lem:scattereddiamondtrees}  to get a $d_*$-scattered $K_r$-diamond tree $\cD_{sc}=(T_{sc},R_{sc},\Sigma_{sc})$ as given by the lemma and define $X\subset R_{sc}$ to be the non-leaves of $\cD_{sc}$. See Figure~\ref{fig:scattered} for an example. In order to bound $|X|$ and prove property \eqref{property:small X}, we appeal to Lemma \ref{lem:mostly leaves} which gives that
\[|X|\leq \frac{|R_{sc}|-2}{d_*-1}\leq \frac{z+d_*-2}{d_*-1}\leq \frac{2z}{d_*},\]
using that $z\geq d_*$ in the final inequality. 

We note that for $n$ large (using Fact~\ref{fact:dense}) we have that $d_*\geq 4$, implying that $|X|\leq z/2$. We fix $Y\subset R_{sc}\setminus X$ to be an arbitrary subset of size $z-|X|$ and claim that the conditions \eqref{property:size right}, \eqref{property:small X} and \eqref{property:flexibility} of the proposition are all satisfied. Indeed it remains only to prove \eqref{property:flexibility} and this follows similarly to above, by taking sub-diamond trees of $\cD_{sc}$.  In detail, fix some $Y'\subset Y$ and let $R=Y'\cup X$. Then if $\rho_{sc}:V(T_{sc})\rightarrow R_{sc}$  and $\sigma_{sc}:E(T_{sc})\rightarrow \Sigma_{sc}$ are the defining bijective maps for $\cD_{sc}$, we have that the set of vertices $\{\rho_{sc}^{-1}(v):v\in R\}$ spans a subtree $T\subset T_{sc}$. Indeed we simply deleted  leaves from $T_{sc}$, namely $\rho_{sc}^{-1}(x)$ for $x\in R_{sc} \setminus Y'$. Taking $\Sigma=\{\sigma_{sc}(e):e\in E(T)\}$, we have that $\cD=(T,R,\Sigma)$ is the desired diamond tree.  \end{proof}

\section{Cascading absorption through orchards} \label{sec:orchards}

In this section we discuss orchards in our $(p,\beta)$-bijumbled graphs. We begin in Section~\ref{sec:absorbing orchards} by proving Lemma~\ref{lem:absorptionbetweenlayers} which details conditions for when one orchard absorbs another. In Section~\ref{sec:shrinkable orchards}, we then discuss the existence of shrinkable orchards, addressing Proposition~\ref{prop:shrinkable orchards} which tells us that we can find shrinkable orchards of all desired orders in the graphs we are interested in. The proof of Proposition~\ref{prop:shrinkable orchards} requires many ideas and two distinct approaches. Therefore, we defer the majority of the work to later sections and simply reduce the proposition here, splitting it into two `subpropositions' which will be tackled separately. Recall that Lemma~\ref{lem:absorptionbetweenlayers} and Proposition~\ref{prop:shrinkable orchards} were the two ingredients we needed to prove the cascading absorption through constantly many orchards in the proof of Theorem~\ref{thm:main}.

\subsection{Absorbing Orchards} \label{sec:absorbing orchards}
Recall the definition (Definition~\ref{def:orchard}) of an orchard and that
we say a~$(K,M)_r$-orchard $\bc{O}$ \emph{absorbs} a $(k,m)_r$-orchard $\bc{R}$ if there is a $((r-1) k,M)_r$-suborchard $\bc{O}'\subset \bc{O}$, such that there is a $K_r$-factor in $G[V(\bc{R})\cup V(\bc{O}')]$.  

In this section we prove Lemma \ref{lem:absorptionbetweenlayers}, restated below for convenience, which is a generalisation of \cite[Lemma 3.5]{Nen18}.  The lemma gives some sufficient conditions for an orchard to be able to absorb another orchard. 

\theoremstyle{plain}
\newtheorem*{lem:absorptionbetweenlayers}{Lemma \ref{lem:absorptionbetweenlayers}}
\begin{lem:absorptionbetweenlayers}{(Restated)}
For any $3\leq r\in \NN$ and $0<\zeta,\eta<1$ there exists an $\eps>0$ such that the following holds for any $n$ vertex $(p,\beta)$-bijumbled graph $G$ with $\beta\leq \eps p^{r-1}n$. 
Let $\bc{O}$ be a $(K,M)_r$-orchard  in $G$ such that $KM\geq\zeta n$. Then there exists a set $B \subset V(G)$ such that $|B|\leq \eta p^{2r-4}n$ and  $\bc{O}$ absorbs any $(k,m)_r$-orchard  $\bc{R}$ in $G$ with
\begin{equation} \label{eq:absorptionconstraints2}
  V(\bc{R})\cap (B \cup V(\bc{O}))=\emptyset, \qquad   k\leq K/(8r) \qquad \mbox{ and } \qquad
  kM\leq mK. \end{equation}
\end{lem:absorptionbetweenlayers}

Our proof scheme follows that of \cite{HKMP18b} which gives a polynomial time two phase algorithm for finding the necessary $K_r$-factor. The algorithm is a simple greedy algorithm and works by absorbing each diamond tree $\cB$ in the small orchard $\bc{R}$, one at a time. In more detail, for each diamond tree $\cB$ in $\bc{R}$, we find $r-1$ diamond trees $\cD_1,\ldots,\cD_{r-1}\in \bc{O}$ such that there is a copy of $K_r$ traversing the sets of removable vertices of $\cB$ and the diamond trees $\cD_1,\ldots,\cD_{r-1}$. This implies that there is a $K_r$-factor in $G[V(\cB)\cup V(\cD_1)\cup \ldots \cup V(\cD_{r-1})]$ (see Observation~\ref{obs:diamondtree}) and so we can add the $\cD_i$ to the suborchard $\bc{O}'$, forbid them from being used again, and move to the next diamond tree $\cB'\in \bc{R}$. Note that typically, we expect to succeed with this process. Indeed, the set of removable vertices of diamond trees in $\bc{O}$ is linear in size (and remains linear even after forbidding diamond trees $\cD\in \bc{O}$ used for previous $\cB\in \bc{R}$) and so a typical vertex has  $\Omega(pn)$ neighbours among this set of removable vertices. Hence, appealing to Corollary~\ref{cor:transversalcliques}~\eqref{cor:forbiddingsubgraph}~\ref{cor:r-1cliqueinpn} which states  that sets of size $\Omega(pn)$ host copies of $K_{r-1}$, we can expect to find a copy of $K_{r-1}$ in the neighbourhood of a typical removable vertex of $\cB\in \bc{R}$ which lies on the removable vertices of diamond trees in $\bc{O}$. As long as this copy of $K_{r-1}$ traverses sets of removable vertices of distinct diamond trees in $\bc{O}$, we will succeed. With a few extra ideas and a bit of preprocessing (for example partitioning $\bc{O}$ into $r-1$ suborchards at the start), this intuition holds true and we can successfully greedily start to build $\bc{O}'$. In fact, if $kM$ is small compared to $pn$, we can fully form $\bc{O}'$ in this way and no second phase is necessary. However, if $kM$ is large compared to $pn$ we may run into trouble as with this greedy approach, it may be the case that the neighbourhood of a removable vertex $v$ of a diamond tree $\cB\in \bc{R}$ has too small a size by the time we come to considering $\cB$. Indeed, as we run this greedy process, we forbid the diamond trees (and their removable vertices) which we add to $\bc{O}'$, from being used again. This could result in  $v$ having much fewer than $pn$ neighbours in the removable vertices of diamond trees in (the remainder of)~$\bc{O}$ and so we have no guarantee of finding a copy of $K_{r-1}$ in this neighbourhood. We resolve this issue by running a two-phase algorithm and reserving half of~$\bc{O}$ for the second phase. The key point is that if a diamond tree $\cB$ fails in the first round then it must be the case that \emph{all} of the removable vertices of $\cB$ have small neighbourhoods amongst the removable vertices of diamond trees in $\bc{O}$. Given that throughout the process,  many diamond trees in~$\bc{O}$ will remain available to use, pseudorandomness (more precisely, Corollary~\ref{cor:transversalcliques}~\eqref{cor:forbiddingsubgraph}~\ref{cor:1small}) tells us that the number of vertices that do not have large enough neighbourhoods, is relatively small. Hence, as each diamond tree $\cB\in \bc{R}$ which failed in the first phase, has a set of removable vertices which are atypical in this way, we can upper bound the number of diamond trees in $\bc{R}$ that fail in the first round. This upper bound will then be used to show that in the second phase, we are successful with each diamond tree, as throughout the second round, the number of removable vertices being forbidden (due to being used to absorb other diamond trees in $\bc{R}$) will be negligible and so the neighbourhoods of vertices amongst the removable vertices of diamond trees in the half of $\bc{O}$ reserved for this second phase, will remain large. 

\begin{proof}[Proof of Lemma \ref{lem:absorptionbetweenlayers}]
We fix $\alpha,\eta'<\frac{\eta\zeta ^2}{2^{3r}r^2}$ and choose $\eps>0$ small enough to apply Lemma \ref{lem:bad degrees} with $\eta_{\ref{lem:bad degrees}}=\eta'$ and Corollary \ref{cor:transversalcliques}  with  $\alpha_{\ref{cor:transversalcliques}}=\alpha$. Let $\bc{O}=\{\cD_1,\ldots,\cD_{K}\}$  be the $(K,M)_r$-orchard with each $\cD_i=(T_i,R_i,\Sigma_i)$ being a $K_r$-diamond tree of order between $M$ and $2M$. We start by arbitrarily partitioning $\bc{O}$ into $2(r-1)$ suborchards of size as equal as possible so that $\bc{O}=\bigcup_{j=1}^{2(r-1)}\bc{O}_j$ and each $\bc{O}_j$ is a $(K_j,M)_r$-orchard with $K_j=\frac{K}{2(r-1)}\pm 1\geq \frac{K}{2r}$. For $j\in[2(r-1)]$, we let \[Y_j:=\bigcup_{i:\cD_i\in \bc{O}_j}R_i\]
be the set of removable vertices of the diamond trees which feature in the $j^{th}$ suborchard. Note that $|Y_j|\geq K_jM\geq KM/(2r) \geq \zeta n/(2r)$ for each $j\in[2(r-1)]$. We define $B$ to be the set of vertices $v\in V\setminus V(\bc{O})$ such that for some $j\in [2(r-1)]$, $\deg_{Y_j}(v)< p |Y_j|/2$. By Lemma \ref{lem:bad degrees}  \ref{lem:badverticessmalldegs}, we have that 
\[|B|<\frac{2(r-1) \eta' p^{2r-4}n^2}{\min_j|Y_j|}\leq \eta p^{2r-4}n,\]
due to our lower bound on the size of the $|Y_j|$ and our upper bound on $\eta'$. 

Now as in the statement of the lemma, consider a $(k,m)_r$-orchard $\bc{R}=\{\cB_1,\ldots,\cB_{k}\}$ of diamond trees whose vertices lie in $V\setminus (B\cup V(\bc{O}))$. For $i'\in[k]$, let $Q_{i'}$ be the set of removable vertices of the diamond tree $\cB_{i'}$. We will show that for each $i'\in[k]$, there exists distinct indices $i_1=i_1(i'),\ldots, i_{r-1}=i_{r-1}(i')\in [K]$ such that there is copy of  $K_r$ which traverses the sets $Q_{i'}$ and $R_{i_1},\ldots, R_{i_{r-1}}$, where $R_{i_1}$ is the set of removable vertices of $\cD_{i_1}$ and likewise for $i_2,\ldots,i_{r-1}$.  Now, from Observation \ref{obs:diamondtree}, we have that for such an $r$-tuple $\cB_{i'}$, $\cD_{i_1},\ldots,\cD_{i_{r-1}}$, there is a  $K_r$-factor in $G[V(\cB_{i'})\cup V(\cD_{i_1}) \cup\ldots\cup V(\cD_{i_{r-1}})]$. We will prove that one can choose such indices $i_1,\ldots,i_{r-1}$ for each $i'\in[k]$ in such a way that no $i\in[K]$ is chosen more than once. That is, for $i'\neq j'\in[k]$, the sets $\{i_1(i'),\ldots,i_{r-1}(i')\}$ and $\{i_1(j'),\ldots,i_{r-1}(j')\}$ are  disjoint. Therefore our  suborchard $\bc{O}'\subset \bc{O}$ can simply be defined to be  the union of all the choices of $\cD_{i_j(i')}$ for $i'\in[k]$ and $j\in[r-1]$. 

\vspace{2mm}

We now show how to find the indices $i_1(i'),\ldots,i_{r-1}(i')$ for each $i'\in[k]$. We will achieve this via the following simple algorithm. We initiate the first round of the algorithm with $\bc{O}'=\emptyset$, $I=[k]$, $\bc{P}_j=\bc{O}_j$  and $Z_j=Y_j$ for $1\leq j\leq r-1$.   
Note that the $\bc{O}_j$ for $r\leq j \leq 2(r-1)$ do not feature in these definitions. This is because we will not use any diamond trees that lie $\cup_{j=r}^{2(r-1)}\bc{O}_j$ in this first round. Now the algorithm runs as follows.  For $i'=1,\ldots, k,$ we check if there exists some set $\{\cD_{i_j}\in \bc{P}_j:j\in[r-1]\}$ 
such that there is a $K_r$ traversing $Q_{i'}$ and the sets of removable vertices $R_{i_1}, \ldots, R_{i_{r-1}}$. If this is the case then we delete $\cD_{i_j}$ from $\bc{P}_j$  and add it to $\bc{O}'$ for $j\in[r-1]$ and we also delete $R_{i_j}$ from $Z_{j}$ for all $j\in[r-1]$. 
Furthermore, we delete $i'$ from $I$ and move to the next index $i'+1$ (or finish this round if $i'=k$). If it is not the case that such diamond trees exist in the  orchards $\bc{P}_j$,  
then we simply leave $i'$ as a member of $I$ and move on to the next index.

At the end of the first round, we have some set $I$ of indices remaining. We define $t:=|I|$ at this point. We will now use diamond trees in the orchards $\bc{O}_j$ with $r\leq j \leq 2(r-1)$  
to absorb these remaining diamond trees $\cB_{i'}$ with $i'\in I$. Thus we reset the process, setting $\bc{P}_j=\bc{O}_{j+r-1}$ and $Z_j=Y_{j+r-1}$ for all $j\in[r-1]$. 
We then follow the same simple process in the second round as we did in the first, running through the (remaining) $i'\in I$ in order and trying to find an appropriate set of diamond trees $\{\cD_{i_j}\in \bc{P}_j:j\in[r-1]\}$ at each step.  
We claim that in this second round, we can find such a set for every $i'\in I$ and so by the end of the second round, the set $I$ is empty and $\bc{O}'$ is such that $G[V(\bc{R})\cup V(\bc{O}')]$ hosts a $K_r$-factor. 

\vspace{2mm}

In order to prove this, our analysis splits into two cases. First consider when 
$kM< \frac{\zeta pn }{16r}$. In this case,  the second round is not even necessary as all indices succeeded in the first round. Indeed, note that every time we are succesful for an index $i'$, we delete at most $2M$ vertices from each of the $Z_j$. Therefore, at any instance in the first round of the process,  we have that any vertex $v$ which is not in $B$ has \[\deg_{Z_j}(v)\geq \frac{p|Y_{j}|}{2}-2kM \geq \frac{\zeta pn }{4r} -2kM \geq \frac{\zeta pn }{8r} ,\]
for all $j\in[r-1]$, using our lower bound on the  $|Y_{j}|$ and our upper bound on $kM$. 
But then, by Corollary \ref{cor:transversalcliques} \eqref{cor:forbiddingsubgraph} \ref{cor:r-1cliqueinpn} (applied in this instance with~$\tilde{G}$ being the empty graph and~$G'=G$), 
there exists a copy of $K_{r-1}$ traversing the sets $N_{Z_j}(v)$ for $1\leq j \leq r-1$.  
 When $v$ is any vertex in the removable set of vertices $Q_{i'}$ for some diamond tree $\cB_{i'}$ in the process, this  gives a copy of $K_{r}$ traversing $Q_{i'}$ and some sets of removable vertices $R_{i_j}$ for diamond trees $\cD_{i_j}\in \bc{P}_j$, $j\in[r-1]$, as desired.  In this way, we see that the process succeeds in every step of the first round to find a suitable $\{i_j(i'):j\in[r]\}$ for each $i'\in  [k]$ and $I$ is empty (i.e. $t=0$) at the end of the round. Note that we used here that the vertices of $Q_{i'}$ are not in $B$.  

\vspace{2mm}

When $ \frac{\zeta pn }{16r}\leq kM\leq mK $, the second round may be neccesary and we start with estimating $t$, the size of $I$ after the first round. 
Now note that at the end of the first round, \emph{before} we reassign the sets $Z_j$ to removable vertices in diamond trees in $\bc{O}_{j+r-1}$ for $j\in[r-1]$, if we take $Q=\bigcup_{i'\in I}Q_{i'}$, we have that there is no $K_r$ traversing $Q$ and the sets  $Z_1,\ldots,Z_{r-1}$. Indeed, otherwise there would be an $i'\in I$ and a vertex $v\in Q_{i'}\subseteq Q$ which is contained in a $K_r$ with a set of vertices $\{v_{i_j}\in Z_j:j\in[r-1]\}$.  
This contradicts that the index $i'$ failed to find a suitable set of $i_j$ in the first round. Thus, at the end of the first round, there is no $K_r$ traversing $Q$, and the $Z_j$, $j\in[r-1]$. Moreover 
we have that \[|Z_j|\geq \frac{KM}{2r}-2kM \geq \frac{KM}{4r} \geq \frac{\zeta n}{4r},\]
using the upper bound on $k$ from \eqref{eq:absorptionconstraints2} and the fact that at most $2M$ vertices are deleted from $Z_j$ every time we are successful with an index $i'\in I$. Thus, 
we can conclude from Corollary \ref{cor:transversalcliques} \eqref{cor:forbiddingsubgraph} \ref{cor:1small} 
that at the end of the first round,  $tm<|Q|<\alpha p^{2r-4}n$. Therefore 
\[t< \frac{\alpha p^{2r-4}n }{m}\leq \frac{\alpha p^2n}{m}\leq \frac{16\alpha r p K}{\zeta}\leq \frac{16\alpha r p n}{\zeta M} \leq \frac{\zeta p n }{16rM}, \] where we used here our lower and upper bounds on $kM$ to give an upper bound on $pn/m$ in the third inequality, the fact that $KM\leq n$ in the fourth inequality and our upper bound on $\alpha$ in the final inequality. 

We now turn to analyse the second round. Using our upper bound on $t$, we can upper bound the number of vertices deleted in each $Z_j$ throughout the second round, and using this we have that for any vertex $v$ not in $B$, any $j\in[r-1]$ and at any point in the second round, 
\[\deg_{Z_j}(v)\geq \frac{p|Y_{j+r-1}|}{2}-2tM\geq 
\frac{\zeta pn}{4r}-\frac{\zeta pn}{8r}\geq \frac{\zeta pn}{8r}.\]
Thus we can repeat the argument used for the case when $kM$ was small, seeing that at every step in the second round we are successful in finding an appropriate set of $i_j$ for $j\in[r-1]$ for each $i'\in I$. This completes the proof.  
\end{proof}

\subsection{Shrinkable Orchards} \label{sec:shrinkable orchards}
Here we are concerned with the existence of shrinkable orchards in pseudorandom graphs and verifying Proposition \ref{prop:shrinkable orchards}, which we restate below for the convenience of the reader. We also encourage the reader to remind themselves of Definitions \ref{def:trianglehypergraph} and \ref{def:shrinkable} as well as Observation \ref{obs:matchings in triangle hypergraph}.

\theoremstyle{plain}
\newtheorem*{prop:shrinkable orchards}{Proposition \ref{prop:shrinkable orchards}}
\begin{prop:shrinkable orchards}{(Restated)}
For any $3\leq r\in \NN$ and $0<\alpha,\gamma<1/2^{12r}$  there exists an $\eps>0$ such that the following holds for any $n$ vertex $(p,\beta)$-bijumbled graph $G$ with $\beta\leq \eps p^{r-1}n$ 
and any vertex subset $U\subseteq V(G)$ with $|U|\geq n/2$. For any $m\in \mathbb{N}$ with $1\leq m \leq n^{7/8}$  there exists a $\gamma$-shrinkable $(k,m)_r$-orchard $\bc{O}$ in $G[U]$ with $k\in \NN$ such that $\alpha n\leq km \leq  2\alpha n$. 
\end{prop:shrinkable orchards}

In order to prove Proposition \ref{prop:shrinkable orchards}, we will appeal to the methods of Sections \ref{sec:fractional matchings lin programming} and \ref{sec:almost perfect matchings}. We will  use Theorem \ref{thm:almost perfect matchings via greedy pfms} to reduce the problem to establishing the existence of perfect fractional matchings in the appropriate $K_r$-hypergraphs and 
we will then employ Lemmas \ref{lem:vertex fans} and \ref{lem:two stage fans} to find these perfect fractional matchings. In order that our hypergraph has the desired properties to apply these lemmas we need to choose the diamond trees which define our orchard carefully.  

It turns out that different arguments are needed for finding shrinkable orchards of different orders. In Section \ref{sec:small order shrinkable} we show how to find shrinkable orchards of small order, establishing the following intermediate proposition. 

\begin{prop}\label{prop:small shrinkable}
For any $3\leq r\in \NN$ and $0<\alpha,\gamma< 1/2^{3r}$  there exists an $\eps>0$ such that the following holds for any $n$ vertex $(p,\beta)$-bijumbled graph $G$ with $\beta\leq \eps p^{r-1}n$  and any vertex subset $U\subseteq V(G)$ with $|U|\geq n/2$. For any $m\in \mathbb{N}$ with \[1\leq m \leq \min\{ p^{r-2}n^{1-2r^3\gamma},n^{7/8}\},\]  there exists a $\gamma$-shrinkable $(k,m)_r$-orchard $\bc{O}$ in $G[U]$ with $k\in \NN$ such that $\alpha n\le km \le 2\alpha n$. 
\end{prop}

In Section  \ref{sec:large order shrinkable} we then address shrinkable orchards with large order, which results in the following.

\begin{prop}\label{prop:large shrinkable}
For any $3\leq r\in \NN$ and $0<\alpha,\gamma< 1/2^{12r}$  there exists an $\eps>0$ such that the following holds for any $n$ vertex $(p,\beta)$-bijumbled graph $G$ with $\beta\leq \eps p^{r-1}n$  and any vertex subset $U\subseteq V(G)$ with $|U|\geq n/2$. For any $m\in \mathbb{N}$ with \[ p^{r-1}n\leq m \leq n^{7/8},\]  there exists a $\gamma$-shrinkable $(k,m)_r$-orchard $\bc{O}$ in $G[U]$ with $k\in \NN $ such that $\alpha n\leq km\leq 2\alpha n$. 
\end{prop}

The proof of Proposition~\ref{prop:shrinkable orchards} is basically immediate from Propositions \ref{prop:small shrinkable} and \ref{prop:large shrinkable} but we spell it out nonetheless. 

\begin{proof}[Proof of Proposition \ref{prop:shrinkable orchards}]
We split into a case analysis based on the density $p$ of our graph $G$. First consider when $p\geq n^{-1/(10r)}$. Then we claim that $p^{r-2}n^{1-2r^3\gamma}\geq n^{7/8}$ and so the desired $\gamma$-shrinkable orchard of all orders up to $n^{7/8}$ can derived from Proposition~\ref{prop:small shrinkable}. Indeed we have that $p^{r-2}n^{1-2r^3\gamma}\geq n^{1-\frac{r-2}{10r}-2r^3\gamma}$ and \[1-\frac{r-2}{10r}-2r^3\gamma> 1-\frac{1}{10}-\frac{1}{40}=7/8,\]
due to our upper bound on $\gamma$ (and lower bound on $r$).

When $p<n^{-1/(10r)}$, we have that $p\leq n^{-2r^3\gamma}$ again due to our upper bound on $\gamma$. Hence we can apply Proposition~\ref{prop:small shrinkable} to find $\gamma$-shrinkable orchards of orders $m\leq n^{7/8}$ such that  $m<p^{r-1}n\leq p^{r-2}n^{1-2r^3\gamma}$ and apply Proposition~\ref{prop:large shrinkable} to find $\gamma$-shrinkable orchards with orders $m$ such that $p^{r-1}n\leq m \leq n^{7/8}$. This settles all cases and so gives the proposition. 
\end{proof}

In both cases, a simpler argument works for the extreme cases, that is, when the order is small in Proposition~\ref{prop:small shrinkable} or when the order is large in Proposition~\ref{prop:large shrinkable}. Extra ideas are then needed to push the approaches, extending the ranges of the two propositions so that they meet and cover all desired orders. In more detail, an easier form of Proposition~\ref{prop:small shrinkable}  can cover orders which get close to $p^{r-1}n$ (see Proposition  ~\ref{prop:verysmallorchards}). Again the separation required depends on $\gamma$, explicitly $m\leq p^{r-1}n^{1-r^3\gamma}$.  This is already enough to cover all desired orchard orders when $p$ is large. On the other hand, a basic form of the argument for large order orchards gives shrinkable orchards of order at least $p^{r-1}n$ when $p$ is large and of order at least $p^{1-r}$ when $p$ is smaller (see Proposition~\ref{prop:easy large order shrinkable}). Interestingly,  Fact~\ref{fact:dense} implies exactly that $p^{1-r}=\Omega(p^{r-2}n)$ always and so proves that when $p$ is small (close to the lower bound of $\Omega(n^{-1/(2r-3)})$) and our bijumbled graph is sparse, both  the simpler arguments for small orders and large orders \emph{as well as} their extensions are needed. Indeed using the simpler version, Proposition~\ref{prop:verysmallorchards}, for small orders and the full power of Proposition~\ref{prop:large shrinkable} leaves a small gap in the orders and so does using Proposition~\ref{prop:small shrinkable} in conjunction with the easier Proposition~\ref{prop:easy large order shrinkable}. In order to help the reader through the next two sections, in both cases we begin by presenting the easier weaker versions of the statements we need. This then lays the foundation for the full proofs and allows us to discuss the more technical aspects needed to push the ranges for which we can prove the existence of shrinkable orchards.

\section{Shrinkable orchards of small order} \label{sec:small order shrinkable}

Our first argument for proving the existence of shrinkable orchards works provided the order of the orchard is not too large, establishing Proposition~\ref{prop:small shrinkable}. Before embarking on this we have to go through several steps. Firstly, in Section~\ref{sec:systems}, we generalise the theory of shrinkable orchards built up in Section~\ref{sec:proofmain}, allowing slightly more flexibility for our consequent proofs. In Section~\ref{sec:smallconditionsforshrinkability}, we then use the theory of perfect fractional matchings to give conditions that guarantee an orchard is shrinkable. In Section~\ref{sec:smallest shrinkable}, we  show how this immediately implies the existence of shrinkable orchards of  small order. However this falls short of Proposition~\ref{prop:small shrinkable} and in the rest of this section we push the ideas to extend the range of orders we can cover, showing how to cleverly choose diamond trees of our orchard in Section~\ref{sec:controlling degrees} which allows us to prove the full Proposition~\ref{prop:small shrinkable}  in Section~\ref{sec:slightly larger}.

\subsection{From orchards to systems} 
\label{sec:systems}
We begin by generalising our definitions slightly, allowing us to work not just with orchards but also set systems. 

\begin{dfn} \label{def:setsystems}
 Given a graph $G$ we say a set of pairwise disjoint subsets $\Lambda\subset 2^{V(G)}$ is a \emph{$(k,m)$-system} if $m\leq|Q|\leq 2m$ for each $Q\in \Lambda$ and $|\Lambda|=k$. That is, a $(k,m)$-system is just a set family of $k$ disjoint vertex sets of size between $m$ and $2m$. 
 
 Now given
 a $(k,m)$-system $\Lambda$ in a graph $G$, the \emph{$K_r$-hypergraph} generated by $\Lambda$, denoted $H=H(\Lambda;r)$ is the $r$-uniform hypergraph with vertex set $V(H)=\Lambda$ and with $\{Q_{i_1},\ldots,Q_{i_r}\}\in \binom{\Lambda}{r}$ forming a hyperedge in $H$ if and only if there is a copy of $K_r$ traversing the sets $Q_{i_1},\ldots,Q_{i_r}$ in $G$.
 
 Finally for  $0<\gamma<1$, we say a $(k,m)$-system $\Lambda$ in a graph $G$ is \emph{$\gamma$-shrinkable} (with respect to $r$) if there exists a subsystem $\Gamma \subset \Lambda$ of size at least $\gamma k$ such for \emph{any} subsystem $\Gamma'\subseteq \Gamma$, we have that there is a matching in $H:=H(\Lambda\setminus \Gamma';r)$ covering all but $k^{1-\gamma}$ of the vertices of $H$.
 \end{dfn}
 
 Note that given a $(k,m)_r$-orchard $\bc{O}$ we can define a $(k,m)$-system $\Lambda$ as the sets of removable vertices of diamond trees in $\bc{O}$. That is,  $\Lambda:=\{R_\cD:\cD\in\bc{O}\}$. Then the $K_r$-hypergraphs generated by $\bc{O}$ and $\Lambda$ coincide i.e. $H(\Lambda;r)=H(\bc{O})$, 
 and $\bc{O}$ is $\gamma$-shrinkable if and only if $\Lambda$ is  $\gamma$-shrinkable. However Definition~\ref{def:setsystems} allows us slightly more flexibility, giving us the ability to focus on subsets of removable vertices. The next observation highlights this and although the result is trivial, it will be important for our proofs.
 
 \begin{obs} \label{obs:setsystemfromorchard}
 Suppose $r\geq 3$, $0<\gamma<1$ and  $\bc{O}=\{\cD_1,\ldots,\cD_k\}$ is a $(k,m)_r$ orchard in a graph $G$ with $R_i$ being the set of removable vertices of $\cD_i$ for $i\in[k]$. Then if $\Lambda=\{Q_1,\ldots,Q_k\}$ is some $(k,m')$-system (for some $m'$) such that $Q_i\subseteq R_i$ for $i\in[k]$ and $\Lambda$ is $\gamma$-shrinkable (with respect to $r$), then $\bc{O}$ is also $\gamma$-shrinkable.  
 \end{obs}

It will become clear why such a relaxation is useful for us and thus why we make this switch to working with set systems. 

\subsection{Sufficient conditions for shrinkability} \label{sec:smallconditionsforshrinkability}

We now explore the conditions on set systems which guarantee shrinkability. We begin by giving some local conditions on a set system which guarantee that it is shrinkable given that it lies in the pseudorandom graphs we are interested in. 

\begin{lem} \label{lem:local to shrink}
For any $3\leq r\in \NN$ and $0<\alpha,\gamma< 1/(2^{r}r^2)$  there exists an $\eps>0$ such that the following holds for any $n$ vertex $(p,\beta)$-bijumbled graph $G$ with $\beta\leq \eps p^{r-1}n$. Suppose $\Lambda\subset 2^{V(G)}$ is a $(k,m)$-system such that $m\leq n^{7/8}$,  $km\geq \alpha n$, $pk\ge n^\gamma$ and:
\begin{enumerate}
 \item \label{cond:local min degree} There exists a subsystem $\Gamma\subset \Lambda$ such that $|\Gamma|\geq \gamma k$ and the following holds with  $Y:=\cup\{P: P\in \Lambda\setminus \Gamma\}$.  For every $Q\in \Lambda$, there exists a vertex $v\in Q$ such that 
    \[\deg^G_{Y}(v)\geq \alpha pkm;\]
\item \label{cond:local upper bound}For any $u\in \cup_{P\in \Lambda} P$ and $Q\in \Lambda$, we have that $\deg^G_Q(u)\leq p^{r-1}n^{1-r^3\gamma}$.

\end{enumerate}
Then $\Lambda$ is $\gamma$-shrinkable with respect to $r$.
\end{lem}

Let us make a few remarks before proving the lemma. Firstly, note that condition \eqref{cond:local min degree}, despite the slight technicality necessary to avoid dependence on sets in $\Gamma$, is a natural condition. Indeed, we are requiring that at least one vertex in each set is well connected to the other sets and has a constant fraction of the degree that we would expect on average. Condition \eqref{cond:local upper bound} is perhaps more mysterious as it is unclear why having an \emph{upper} bound on the degree of a vertex  to another set in the system is advantageous. The point is that this guarantees that each of the vertices has a neighbourhood that is well-spread across the other sets of the set system, without being too concentrated on any one other set. Within the proof this necessity manifests itself as we appeal to Theorem~\ref{thm:almost perfect matchings via greedy pfms} and so will need that when we disallow edges between certain pairs of sets  from being used (dictated by the graph $J$), we do not significantly alter the graph in which we work. The details follow in the proof. 

\begin{proof}[Proof of Lemma~\ref{lem:local to shrink}]

Fix $\eps>0$ small enough to apply Corollary~\ref{cor:transversalcliques} with $\alpha_{\ref{cor:transversalcliques}}=\alpha'<\alpha^3/2r$ and small enough to force $n$ to be sufficiently large in what follows. Fixing $\Gamma\subset \Lambda$ as in condition \eqref{cond:local min degree}, we have to show that for any $\Gamma'\subset \Gamma$, the $K_r$-hypergraph $H=H(\Lambda\setminus \Gamma';r)$ has a matching covering all but $k^{1-\gamma}$ vertices  of $H$. So fix such a $\Gamma'$, let $\Lambda^*:=\Lambda\setminus \Gamma'$ and let $H:=H(\Lambda^*;r)$. 

In order to show the existence of a large matching in $H$, we appeal to Theorem~\ref{thm:almost perfect matchings via greedy pfms}. So let us fix  $N=|V(H)|$ and note that as $N\geq (1-\gamma)k$ and $k\geq \alpha n^{1/8}$ due to our conditions on $k$ and $m$, we can assume that $N$ is sufficiently large in what follows. Now fix some $2$-uniform graph $J$ on $V(H)$ of maximum degree at most $N^{r^2\gamma}$.  If we can show that $H\setminus H_J$ contains  a perfect fractional matching, then we are done by Theorem~\ref{thm:almost perfect matchings via greedy pfms} as, because $J$ was arbitrary, the theorem guarantees a matching covering all but at most $N^{1-\gamma}\leq k^{1-\gamma}$ vertices of $H$. 

In order to study $H\setminus H_J$, we look at the forbidden edges of $G$ which $J$ imposes. That is, we define  \[\tilde{G}_J:= \bigcup_{\{Q_1,Q_2\}\in E(J)} G\left[Q_1,Q_2\right]\bigcup_{Q\in \Lambda^*} G[Q]\]
where $G\left[Q_1,Q_2\right]$ denotes the set of all edges in $G$ between the sets $Q_1$ and $Q_2$ and $G[Q]$ denotes all the edges induced by $G$ in the set $Q$. We then have that for any vertex $v\in V(G)$, $\deg^{\tilde{G}_J}(v)=0$ if $v\notin \cup_{P\in \Lambda^*}P$ and if $v\in Q\in \Lambda^*$, then
\begin{equation} \label{eq:tilde G max deg} \deg^{\tilde{G}_J}(v)\leq \sum_{P\in N^J(Q)\cup \{Q\}}\deg^G_P(u)\leq  (N^{r^2\gamma}+1)p^{r-1}n^{1-r^3\gamma}\leq p^{r-1}n^{1-\gamma},\end{equation}
using \eqref{cond:local upper bound}, the upper bound on the degrees in $J$ and the fact that $N\leq n$.

Now defining $G'_J:=G \setminus \tilde{G}_J$, we have that $H \setminus H_J$ is precisely the hypergraph obtained by viewing $\Lambda^*$ as a $(N,m)$-system in $G'_J$ and taking the $K_r$-hypergraph $H^*=H(\Lambda^*;r)$ in $G'_J$. Indeed, as there are no edges of $G'_J$ between two sets,  $Q_1$ and $Q_2$ say, which form an edge in $J$, there can be no edge of $H$ in $H^*$ which contains both $Q_1$ and $Q_2$. We therefore switch from now on to considering  $H^*$ as the $K_r$-hypergraph generated by $\Lambda^*$ in $G'_J$.

In order to prove the existence of a perfect fractional matching in $H^*$, we will appeal to Lemma~\ref{lem:two stage fans}, fixing $M=\alpha^2pk$. Note that due to our lower bound on $pk$, we certainly have that $M_1\ge r$.  We thus need to check that conditions \ref{cond:small vertex fan} and \ref{cond:expansion}  
of that lemma hold. For \ref{cond:small vertex fan}, fix some $Q\in \Lambda^*$. From \eqref{cond:local min degree} we have that there exists some vertex $v$ in $Q$ such that $\deg^G_{W}(v)\geq \alpha pkm$ where $W:=\cup_{P\in\Lambda^*}P$ and so taking $U:=N^{G_J'}_{W}(v)$ we have that \[|U| \geq \alpha pkm - \deg^{\tilde{G}_J}(v) \geq \alpha pkm/2 ,\]
using \eqref{eq:tilde G max deg}. Moreover, due to \eqref{cond:local upper bound}, we can spilt $U$ into disjoint sets $U_1,\ldots, U_{r-1}$ such that $|U_i|\geq \alpha pkm /(2r)$ for each $i$ and we have that for any $P\in \Lambda^*$, there exists an $i\in[r-1]$ such that $|P\cap U| \subset U_i$. 
That is,  we simply partition $U$ into $r-1$ roughly equal size parts such that vertices which lie in the same $P$ end up in the same part.   Condition~\eqref{cond:local upper bound} of the lemma guarantees that $U\cap P$ is small enough for each $P\in { \Lambda}^*$ and so we can do this partition in such a way that each of the $U_i$ are roughly equal in size. We will now repeatedly find $(r-1)$-cliques  in $G_J'$ traversing the $U_i$ and build a fan $F_Q$ of size $M$ in $H^*$ focused at $Q$. We start with $F_Q$ being empty  and each time we find a  copy $S=\{u_1,\ldots,u_{r-1}\}\in K_{r-1}(G'_J)$ in $G'$  with $u_i\in U_i$ for  $i\in[r-1]$   we have that there exists some sets $P_{1},\ldots, P_{r-1}\in \Lambda^*$ such that $u_i\in P_i$. 
We add the hyperedge between  $P_1, \ldots, P_{r-1}$ and $Q$ to the fan $F_Q$ and delete  any vertices in $P_i$ from $U_i$ for $i\in[r-1]$. We repeat this process and note that we are successful in every step until $F_Q$ has size $M$. Indeed this follows from  Corollary \ref{cor:transversalcliques} \eqref{cor:forbiddingsubgraph} \ref{cor:r-1cliqueinpn} as while $|F_Q|<M$, we have deleted at most\footnote{Here  we use that every set in $\Lambda$ has size at most $2m$. } $M2m\leq \alpha pkm/(4r)$ vertices from each $U_i$ and so $U_i\geq \alpha p km/(4r)\geq \alpha'pn$, using our upper bound on $\alpha'$ and our lower bound on $km$. 

We now turn to verifying  \ref{cond:expansion} 
of Lemma \ref{lem:two stage fans}.  We will show that given any $r$-tuple of disjoint subsystems $\Gamma_1,\Gamma_2, \ldots, \Gamma_{r}\subset \Lambda^*$ such that $|\Gamma_1|=M$ and $|\Gamma_i|\geq \alpha k/r$ for $2\leq i \leq r$, there exists a hyperedge of  $H^*$ with one endpoint in  each of the $\Gamma_i$. Indeed, this follows from Corollary \ref{cor:transversalcliques} \eqref{cor:forbiddingsubgraph} \ref{cor:1small} 
as taking $U_i:= \bigcup_{P\in \Gamma_i}P$ for $i\in[r]$, we have that there exists an $r$-clique $S\in K_{r}(G'_J)$ traversing the $U_i$ which in turn gives the hyperedge.  
The condition \ref{cond:expansion} then clearly follows as any subsystem of $\Gamma^*$ of size $N/(2r)\geq  \alpha k$ can be split into $(r-1)$ subsystems of size at least $\alpha k/r$.  Note that in both applications of Corollary~\ref{cor:transversalcliques} \eqref{cor:forbiddingsubgraph} above  we used \eqref{eq:tilde G max deg} to show that we could find cliques that avoid using edges of  $\tilde{G}_J$.  
 The lemma  now follows from Lemma \ref{lem:two stage fans}.  
\end{proof}

As previously noted, condition \eqref{cond:local min degree} of Lemma~\ref{lem:local to shrink} is somewhat weak and just requires that each set in the set system  contains a vertex that acts typically. We now show how we can `clean up' a set system; losing sets which do not have a typical vertex in order to recover condition \eqref{cond:local min degree}. This allows us to focus on finding systems which satisfy condition \eqref{cond:local upper bound} of Lemma \ref{lem:local to shrink}.

\begin{lem} \label{lem:recover min deg}
For any $3\leq r\in \NN$ and $0<\alpha,\gamma< 1/(2^{r}r^2)$  there exists an $\eps>0$ such that the following holds for any $n$ vertex $(p,\beta)$-bijumbled graph $G$ with $\beta\leq \eps p^{r-1}n$. Suppose $\Lambda'\subset 2^{V(G)}$ is a $((1+\gamma)k,m)$-system such that $m\leq n^{7/8}$,  $\alpha n\leq km\leq 2\alpha n$, $pk\geq n^{\gamma}$ and for any $u\in \cup_{P\in \Lambda'} P$ and $Q\in \Lambda'$, we have that $\deg^G_Q(u)\leq p^{r-1}n^{1-r^3\gamma}$. Then there exists a $(k,m)$-system  $\Lambda\subset \Lambda'$ 
which is $\gamma$-shrinkable with respect to $r$.
\end{lem}

\begin{proof}
We fix $\eps>0$ small enough to apply Lemma~\ref{lem:local to shrink} and to apply Lemma \ref{lem:bad degrees}  with $\eta<\gamma \alpha^2/2$. The method of the proof is simple; we aim to apply Lemma~\ref{lem:local to shrink} and so obtain $\Lambda$ from $\Lambda'$ by losing the sets which violate condition \eqref{cond:local min degree} of that lemma. By Lemma \ref{lem:bad degrees} \ref{lem:badverticessmalldegs}, there are few vertices which have small degree ($\leq p|Y|/2$) to any set  $Y$ which is large enough and so we can expect that we do not lose many sets when transitioning from $\Lambda'$ to $\Lambda$. One complication is that the definition of $Y$ in condition \eqref{cond:local min degree} of Lemma \ref{lem:local to shrink} \emph{depends} on the sets in the system and so we cannot guarantee that a set satisfying \eqref{cond:local min degree} continues to satisfy the condition once other sets have been removed. In order to handle this, we delete sets in the system one by one, creating a process which will terminate with a system which has the desired minimum degree condition. The details now follow.

 We begin by fixing  $ \Gamma_0 \subset \Lambda'$ to be some arbitrary subsystem of size $(1-\gamma)k$ and we initiate the process by setting $\Theta=\Lambda'$ and setting a `bin' system $\Phi$ which we initiate as being empty, that is,   we set $\Phi=\emptyset$.  Throughout the process 
we also define $W$ so that \[W:=   \bigcup_{Q\in \Gamma_0\cap \Theta} Q\] is the subset of vertices that lie  in (sets that belong to) the current system $\Theta$ as well as the system $\Gamma_0$. Now the process runs as follows. If there is a set $P$ in $\Theta$ such that    $\deg_{W}(v)< \alpha pkm$ for all $v\in P$,  then we delete $P$ from $\Theta$ and add it to $\Phi$. Hence if $P\in \Gamma_0$ then we also delete $P$ from $W$. 
We claim that this process terminates with $|\Phi| \leq  \gamma k$. Indeed if this were not the case  then consider the process at the point where $|\Phi|= \gamma k$. At this point we have that 
\[|W|= \left|\bigcup_{Q\in \Gamma_0\setminus \Phi}Q\right| = \left(|\Gamma_0|-|\Phi|\right)m= (k-2\gamma k)m\geq km/2\geq  \alpha n/2.\]
 Now  Lemma \ref{lem:bad degrees}  \ref{lem:badverticessmalldegs} implies that at most  \[\frac{\eta p^{2r-4}n^2}{|W|}\leq \frac{2\eta p^{2r-4}n}{\alpha}<\gamma \alpha n\leq \gamma km= |\Phi|m,\]
vertices can have degree less than $p|W|/2$ to $W$. This leads to a contradiction. 
Indeed, it follows from how $\Phi$ is defined that at this point in the process, $|\cup_{P\in \Phi}P|\geq |\Phi|m$ and  for all $v\in \cup_{P\in \Phi}P$, we have  $\deg_{W}(v) < \alpha pkm < p |W|/2$. Indeed, if a vertex $v\in P\in \Phi$ had a larger degree to $W$ then~$P$ would not have been added to $\Phi$ in the process.

Hence when the process terminates we have that $|\Phi|\leq \gamma k$ and we have  that $|\Theta|=|\Lambda'\setminus \Phi| \geq k$. We fix~$ \Lambda\subseteq  \Theta $ of size $k$ so that ${ \Gamma}_0 \cap { \Theta}\subseteq  \Lambda$. We also fix ${ \Gamma}\subseteq { \Lambda} \setminus ({ \Gamma}_0 \cap { \Theta}) $ of size $\gamma k$ (which is possible as $|{ \Gamma}_0 \cap { \Theta}|\le (1-\gamma)k$). 
We claim that~${ \Lambda}$ is~$\gamma$-shrinkable with respect  to~$r$. Indeed, this follows directly from Lemma~\ref{lem:local to shrink} noting that  condition~\eqref{cond:local min degree} is satisfied in~${ \Lambda}$ with respect to~${ \Gamma}$ due to how we constructed~${ \Lambda}$. 
\end{proof}

\subsection{The existence of shrinkable orchards of small order} \label{sec:smallest shrinkable}
We are now ready to prove the existence of shrinkable orchards by appealing to Lemma~\ref{lem:recover min deg}. Indeed, we simply need to find orchards which satisfy the maximum degree condition given there. This condition is immediate when the order of the orchards which we aim for is sufficiently small, leading to the following easy consequence. 

\begin{prop} \label{prop:verysmallorchards}
For any $3\leq r\in \NN$ and $0<\alpha,\gamma< 1/2^{3r}$  there exists an $\eps>0$ such that the following holds for any $n$ vertex $(p,\beta)$-bijumbled graph $G$ with $\beta\leq \eps p^{r-1}n$  and any vertex subset $U\subseteq V(G)$ with $|U|\geq n/2$. For any $m\in \mathbb{N}$ with \[1\leq m \leq \min\{p^{r-1}n^{1-r^3\gamma},n^{7/8}\},\]  there exists a $\gamma$-shrinkable $(k,m)_r$-orchard $\bc{O}$ in $G[U]$ with $k$ such that $\alpha n\leq km\leq 2\alpha n$. 
\end{prop}

\begin{proof}
We fix $\eps>0$ small enough to apply Proposition~\ref{prop: choosing removable vertices} and  Lemma~\ref{lem:recover min deg} with $\alpha,\gamma$ as defined here and $k\in \mathbb{N}$ such that $\alpha n\le km \le 2\alpha n$. Note that due to our upper bound on $m$, we certainly have that $pk\ge n^\gamma$.  We begin by finding a $((1+\gamma)k,m)_r$-orchard $\bc{O}'$ in $G[U]$. This can be done by repeated applications of Proposition~\ref{prop: choosing removable vertices}. Indeed we initiate a process by fixing $U'=U$ and $\bc{O}'=\emptyset$ and at each step we find some  $K_r$-diamond tree $\cD$ of order $m$  in $U'$, add it to $\bc{O}'$ and delete its vertices from  $U'$. We claim that we can do this until $\bc{O}'$ has size $(1+\gamma)k$. Indeed this follows because at any point in the process, $|V(\bc{O}')|\leq (1+\gamma)kmr=(1+\gamma)2\alpha rn\leq n/4$ due to our upper bounds on $\alpha$ and $\gamma$. Therefore, throughout the process, we have that $|U'|\geq n/4$ and so it can be split into two disjoint sets of size at least $4\alpha r n$. Therefore applying Proposition~\ref{prop: choosing removable vertices} with $z=m$ (and taking $Y'=Y$ in \eqref{property:flexibility}) gives us the existence of the diamond tree at each step of this process. 

Now defining $\Lambda'=\{R_\cD:\cD\in\bc{O}'\}$, to be the $((1+\gamma)k,m)$-system generated by taking the sets of removable vertices of diamond trees that lie in $\bc{O}'$, we have that $\Lambda'$ satisfies the hypothesis of Lemma~\ref{lem:recover min deg} due to the fact that $m\leq p^{r-1}n^{1-r^3\gamma}$. Hence Lemma~\ref{lem:recover min deg} implies the existence of a subsystem $\Lambda\subset \Lambda'$ of size $k$ which is $\gamma$-shrinkable with respect to $r$. Finally taking $\bc{O}:=\{\cD\in \bc{O}':R_\cD\in \Lambda\}$, we have that $\bc{O}$ is the required $\gamma$-shrinkable $(k,m)_r$-orchard by Observation~\ref{obs:setsystemfromorchard}. 
\end{proof}

For dense graphs (that is, when $p$ is large), Proposition~\ref{prop:verysmallorchards} is already enough to establish Proposition~\ref{prop:shrinkable orchards}. On the other hand, for sparse graphs Proposition~\ref{prop:verysmallorchards} can only be used for orchards of very small order and becomes redundant as the order $m$ approaches $p^{r-1}n$. However, in deriving Proposition~\ref{prop:verysmallorchards}, we were quite naive in our application of Lemma~\ref{lem:recover min deg}, using the order of a diamond tree as an upper bound on the degrees of vertices to the removable set of vertices of the diamond tree. For a set $Q$ we expect a typical vertex $v\in V(G)$ to have $\deg_Q(v)\leq p|Q|$ and so we can hope that Lemma~\ref{lem:recover min deg} can be applied to imply the existence of shrinkable orchards whose orders approach $p^{r-2}n$, gaining an extra power of $p$ over Proposition~\ref{prop:verysmallorchards}. This is the content of the rest of this section. 

\subsection{Controlling degrees to removable sets of vertices} \label{sec:controlling degrees}

A reasonable approach when trying to  apply Lemma~\ref{lem:recover min deg} to deduce the existence of larger shrinkable orchards is to start with a larger (in size) orchard than we desire and crop diamond trees which fail the  bounded degree condition. This approach is reminiscent of how we derived Lemma~\ref{lem:recover min deg} from Lemma~\ref{lem:local to shrink}, where we greedily lost diamond trees which violated condition \eqref{cond:local min degree} of Lemma \ref{lem:local to shrink}. In this case though, our condition is harder to satisfy. Indeed, we require that \emph{all} vertices in a set in our system satisfy the  degree condition and not just a single vertex.  In order to achieve this, we will need to appeal to (the full power of) Proposition~\ref{prop: choosing removable vertices} to choose our diamond trees. 
As Proposition~\ref{prop: choosing removable vertices} does not give full control over the set of vertices which end up as the set of removable vertices, we have to settle with being able to conclude our desired upper bound on the degrees of vertices to a \emph{subset} of the removable vertices. The detailed statement is as follows.

\begin{lem}
\label{lem:diamond trees with low degrees}
For any $3\leq r\in \NN$ and $0<\gamma, \eta<1/r^2 $  there exists an $\eps>0$ such that the following holds for any $n$ vertex $(p,\beta)$-bijumbled graph $G$ with $\beta\leq \eps p^{r-1}n$  and any vertex subset $U\subseteq V(G)$ with $|U|\geq n/4$. For any $m\in \mathbb{N}$ with $p^{r-1}n^{1-\gamma}\leq m \leq p^{r-2}n^{1-2\gamma},$ 
 there exits a $K_r$-diamond tree $\cD=(T,R,\Sigma)$, of order at most $2m$ such that $V(\cD)\subset U$ and there exists a subset $Q\subseteq R$ of removable vertices such that $|Q|=m$ and  all but at most $\eta m$ vertices $v\in V(G)$ have $\deg_Q(v)\leq p^{r-1}n^{1-\gamma}$.  
\end{lem}

\begin{proof}
Fix $\eps>0$ small enough to apply Proposition~\ref{prop: choosing removable vertices} with $\alpha= 1/(2^{2r}r)$, small enough to apply Lemma~\ref{lem:bad degrees} with $\eta_{\ref{lem:bad degrees}}=\eta'<\alpha^2\eta/2^6$ and small enough to force $n$ to be sufficiently large in what follows. We begin by splitting $U$ into disjoint subsets $U'$ and $W'$ arbitrarily so that $|U'|,|W'|\geq n/8\geq 4 \alpha r n $, noting that this is possible due to our definition of $\alpha$. 

Now fix some $m$ with $p^{r-1}n^{1-\gamma}\leq m \leq p^{r-2}n^{1-2\gamma}$ and define $q:=p^{r-1}n^{1-\gamma}m^{-1}$.  Note that
$8p\leq p n^\gamma\leq q \leq 1$ due to our conditions on $m$. As we aim to find a set $Q$ of size $m$, the condition that $\deg_Q(v)> p^{r-1}n^{1-\gamma}$ is equivalent to having that $\deg_Q(v)>q|Q|$. As discussed above, given a diamond tree of the correct order and a subset $Q$ of $m$ removable vertices, 
we can appeal to Lemma \ref{lem:bad degrees}  \ref{lem:bad vertices high degs} 
to bound the number of vertices which have high degree to  $Q$. However, the bound is not strong enough for our purposes so we instead appeal to the full power of Proposition \ref{prop: choosing removable vertices}. The idea is to take $Y$ to be much bigger than $m$. Therefore applying Lemma \ref{lem:bad degrees}  \ref{lem:bad vertices high degs} with respect to $Y$ gives a much stronger upper bound on the number of vertices which have large degree (at least $q|Y|$, say) to $Y$. 
If we then take $Q$ to be a \emph{random subset} of $Y$ then we expect the density of the neighbourhood of a vertex in $Q$ to have roughly the same density as the neighbourhood of that vertex in $Y$. Hence, we can bound the number of vertices which have large degree to $Q$ by `carrying over' the bound on the number of vertices which had large degree to $Y$. 
The details follow. 

First we fix $d_*:=\alpha^2p^{r-1}n$ and $z:=\min\{\alpha n,d_* m/2\}$. Note that for $n$ large, due to Fact~\ref{fact:dense}, $d_*$ will also be large. Now apply Proposition \ref{prop: choosing removable vertices} to obtain disjoint subsets $X,Y\subset U'$ as in the statement of Proposition \ref{prop: choosing removable vertices}. Note that $|X|\leq 2z/d_* \leq m$ and  $|Y|\geq z-|X|\geq  z/2$ for $n$ sufficiently large. Fix a subset $Z\subset Y$ of size $z/2$ and let $B\subset V(G)$  be the set of vertices $v\in V(G)$ such that $\deg_{Z}(v)>q|Z|/4$. We claim that $|B|$ has size at most $\eta m$. 
Indeed, noting that $q/4\geq 2p$, 
Lemma \ref{lem:bad degrees} \ref{lem:bad vertices high degs} gives that  
\begin{equation} \label{eq:high deg vertex count}|B| \leq \frac{2^4\eta'p^{2r-2}n^2}{q^2|Z|}=\frac{2^5\eta'n^{2\gamma}m^2}{z}\leq \begin{cases}
    \eta  n^{2\gamma-1}p^{1-r}m  & \text{if }  z=d_* m/2, \\
    \eta n^{2\gamma-1}m^2   & \text{if } z=\alpha n. 
\end{cases} \end{equation}
In the case that $z=d_* m/2$, the estimate in \eqref{eq:high deg vertex count} is less than $\eta m$ for large $n$ due to the condition that $\gamma<1/r^2$ and the fact that $p\geq n^{-1/(2r-3)}$ (Fact \ref{fact:dense}). In the case that $z=\alpha n$, the estimate in \eqref{eq:high deg vertex count} is less than $\eta m$ due to the fact that $m\leq  p^{r-2} n^{1-2\gamma}\leq n^{1-2\gamma}$. 

For each $v\notin B$, we have that $\deg_{Z}(v)\leq q|Z|/4$ and so we let $N_v\subset Z$ be a subset of exactly $q|Z|/4$ vertices in $Z$ such that $N_v$ contains all the neighbours of $v$ which lie in $Z$. 
Now consider a random subset $Q_1\subset Z$ where we keep each vertex independently with probability $p'=4m/z$, noting that $0 \leq p' \leq 1$ for large enough $n$. Clearly $\EE \big[|Q_1|\big]=p'|Z|=2m$ and for each $v\in V(G)\setminus B$, have that $\EE\big[|Q_1\cap N_v|\big]=p'|N_v|=qm/2$. We get concentration for these random variables from Theorem \ref{thm:chernoff} which is strong enough to do a union bound and conclude that whp as $n$ (and hence $m$ and $qm$) tend to infinity, we have that $|Q_1|\geq m$ and $|Q_1\cap N_v|\leq qm$ for all $v\in V(G)\setminus B$. Therefore, for sufficiently large $n$, we can fix such an instance of $Q_1$ and take $Q$ to be a subset of $Q_1$ such that $|Q|=m$. Therefore for all vertices $v\in V(G)\setminus B$, we have that \[\deg_Q(v) \leq |Q\cap N_v| \leq |Q_1 \cap N_v|\leq qm= p^{r-1}n^{1-\gamma}.\] 
We have that $|B|\leq \eta m$ from above and we use the conclusion of Proposition \ref{prop: choosing removable vertices} to give a $K_r$-diamond tree $\cD=(T,R,\Sigma)$ with  removable vertices $R:=X \cup Q\subset U'\subset U$ and $\Sigma$ a matching of $(r-1)$-cliques in $W'\subset U$. We thus  have that $V(\cD) \subset U$ as required and the order of $\cD$ is $|Q|+|X|\leq 2m$.   
\end{proof}

\subsection{The existence of shrinkable orchards of larger order} \label{sec:slightly larger}
Lemma~\ref{lem:diamond trees with low degrees} gives us the key to being able to push the methods above (which culminated in Lemma~\ref{lem:recover min deg}) to be able to handle orchards with larger order. We remark that the flexibility given by  dealing with $(k,m)$-systems and  Observation~\ref{obs:setsystemfromorchard} is necessary in order to handle this extension. Indeed, this is due to Lemma~\ref{lem:diamond trees with low degrees} only giving control over the degree to a subset of the removable vertices of the diamond tree generated.

\begin{proof}[Proof of Proposition \ref{prop:small shrinkable}]
By Proposition~\ref{prop:verysmallorchards} we can focus on the case that \[ p^{r-1}n^{1-r^3\gamma}\leq m \leq \min\{ p^{r-2}n^{1-2r^3\gamma},n^{7/8}\}.\] We fix $\eps>0$ small enough to apply Lemma~\ref{lem:recover min deg} with $\alpha_{\ref{lem:recover min deg}}=\alpha':=\alpha/4$ and $\gamma_{\ref{lem:recover min deg}}=\gamma$ and small enough to apply Lemma~\ref{lem:diamond trees with low degrees} with $\gamma_{\ref{lem:diamond trees with low degrees}}=\gamma':=r^3\gamma$ and $\eta_{\ref{lem:diamond trees with low degrees}}=\eta<\alpha/8$.  Finally we fix some $k\in \mathbb{N}$ such that $\alpha n \le km \le 2 \alpha n$ and note that we have $pk\geq n^\gamma$ due to our upper bound on $m$. By repeatedly applying Lemma  \ref{lem:diamond trees with low degrees}, we find a $(2k,m)$-orchard $\bc{O}_0$ with $V(\bc{O}_0)\subset U$ and each $\cD=(T,R,\Sigma)\in \bc{O}_0$ has the property that there exists some distinguished subset $Q_{\cD}\subset R$ of removable vertices such that $|Q_{\cD}|=m$ and all but at most $\eta m$ vertices $v$ in $V(G)$ have that $\deg_{Q_{\cD}}(v)\leq p^{r-1}n^{1-\gamma'}$. Indeed, we can find $\bc{O}_0$ by sequentially choosing diamond trees and deleting their vertices from $U$, using that $|V(\bc{O}_0)|\leq 4r \alpha n$ at all times in this process and so $|U\setminus V(\bc{O}_0)|\geq   n/4 $ and we can apply Lemma \ref{lem:diamond trees with low degrees}. 

Now we will crop our orchard  $\bc{O}_0$ to arrive at an orchard for which we can apply Lemma~\ref{lem:recover min deg} to subsets of removable vertices. Similarly to the proof of Lemma~\ref{lem:recover min deg}, we do this by a process of `cleaning up';  losing diamond trees in the orchard which have lots of removable vertices which are atypical. So let $B_1\subset V(G)$ be the set of vertices $v\in V(G)$ such that $\deg_{Q_{\cD}}(v)>p^{r-1}n^{1-\gamma'}$ for some 
$\cD\in \bc{O}_0$ as above. It follows that $|B_1|\leq \eta m \cdot 2k\leq  \alpha n/4$. Next we delete $\cD'$ 
from $\bc{O}_0$ if $|B_1\cap Q_{\cD'}|\geq m/2$. Due to our upper bound on $|B_1|$, we delete at most   $ k/2$ of the diamond trees $\cD'$ from $\bc{O}_0$. Let the resulting suborchard be $\bc{O}_1\subseteq \bc{O}_0$ and for each diamond tree $\cD=(T,R,\Sigma)\in \bc{O}_1$ define a distinguished subset $S_{\cD}\subset Q_{\cD}\subseteq R$  of removable vertices such that $|S_{\cD}|=m/2$ and 
 \begin{equation} \label{eq: small deg to Ss}\deg_{S_{\cD}}(v)\leq p^{r-1}n^{1-\gamma'} = p^{r-1}n^{1-r^3\gamma}\mbox{ for all } \cD\in \bc{O}_1 \mbox{ and all }v\in \bigcup_{\cD'\in \bc{O}_1}S_{\cD'}.
\end{equation}
Let $\bc{O}_2$ be an arbitrary suborchard of $\bc{O}_1$ so that $|\bc{O}_2|=(1+\gamma)k$. Moreover, let $\Lambda'=\{S_\cD:\cD\in\bc{O}_2\}$ be the $((1+\gamma)k,m/4)$-system defined by the distinguished subsets of removable vertices for the $K_r$-diamond trees in $\bc{O}_2$. Now due to \eqref{eq: small deg to Ss}, we have that Lemma \ref{lem:recover min deg} gives the existence of some $\gamma$-shrinkable (with respect to $r$) subsystem $\Lambda\subset \Lambda'$. Taking $\bc{O}:=\{\cD\in \bc{O}_2:S_\cD\in \Lambda\}$ thus gives a $\gamma$-shrinkable $(k,m)_r$-orchard as required, appealing to Observation~\ref{obs:setsystemfromorchard}.  
\end{proof}

\section{Shrinkable orchards of large order} \label{sec:large order shrinkable}
In this section, we establish the existence of shrinkable orchards with large order, proving Proposition~\ref{prop:large shrinkable}. Our approach is to find an orchard such that the $K_r$-hypergraph $H$ generated by the orchard is very dense. This allows us to apply Lemma~\ref{lem:vertex fans} in many subhypergraphs of $H$. Coupled with Theorem~\ref{thm:almost perfect matchings via greedy pfms}, this will imply that the orchard is shrinkable. As in the previous section we begin in Section~\ref{sec:large conditions for shrinkability} by using these results on fractional matchings to deduce conditions on an orchard which guarantee shrinkability. We will then show in Section~\ref{sec:popular} that we can appeal to Proposition~\ref{prop: choosing removable vertices} to generate diamond trees whose removable vertices are contained in many copies of $K_r$. This will then allow us to prove the existence of shrinkable orchards of large order in Section~\ref{sec:largest shrinkable}. As in Section~\ref{sec:small order shrinkable} however, this first  argument will fall short of the range of orders needed in Proposition~\ref{prop:large shrinkable}. The rest of the section is thus concerned with extending our methods to capture more orders. This leads us to a process which generates an orchard in two rounds. The outcome of the first round is discussed in Section~\ref{sec:first round} and building on this, in Section~\ref{sec:second round} we detail properties of the orchard after a second round of generation. Finally in Section~\ref{sec:slightly smaller}, we show that by generating orchards via this two phase process, we end up with orchards which are shrinkable. This allows us to complete the proof of Proposition~\ref{prop:large shrinkable}.

\subsection{A density condition which guarantees shrinkability} \label{sec:large conditions for shrinkability}

We begin by applying Lemma \ref{lem:vertex fans} and Theorem~\ref{thm:almost perfect matchings via greedy pfms} to give a density condition which we can use to show that an orchard is shrinkable. This transforms our problem into finding orchards which satisfy this condition. 

\begin{lem} \label{lem:density condition for shrinkable}
 For all   $3\leq r\in \NN$ and   $0<\gamma< 1/(2r^3) $, there exists a $k_0\in \NN$ such that the following holds. Suppose that $\bc{O}$ is a $(k,m)_r$-orchard in a graph $G$ with $k_0\leq k\in \NN$ and $m\in \NN$.  For a diamond tree $\cD\in\bc{O}$, let $R_\cD$ denote its removable vertices and for a suborchard $\bc{O}'\subset \bc{O}$, let
$R(\bc{O}'):=\cup_{\cD\in \bc{O}'}R_\cD$ denote the union of the sets of removable vertices of diamond trees in $\bc{O}'$. Suppose that the following condition holds:
\begin{equation} \label{density cond}
      \begin{minipage}{0.85\textwidth}
{For any }$\cD\in \bc{O}$ { and } $\bc{P}\subset \bc{O}\setminus \{\cD\}$ { such that } $|\bc{P}|\geq k/(4r)$,  { there exists a  suborchard } $\bc{P}^*=\bc{P^*}(\cD,\bc{P})\subset \bc{P}$ { such that } $|\bc{P}^*|\leq k^{1-r^3\gamma}$ 
    { and for \emph{any} disjoint suborchards } $\bc{O}_1, \ldots, \bc{O}_{r-2}\subset \bc{P}\setminus \bc{P}^*$, 
    { with }  $|\bc{O}_i|\geq k^{1-r^3\gamma}$ for $i \in [r-3]$ and $|\bc{O}_{r-2}|\geq \gamma k$, 
    there is a copy of $K_r$ in $G$ traversing $R_\cD$, $R(\bc{P^*})$ and  $R(\bc{O}_i)$ for $i\in[r-2]$. 
    \end{minipage}
\end{equation}
Then $\bc{O}$ is $\gamma$-shrinkable. 
\end{lem}

Let us take a moment to digest the density condition \eqref{density cond}. For simplicity, one can think of $\bc{P}^*$ being a single diamond tree $\cD^*=\cD^*(\cD,\bc{P})$. Indeed this is the setting that we will work in first when applying Lemma~\ref{lem:density condition for shrinkable}. Simplifying further and just focusing on the case that $r=3$, the condition \eqref{density cond} then translates as having that for any $K_3$-diamond tree $\cD$ in the orchard and large suborchard $\bc{P}\subset \bc{O}$, there is some  diamond tree $\cD^*\in\bc{P}$ so that the pair $\{\cD,\cD^*\}$ has high degree in the $K_3$-hypergraph generated by $\bc{P}$. Indeed for any small linear sized $\bc{O}_1\subset \bc{P}$, there is a hyperedge in $H(\bc{O})$ containing $\cD, \cD^*$ and a diamond tree in $\bc{O}_1$.  In general, when $r\geq 4$, we need to guarantee traversing $K_r$s when some of the sets we look to traverse are smaller than linear (size $k^{1-r^3\gamma})$. Also later on we will need the full power of Lemma \ref{lem:density condition for shrinkable} which allows us to choose the $\bc{P}^*$ as  a small suborchard as opposed to a single diamond tree. We now prove the lemma. 
\begin{proof}[Proof of Lemma~\ref{lem:density condition for shrinkable}]
Let $\bc{Q}\subset \bc{O}$ be an arbitrary suborchard of $\bc{O}$ of size $\gamma k$. We will show that $\bc{O}$ is shrinkable with respect to $\bc{Q}$. So fix some arbitrary suborchard $\bc{Q}'\subset \bc{Q}$ and let $H:=H(\bc{O}\setminus \bc{Q}')$ be the $K_r$-hypergraph generated by $\bc{O}\setminus \bc{Q}'$. We have to show that $H$ has a matching covering all but at most $k^{1-\gamma}$ vertices of $H$.

In order to show the existence of a large matching in $H$, as we did in  Lemma~\ref{lem:local to shrink}, we appeal to Theorem~\ref{thm:almost perfect matchings via greedy pfms}. So let us fix  $N=|V(H)|$ and note that as $N\geq (1-\gamma)k$, by choosing $k_0$ to be large, we can assume that $N$ is sufficiently large in what follows. Now fix some $2$-uniform graph $J$ on $V(H)$ of maximum degree at most $N^{r^2\gamma}$.  If we can show that $H\setminus H_J$ contains  a perfect fractional matching, then we are done by Theorem~\ref{thm:almost perfect matchings via greedy pfms} as, because $J$ was arbitrary, the theorem guarantees a matching covering all but at most $N^{1-\gamma}\leq k^{1-\gamma}$ vertices of $H$. 

In order to prove the existence of a perfect fractional matching in $H\setminus H_J$, we appeal to Lemma~\ref{lem:vertex fans}, fixing $M:=N/(2r)$. Thus, we need to show that given any $K_r$-diamond tree $\cD\in V(H)=\bc{O}\setminus \bc{Q}$ and suborchard $\bc{P}_0\subset V(H)\setminus \{\cD\}$ with $|\bc{P}_0|\geq M$, there is an edge in $H\setminus H_J$ containing $\cD$ and $r-1$ $K_r$-diamond trees in $\bc{P}_0$. So fix such a $\cD$ and $\bc{P}_0$. Let $\bc{P}:=\bc{P}_0\setminus N^J(\cD)$. Therefore, we have that 
\[|\bc{P}|\geq |\bc{P}_0|-|N^J(\cD)|\geq \frac{N}{2r}-N^{r^2\gamma}\geq \frac{(1-\gamma)k}{2r}-k^{r^2\gamma}\geq \frac{k}{4r},\]
for $k$ sufficiently large. Hence by condition \eqref{density cond}, we have the existence of some $\bc{P}^*=\bc{P}^*(\cD,\bc{P})\subset \bc{P}$ as in the hypothesis. Now we will iteratively define $\bc{O}_i$ for $1\leq i\leq r-2$ as follows. We begin by fixing $\bc{P}'=\bc{P}$ and defining $\bc{Q}_0:=\cup_{\cC\in \bc{P}^*}(N^J(\cC)\cup\{\cC\})$. Now for $1\leq i \leq r-2$, we update $\bc{P}'$ by removing any diamond trees in $\bc{Q}_{i-1}$ from $\bc{P}'$ and then define $\bc{O}_i$ to be an arbitrary suborchard of $\bc{P}'$ of size $k^{1-r^3\gamma}$ if $i\in[r-3]$, and of size $\gamma k $ if $i=r-2$. If $i=r-2$ we then end this process. If $i<r-2$, we define $\bc{Q}_i:=\cup_{\cC\in \bc{O}_i}(N^J(\cC)\cup\{\cC\})$ and move to the next index. 

Let us check that we are successful in each round. Indeed this follows because at the beginning of step $i$ in the process, $\bc{P}'$ has size 
\[|\bc{P}'|\geq \frac{k}{4r}-ik^{1-r^3\gamma}(1+N^{r^2\gamma})\geq \frac{k}{4r}-rk^{1-r^3\gamma+r^2\gamma}\geq \gamma k\geq k^{1-r^3\gamma}, \]
for large $k$. Therefore there is always space in $\bc{P}'$ to choose our suborchard $\bc{O}_i$ at each step $i$. Now the condition~\eqref{density cond} gives a copy of $K_r$ in $G$ traversing $R_\cD$, $R(\bc{P^*})$ and  $R(\bc{O}_i)$ for $i\in[r-2]$. This thus gives a hyperedge $e$ in the $K_r$-hypergraph $H=H(\bc{O}\setminus \bc{Q}')$ which has one vertex as $\cD$, one vertex in $\bc{P}^*\subset \bc{P}_0$ and one vertex in each of the $\bc{O}_i\subset \bc{P}_0$. Moreover this edge $e$ lies in $H\setminus  H_J$. Indeed, by our construction of $\bc{P}^*$  and the $\bc{O}_i$, there is no edge in $J$ between any pair of distinct sets in the family $\{\{\cD\},\bc{P}^*, \bc{O}_1,\ldots, \bc{O}_{r-2}\}$.  We have therefore established the existence of a perfect fractional matching in $H\setminus H_J$ due to Lemma~\ref{lem:vertex fans} which implies that $\bc{O}$ is $\gamma$-shrinkable as detailed above. 
\end{proof}

Lemma~\ref{lem:density condition for shrinkable} gives a route to proving the existence of shrinkable orchards. Indeed, if the sets of vertices which arise as pools of removable  vertices of suborchards are sufficiently large, then appealing to Corollary~\ref{cor:transversalcliques} can give the required transversal copy of $K_r$ in $G$, so that \eqref{density cond} is satisfied. However, we cannot immediately derive such results because the size of the sets required in \eqref{density cond} are too small. In particular, \eqref{density cond} forces only one set (namely $R(\bc{O}_{r-2})$) to be linear in size whilst all other sets that feature can have sublinear size. This is troublesome because the examples we have from Corollary~\ref{cor:transversalcliques} to generate transversal copies of $K_r$, require at least two of the sets involved to be linear. Indeed, it can be seen from the more general Lemma~\ref{lem:generaltransversalcliques} that we cannot do any better than this. That is, in order to use Definition~\ref{def:bijumbled} and our condition on $\beta$ to derive the existence of a copy of $K_r$ that traverses a family of sets, at least two of the sets in the family must be linear in size. Therefore in order to apply Lemma~\ref{lem:density condition for shrinkable} and derive the existence of shrinkable orchards, we have to obtain orchards with some additional structure. We start by exploring properties of singular diamond tress that we can guarantee.

\subsection{Popular diamond trees}
\label{sec:popular}

As was the case when we were interested in proving the existence of shrinkable orchards with small order, Proposition~\ref{prop: choosing removable vertices}  gives a powerful tool for proving the existence of diamond trees with additional desired properties. Here 
 we show that we can choose a diamond tree  so that there are many copies of $K_r$ formed with its removable vertices. 

\begin{lem} \label{lem:diamond trees with many transversal triangles}
For any $3\leq r\in \NN$ and $0<\alpha < 1/2^{12r}$  there exists an $\eps>0$ such that the following holds for any $n$ vertex $(p,\beta)$-bijumbled graph $G$ with $\beta\leq \eps p^{r-1}n$  and any vertex subset $U\subseteq V(G)$ with $|U|\geq n/4$. Suppose that $m\in \mathbb{N}$ with \[\max\{p^{1-r},p^{r-1}n\}\leq m \leq n^{7/8}\]  and we have set families $\cW_0,\cW_1,\ldots, \cW_{r-1}\subset 2^{V(G)}$ such that:
\begin{enumerate}
    \item \label{cond:W_0 size} $|W_0|\geq \alpha p^{r-1} n$ for all $W_0\in \cW_0$;
    \item \label{cond:W_i size} $|W_i|\geq \alpha p n$ for all $W_i\in \cW_i$, $1\leq i\leq r-3$;
    \item \label{cond:W_r-2 size} $|W_{r-2}|\geq \alpha n$ for all $W_{r-2}\in \cW_{r-2}$;
    \item \label{cond:small W sets} $\prod_{i=0}^{r-2}|\cW_i|\leq 2^{ \frac{m}{4}}$.
\end{enumerate}
Then there exists a $K_r$-diamond tree $\cD=(T,R, \Sigma)$ in $G[U]$ of order at least $m$ and at most $2m$   such that for \emph{any} choice of sets  $\mathbf{W}=(W_0,\ldots,W_{r-2})\in \cW_0\times \cdots\times\cW_{r-2}$, there is a copy of $K_r$ in $G$ traversing $R$ and the sets  $W_0,\ldots, W_{r-2}$.

\end{lem}
\begin{proof} 
Let us fix $\eps>0$ small enough to apply Proposition~\ref{prop: choosing removable vertices} with $\alpha_{\ref{prop: choosing removable vertices}}=\alpha':=1/2^{3r}$ and Corollary \ref{cor:transversalcliques} with $\alpha_{\ref{cor:transversalcliques}}=\alpha$ as well as being small enough to force $n$ to be sufficiently large. Note that our  lower bound of $\Omega(p^{r-1}n)$ on $m$ and Fact~\ref{fact:dense} imply that $m$ tends to infinity as $n$ tends to infinity and so we can also assume $m$ is sufficiently large in what follows.  We begin by splitting $U$ into disjoint subsets $U'$ and $W'$ arbitrarily so that $|U'|,|W'|\geq n/8=4 \alpha' r n $, noting that this is possible due to our definition of $\alpha'$. We further fix  $d_* :=\alpha'^2p^{r-1} n$.

  Now we apply Proposition \ref{prop: choosing removable vertices} with $z:=\frac{\alpha'^2n}{4}=\frac{n}{2^{6r+2}}$ and fix the sets $X\subset U'$ and $Y\subset U'$ which are output. Note that 
  \begin{equation*}
      |X|\leq \max \left.\begin{cases}\frac{2z}{d_*} &=\frac{1}{2p^{r-1}}  \\ 
      &1\end{cases}\right\}\leq \frac{m}{2} \quad \mbox{ and } \quad |Y| = z-|X| \geq  \frac{z}{2}=\frac{n}{2^{6r+3}}, \end{equation*} for $n$ large. 

Now for each choice of $\mathbf{W}=(W_0,\ldots,W_{r-2})\in \cW_0\times \cdots\times\cW_{r-2}$, we find some subset $Y(\mathbf{W}) \subset Y$ of size $|Y|/2$ such that for every $v\in Y(\mathbf{W})$, there is a copy of $K_{r-1}$ in the neighbourhood of $v$ which traverses $W_0, \ldots,W_{r-2}$. In other words, for every $v\in Y(\mathbf{W})$, there is a copy of $K_r$ traversing $W_0, \ldots,W_{r-2}$ and $\{v\}$.   
We can find $Y(\mathbf{W})$ by repeated applications of Corollary \ref{cor:transversalcliques} \eqref{cor:r-2small}. 
In more detail, we  initiate with $Y_0=Y$ and $Y(\mathbf{W})$ empty and in each step we find a copy of $K_r$ traversing $W_0,\ldots, W_{r-2}$ and $Y_0$. Taking $v$ to be the\footnote{Here we refer to \emph{the} vertex that lies in $Y_0$ although there may be several (if the $W_i$ intersect the $Y_0$). What we mean here is the vertex $v$ in the copy of $K_r$ which is assigned to $Y_0$ by virtue of the copy being traversing.} vertex of this $K_r$ that lies in $Y_0$, we add $v$ to $Y(\mathbf{W})$, delete it from $Y_0$ and move to the next step. We continue for $|Y|/2$ steps using that the conditions of Corollary \ref{cor:transversalcliques} \eqref{cor:r-2small} 
are satisfied at each step. Indeed this is  due to the lower bounds on the sizes of $W_i$ in conditions \eqref{cond:W_0 size}, \eqref{cond:W_i size} and \eqref{cond:W_r-2 size} of this lemma and the fact that $|Y_0|\geq |Y|-|Y(\mathbf{W})|\geq |Y|/2 \geq  \alpha n$ throughout, using our upper bound on $\alpha$ and our lower bound on $|Y|$ here. 

\vspace{2mm}

Similarly to the proof of Lemma  \ref{lem:diamond trees with low degrees}, we now take $Q$ to be a random subset of $Y$ by taking each vertex of $Y$ into $Q$ independently with probability $p':=\frac{5m}{4|Y|}$. Thus $\EE \big[|Q|\big]=5m/4$ and by Theorem \ref{thm:chernoff}, we have that $m\leq|Q| \leq 3m/2$ with probability at least $(1-2e^{-m/60})$. Furthermore, for any fixed $\mathbf{W}\in \cW_0\times \cdots\times\cW_{r-2}$, we have that \[\EE \big[|Q \cap Y(\mathbf{W})|\big]= p'|Y(\mathbf{W})|= \frac{5m }{8}.\]
Applying Theorem \ref{thm:chernoff} again implies that the probability that $|Q \cap Y(\mathbf{W})|=0$ is less than $e^{-5m/16}$. Therefore using that $\prod_{i=0}^{r-2}|\cW_i|\leq 2^{ \frac{m}{4}}$  and appealing to a union bound, we  can conclude that whp as $n$ (and hence $m$) tend to infinity, we have that $m\leq |Q| \leq 3m/2$ and $Q \cap Y(\mathbf{W})\neq \emptyset$ for \emph{all} choices of $\mathbf{W}\in \cW_0\times \cdots\times\cW_{r-2}$. So for sufficiently large $n$  we can fix such an instance $Q\subset Y$ and taking $R:=X \cup Q$ we have that  a $K_r$-diamond tree $\cD=(T,R,\Sigma)$ with removable set of vertices $R$ is guaranteed by Proposition \ref{prop: choosing removable vertices}.  We claim that $\cD$ satisfies all the necessary conditions. Indeed, the fact that the order of $\cD$ lies between $m$ and $2m$ follows from the fact that $m \leq |Q| \leq 3m/2$ and $|X| \leq m/2$ whilst the fact that $Q \cap Y(\mathbf{W})\neq \emptyset$ for each choice of $\mathbf{W}=(W_0,\ldots,W_{r-2})$  guarantees that we have a copy of $K_r$ traversing $Q \subset R$ and the sets  $W_0, \ldots, W_{r-2}$.
\end{proof}

\subsection{The existence of shrinkable orchards of large order} \label{sec:largest shrinkable}

Using Lemma~\ref{lem:diamond trees with many transversal triangles} to generate the diamond trees that form our orchard, we can prove that the orchard generated satisfies the condition of Lemma~\ref{lem:density condition for shrinkable} and hence is shrinkable. This gives the following. 

\begin{prop} \label{prop:easy large order shrinkable}
For any $3\leq r\in \NN$ and $0<\alpha,\gamma <1/2^{12r} $  there exists an $\eps>0$ such that the following holds for any $n$ vertex $(p,\beta)$-bijumbled graph $G$ with $\beta\leq \eps p^{r-1}n$  and any vertex subset $U\subseteq V(G)$ with $|U|\geq n/2$. For any $m\in \mathbb{N}$ with  \[\max\{p^{1-r},p^{r-1}n\}\leq m \leq n^{7/8},\]  there exists a $\gamma$-shrinkable $(k,m)_r$-orchard $\bc{O}$ in $G[U]$ with $k\in\NN$ such that $\alpha n\le km\le 2\alpha n$. 
\end{prop}

\begin{proof}
Fix $\eps>0$ small enough to apply Lemma~\ref{lem:density condition for shrinkable} with $\gamma_{\ref{lem:density condition for shrinkable}}=\gamma$ and Lemma~\ref{lem:diamond trees with many transversal triangles} with $\alpha_{\ref{lem:diamond trees with many transversal triangles}}=\alpha'=\alpha \gamma$. Fix some $k\in \NN$ such that $\alpha n\le k \le 2\alpha n$. We also ensure that $\eps$ is small enough to force $n$ (and hence $k$, due to our upper bound on $m$) to be sufficiently large in what follows. Now we begin by noticing that $k\leq m/8r$. Indeed we have that if $p\geq n^{-1/(2r-2)}$, then 
\[p^{1-r}\leq \sqrt{n}\leq p^{r-1}n\leq m,\]
while if $p\leq n^{-1/(2r-2)}$, then  
\[p^{r-1}n\leq \sqrt{n}\leq p^{1-r}\leq m. \]
Therefore, for any $p$ we have that $m\geq \sqrt{n}$ and $k\le 2\alpha n/m\leq 2^{-11r}\sqrt{n}\leq m/8r$.

Now we turn to finding our $(k,m)_r$-orchard in $G[U]$. We do this by finding one diamond tree at a time as follows. For $1 \leq i \leq k$, fix $U_i:=U \setminus (\bigcup_{i'<i}V(\cD_{i'}))$ and note that $|U_i|\geq  |U|-2\alpha r n\geq n/4$ throughout due to our condition on $\alpha$. We then apply Lemma  \ref{lem:diamond trees with many transversal triangles} to find a diamond tree $\cD_i=(T_i,R_i,\Sigma_i)$ such that $V(\cD_i)\subset U_i$ and for any choice of  $i'\in [i-1]$ and disjoint subsets $I_1,\ldots I_{r-2} \subset [i-1] \setminus \{i'\} $ with $|I_j|\geq p k$ for $1\leq j\leq r-3$, and $|I_{r-2}|\geq \gamma k $ we have that there is a copy of $K_r$ traversing $R_{i}$, $R_{{i'}}$ and the sets  $\bigcup_{\ell\in I_j}R_{\ell}$ for $j\in [r-2]$. The existence of such a $\cD_i$ follows from Lemma \ref{lem:diamond trees with many transversal triangles}. 
Indeed, we  define $\cW_0=\{R_{i'}:i'\in[i-1]\}$, $\cW_j=\{\bigcup_{\ell\in I'}R_{\ell}:I'\subset [i-1],|I'|\geq pk\}$ for $1\leq j \leq r-3$ and finally we define $\cW_{r-2}=\{\bigcup_{\ell\in I'}R_{\ell}:I'\subset [i-1],|I'|\geq \gamma k\}$. We need to check that conditions \eqref{cond:W_0 size}-\eqref{cond:small W sets} of Lemma~\ref{lem:diamond trees with many transversal triangles} are satisfied. Indeed condition \eqref{cond:W_0 size} follows from our lower bound on $m$ whilst \eqref{cond:W_i size} and \eqref{cond:W_r-2 size} follow from the fact that $km\geq \alpha n$ and our definition of $\alpha'$. Finally note that each choice of a set in any of the $\cW_j$ comes from a subset of $[i-1]$. Hence we can upper bound $\prod_{j=0}^{r-2}|\cW_j|$ by $(2^{i})^{r-1}\leq 2^{rk}$. As discussed in the opening paragraph, we have that $k\leq m/(8r)$ and so condition \eqref{cond:small W sets} of Lemma \ref{lem:diamond trees with many transversal triangles} is also satisfied. Thus Lemma~\ref{lem:diamond trees with many transversal triangles} succeeds in finding the necessary  $K_r$-diamond tree at every step of this process. 

Let $\bc{O}=\{\cD_1,\ldots, \cD_k\}$ be the orchard obtained by this process. We claim that $\bc{O}$ is $\gamma$-shrinkable and to show this we appeal to Lemma~\ref{lem:density condition for shrinkable} and so need to show that the density condition \eqref{density cond} is satisfied by $\bc{O}$.  So fix some arbitrary $\cD_i\in \bc{O}$ and $\bc{P}\subset \bc{O}\setminus \{\cD_i\}$ with $|\bc{P}|\geq k/(4r)$. We then define $\cD^*=\cD^*(\cD_i,\bc{P})$ (this plays the role of $\bc{P}^*$ in \eqref{density cond}) to be the diamond tree in $\bc{P}$ with the highest index. That is we define $i^*:=\max\{i':\cD_{i'}\in \bc{P}^*\}$ and set $\cD^*=\cD_{i^*}$. Note that we may have that $i^*<i$ but this will not be a problem. We claim that condition \eqref{density cond} is satisfied with this choice of $\cD^*$. Indeed, let $\bc{O}_1,\ldots,\bc{O}_{r-2}\subset \bc{P}^*\setminus \{\cD^*\}$ be disjoint suborchards satisfying the lower bounds on the sizes given by \eqref{density cond}. For each $j\in[r-2]$, define $I_j:=\{i':\cD_{i'}\in \bc{O}_j\}$. Therefore we have that $|I_{r-2}|\geq \gamma k$. For $1\leq j \leq r-3$ we have that $|I_{j}|\geq k^{1-r^3\gamma}\geq pk$.  This follows from the fact that \[k^{-r^3\gamma}\geq k^{-\frac{1}{2^r}}\geq n^{-\frac{1}{2^{(r+1)}}}\geq n^{-\frac{1}{8(r-1)}}\geq p,\]
where we used the upper bound on $\gamma$ in the first inequality, the fact that $k\leq \sqrt{n}$ in the second inequality (see the opening paragraph of the proof), and the fact that $p^{r-1}n\leq m \leq n^{7/8}$ in the last inequality. Now relabelling $\{i,i^*\}$ as $\{\ell_0,\ell_1\}$ so that $\ell_0<\ell_1$, we have that at the point of choosing $\cD_{\ell_1}$, we guaranteed that there was a $K_r$ traversing $R_{\ell_1}$, $R_{\ell_0}$ and the sets $R(\bc{O}_j)=\cup_{i'\in I_j}R_{i'}$ for $j\in[r-2]$. By Lemma~\ref{lem:density condition for shrinkable} this completes the proof that $\bc{O}$ is $\gamma$-shrinkable. 
\end{proof}

Proposition~\ref{prop:easy large order shrinkable} establishes Proposition~\ref{prop:large shrinkable} when $G$ is very dense. However when $G$ is sparse (when $p\leq n^{-1/(2r-2)}$ to be specific), the lower bound of $m\geq p^{1-r}$ takes over and we are left with a gap between the range covered by Proposition~\ref{prop:easy large order shrinkable} and the desired range of Proposition~\ref{prop:large shrinkable}. Tracing the condition that $m=\Omega(p^{1-r})$ back through the proof, we can see that this was necessary in order to prove Lemma~\ref{lem:diamond trees with many transversal triangles}. There, we used our key Proposition~\ref{prop: choosing removable vertices} to generate a diamond tree where we had a large pool $Y$ of vertices which were candidates for being removable vertices. In order to establish the existence of the cliques we need in Lemma~\ref{lem:diamond trees with many transversal triangles}, we needed $Y$ to be linear in size. The sticking point then comes from the fact that Proposition~\ref{prop: choosing removable vertices} can only guarantee a maximum factor of $O(p^{r-1}n)$ between the size of the pool of vertices $Y$ and the order of the diamond tree that we generate. Indeed, in Proposition~\ref{prop: choosing removable vertices} we are forced to include the set $X$ in the removable vertices of the diamond tree we generate and when $Y$ is linear in size, $X$ could have size as large as $\Omega(p^{1-r})$. It is unclear how one would improve on this and find diamond trees with smaller order that are still contained in sufficiently many copies of $K_r$. 

Thankfully, there is  a way to circumvent this issue and apply our methods to close the gap in the range of orders nonetheless. The key idea is to replace the diamond tree generated by Lemma~\ref{lem:diamond trees with many transversal triangles} with a \emph{set of diamond trees}, that is, a small suborchard. Indeed, by grouping together diamond trees, we can decrease their order but guarantee that the collective pool of potential removable vertices for the group is still linear in size. Through following a similar proof to that of Lemma~\ref{lem:diamond trees with many transversal triangles}, this has the outcome of being able to guarantee many copies of $K_r$ which contain a vertex in the removable vertices of \emph{one of} the diamond trees in the group. Moreover, in the proof of Proposition~\ref{prop:easy large order shrinkable}, we crucially used that we could generate diamond trees from Lemma~\ref{lem:diamond trees with many transversal triangles} to establish the density condition~\eqref{density cond} of Lemma~\ref{lem:density condition for shrinkable}. We chose an appropriate $\cD^*$ and used that it had been generated by Lemma~\ref{lem:diamond trees with many transversal triangles} to prove the required existence of transversal copies of $K_r$. However, Lemma~\ref{lem:density condition for shrinkable} allows for us to use a much larger suborchard $\bc{P}^*$ for this condition as opposed to a single diamond tree. Therefore there is hope to incorporate the idea of using a suborchard instead of a single diamond tree in Lemma~\ref{lem:diamond trees with many transversal triangles} whilst maintaining the  overall scheme of the proof. There are some further difficulties to overcome but on a high level, this is the approach we follow in the next sections to establish Proposition~\ref{prop:large shrinkable}.

\subsection{Preprocessing the orchard}
\label{sec:first round}

As discussed above, in order to prove Proposition~\ref{prop:large shrinkable} and remove the condition that $m=\Omega(p^{1-r})$ from Proposition~\ref{prop:easy large order shrinkable}, we need to replace the role played by $\cD^*$ in the proof by a small suborchard $\bc{P}^*$. This allows us to prove an analogue of Lemma~\ref{lem:diamond trees with many transversal triangles}, where one now finds an orchard whose collective set of removable vertices lie in many copies of $K_r$. Our shrinkable orchard then, will be formed as the union of  many of these  smaller orchards.  Indeed, in what follows we will split $k$ as $k=\ell t$ and will aim to have $t$ smaller $(\ell,m)_r$-orchards contributing to our shrinkable orchard $\bc{O}$.  Each of the $(\ell,m)_r$-orchards will have strong connectivity to the rest of the  orchard $\bc{O}$. 

In order to work with the fact that we are splitting $k$ into $t$ sets of size $\ell$, we introduce a two-coordinate index system, with $(i,j)\in[t]\times [\ell]$ indicating that we  are referring to the $j^{th}$ object in the $i^{th}$ subset and we will work through these indices  lexicographically. In more detail,  we let $<_{L}$ denote the lexicographic order on the pairs $(i,j)\in [t] \times [\ell]$. 
That is $(i',j')<_L(i,j)$ if and only if either $ 1 \leq i'\leq i-1$ and $1 \leq j' \leq \ell$ or $i'=i$ and $1 \leq j'\leq j-1$. Furthermore for each  $1 \leq i \leq t$ and $1 \leq j \leq \ell$, we define \[I_{<ij}:=\{(i',j')\in [t] \times [\ell]: (i',j')<_L (i,j) \},\]
to be the indices $(i',j')$ which come before $(i,j)$ in the lexicographic order.

\vspace{2mm}

A hurdle that arises with our new approach is that we lose the symmetry provided by the fact that both $\cD$ and $\cD^*$ in our applications of Lemma~\ref{lem:density condition for shrinkable} were given by singular diamond trees. Indeed, in our proof of Proposition~\ref{prop:easy large order shrinkable}, when verifying the condition~\eqref{density cond} of Lemma~\ref{lem:density condition for shrinkable}, we use that both the arbitrary diamond tree $\cD=\cD_i$ and the diamond tree $\cD^*=\cD^*(\cD_i,\bc{P})$ that we can  choose, were generated using Lemma~\ref{lem:diamond trees with many transversal triangles}. We now hope to generate our suborchards $\bc{P}^*$ using an equivalent to Lemma~\ref{lem:diamond trees with many transversal triangles} and this will  mean that we can no longer switch the roles of $\cD$ and $\bc{P}^*$ when appealing to the conclusion of (the proof method of) Lemma~\ref{lem:diamond trees with many transversal triangles}. In particular, this  places a higher demand on the properties we need to conclude of our $(\ell,m)_r$-suborchards. 

In more detail, we need to generate suborchards which are highly connected to \emph{all} the other vertices of the $K_r$-hypergraph $H(\bc{O})$. Therefore it no longer suffices to build our orchard in a linear fashion, choosing diamond trees (or indeed suborchards) to be well connected (in terms of the $K_r$-hypergraph) with previously chosen diamond trees. We will instead generate our orchard in two rounds.  In the first round we fix a part of each diamond tree and using Proposition~\ref{prop: choosing removable vertices}, provide large pools of vertices which can extend the parts of the diamond trees chosen so far, which we will then do in the second round.  Lemma~\ref{lem: first round} details the outcome we draw from this preprocessing first round. 

\begin{lem} \label{lem: first round}
For any $3\leq r\in \NN$ and $0<\alpha<1/2^{12r}$  there exists an $\eps>0$ such that the following holds for any $n$ vertex $(p,\beta)$-bijumbled graph $G$ with $\beta\leq \eps p^{r-1}n$, any vertex subset $U\subseteq V(G)$ with $|U|\geq n/2$ and any $k,m,t, \ell\in \mathbb{N}$ such that $k=t\ell$, \[ \alpha n\leq km\leq 2\alpha n \quad \mbox{ and } \quad \ell m\geq p^{1-r}. \]

There exists vertex sets $Z_{ij}, Y_{ij} \subset U$  and matchings of  $(r-1)$-cliques,  $ \Pi_{ij}, \Upsilon_{ij} \subset K_{r-1}(G[U])$ for each  $i\in [t]$ and $j\in [\ell]$  such that   the copies of $K_{r-1}$ in each  $\Upsilon_{ij}=:\{S_v:v\in Y_{ij}\}$  are indexed by the vertices in $Y_{ij}$ and such that the  conditions \ref{cond: Z and Pi size} through \ref{cond: finding diamond trees} below are satisfied for all $1\leq i\leq t$ and $1 \leq j \leq \ell$.  
\begin{enumerate}[start=1,label={\bfseries$(\arabic*_{ij})$ }]
    \item \label{cond: Z and Pi size} 
    We have that $ |Z_{ij}|=m$ and $|\Pi_{ij}|=|Z_{ij}|-1$. 
    \item \label{cond: Y size}
    We have that  $ |Y_{ij}|= |\Upsilon_{ij}| =\sqrt{\alpha}n/\ell$. 
    \item \label{cond: internally disjoint} We have that the vertex sets $Z_{ij}$, $Y_{ij}$, $V(\Pi_{ij})$ and $V(\Upsilon_{ij})$ are all disjoint from each other.  
     \item \label{cond: external disjoint}
     We have that $A \cap A'=\emptyset$ for any choice of $A\in \{Z_{ij},V(\Pi_{ij}), Y_{ij}, V(\Upsilon_{ij})\}$ and \footnote{Crucially, we do not require that $A$ is disjoint from \emph{all} $Y_{i'j'}$ and $V(\Upsilon_{i'j'})$, only those that are in the same subfamily indexed by $i$. } \[A'\in \{Z_{i'j'}, V(\Pi_{i'j'}):(i',j')\in I_{<ij}\} \cup \{Y_{ij'}, V(\Upsilon_{ij'}): 1 \leq j' \leq j-1\}.\]
      \item \label{cond: finding diamond trees}
      For any choice of $\tilde{Y}$ such that $ \tilde{Y}\subseteq Y_{ij}$, there exists a $K_r$-diamond tree $\cD=(T,R,\Sigma)$ such that $R=Z_{ij} \cup \tilde{Y}$  and  $\Sigma=\Pi_{ij}\cup\tilde{\Upsilon}_{ij}$, where $\tilde{\Upsilon}_{ij}\subset \Upsilon_{ij}$ is defined to be \[\tilde{\Upsilon}_{ij}:=\{S_{\tilde{v}}:\tilde{v}\in \tilde{Y}\subset Y_{ij}\}.\] 
\end{enumerate}
\end{lem}

\vspace{2mm}

As mentioned above, in this first round we  put aside part of every single diamond tree in the $(k,m)_r$-orchard we are going to generate, thus partially defining the orchard. We also put aside large pools of vertices which will be used to extend these diamond trees in second round of generating our orchard. 
The fixed parts of the diamond trees chosen in Lemma~\ref{lem: first round}  are the sets $Z_{ij}$ and the interior cliques $\Pi_{ij}$ whilst the pools of potential removable vertices and interior cliques that can be used to extend the diamond trees chosen are given by the sets $Y_{ij}$ and $\Upsilon_{ij}$, respectively. 
 We make sure through the conditions  \ref{cond: Z and Pi size}  that these fixed sub-diamond trees contribute a substantial portion of the final diamond trees that we are shooting for (which will have order between $m$ and $2m$). We also guarantee through the conditions \ref{cond: external disjoint}, that the parts of the diamond trees that we put aside in this preprocessing round do not interfere with each other, in that they are vertex disjoint. Notice also that if we fix $i\in[t]$, then the conditions \ref{cond: external disjoint} for all $j\in[\ell]$ guarantee that the sets $Y_{ij}, V(\Upsilon_{ij})$, $j\in[\ell]$ do not intersect each other. This is important because in the second round of generating our orchard, we will want to extend all the diamond trees in the $i^{th}$ $(\ell,m)_r$-suborchard simultaneously and so we do not want any interference between the choices of the extensions within such a suborchard. Also note that the conditions \ref{cond: Y size} for fixed $i\in[t]$ and all  $j\in[\ell]$, guarantee that the collective pool of potential removable vertices for the $i^{th}$ $(\ell,m)_r$-suborchard (the set $\cup_{j\in[\ell]}Y_{ij}$) is linear in size, as required. Finally, the conditions \ref{cond: finding diamond trees} contain the heart of Proposition \ref{prop: choosing removable vertices}, allowing us to arbitrarily extend any of the diamond trees we have so far using any subsets of the  pools (the $Y_{ij}$) of potential removable vertices and interior cliques (the $\Upsilon_{ij}$) we have put aside. 

Our final remark on the statement of Lemma~\ref{lem: first round} is that we do not require e.g. $Y_{ij}$ and $Y_{i'j'}$ for $i\neq i'$, to be disjoint. Indeed as we have $t$ suborchards and each has a linear collective pool of potential removable vertices, there would not be enough space in the graph to keep these pools disjoint. However, by requiring that the collective pool is much larger than all the vertices in our orchard (that is, much larger than $km$), we guarantee that we will be able to proceed greedily in our second round (Lemma~\ref{lem:orchard with dense triangle hypergraph}) of defining the orchard, always having a large enough set of potential removable vertices at each step.

\begin{proof}[Proof of Lemma~\ref{lem: first round}]
Let us fix $\eps>0$ small enough to apply Proposition~\ref{prop: choosing removable vertices} with $\alpha_{\ref{prop: choosing removable vertices}}=\alpha':=1/2^{2r+1}$. We will find these vertex sets  and matchings of $(r-1)$-cliques algorithmically working through the pairs $(i,j)\in[t]\times[\ell]$ in lexicographic order. 
So let us fix some $(i^*,j^*)\in [t]\times [\ell]$ and suppose that we have already found $Z_{ij}, Y_{ij}, \Pi_{ij}$ and $\Upsilon_{ij}$  such that 
the conditions \ref{cond: Z and Pi size} through \ref{cond: finding diamond trees} are satisfied for all $(i,j)\in I_{<i^*j^*}$. We fix $W^* \subset U$ to be \[W^*:=\left(\bigcup\left\{Z_{ij}\cup V(\Pi_{ij}): (i,j)\in I_{<i^*j^*}
\right\} 
\right) \bigcup \left( \bigcup \left\{Y_{i^*j} \cup V(\Upsilon_{i^*j}): 1 \leq j \leq j^*-1\right\} \right), \]
and let $U^*:= U \setminus W^*$. We use conditions \ref{cond: Z and Pi size} and \ref{cond: Y size} to upper bound the size of $W^*$ as follows. We have that
\[  |W^*|\leq rm((i^*-1)\ell+j^*-1) +\frac{\sqrt{\alpha}rn}{\ell}(j^*-1) 
    \leq rmt\ell+ \sqrt{\alpha}rn 
    \leq (2\alpha +\sqrt{\alpha})rn,\]
using that $mt\ell=mk\leq 2\alpha n$. Hence $|U^*|\geq  n/4$ from our upper bound on $\alpha$. We will find $Z_{i^*j^*},Y_{i^*j^*}\subset U^*$ and $\Pi_{i^*j^*}, \Upsilon_{i^*j^*}\subset K_{r-1}(G[U^*])$ and so condition  {\bfseries{($4_{i^*j^*}$)}} 
will be satisfied. The required vertex  sets $Z_{i^*j^*}$ and $Y_{i^*j^*}$ are  found by an application of Proposition \ref{prop: choosing removable vertices}. So let us 
split $U^*$ into disjoint subsets $U'$ and $W'$ arbitrarily so that $|U'|,|W'|\geq n/8\geq 4 \alpha' r n $, noting that this is possible due to our definition of $\alpha'$. We further fix  $d_* :=\alpha'^2p^{r-1} n$ and $z:=m+\sqrt{\alpha}n/\ell$ and note that $z\leq \alpha'n$ due to the fact that $m\leq 2\alpha n/k\leq 2\alpha n$ and our upper bound on $\alpha$.

So Proposition \ref{prop: choosing removable vertices} gives us that there exists disjoint vertex subsets $X,Y\subset U'\subset U^*$ such that $|X|+|Y|=z$ and $|X|=1\leq m$  or 
\[|X|\leq 2z/d_*\leq \frac{2m}{d_*}+\frac{2\sqrt{\alpha}n}{d_*\ell}\leq \frac{m}{2}+\frac{2\sqrt{\alpha}}{\alpha'^2p^{r-1}\ell}\leq m,\] 
using our upper bound on  $\alpha$ and lower bound on  $\ell m$  in the last inequality. 
As $|X|\leq m$, we can fix some $Z_{i^*j^*}\subset X\cup Y$ such that $X\subseteq Z_{i^*j^*}$ and $|Z_{i^*j^*}|=m$. Therefore letting $Y_{i^*j^*}:=Y\setminus Z_{i^*j^*}$, we have that $|Y_{i^*j^*}|=z-m=\sqrt{\alpha}n/\ell$ and so the size 
requirements on $Z_{i^*j^*}$ in {\bfseries{($1_{i^*j^*}$)}} 
and on $Y_{i^*j^*}$ in  {\bfseries{($2_{i^*j^*}$)}} are both satisfied.  
Moreover, we also have that  part of condition {\bfseries{($5_{i^*j^*}$)}} is  satisfied. Indeed, for some $\tilde{Y}\subset Y_{i^*j^*}$, taking $Y'=\tilde{Y}\cup (Z_{i^*j^*}\setminus X)$, Proposition~\ref{prop: choosing removable vertices} gives that  there is a diamond tree $\cD=(T,R,\Sigma)$ with removable vertices $R=X \cup Y'= Z_{i^*j^*} \cup \tilde{Y}$ and $\Sigma$ a matching of $(r-1)$-cliques in $G[U^*]$. 

Now in order to complete the proof of the lemma, we need to define the matchings of $(r-1)$-cliques $\Pi_{i^*j^*}$ and $\Upsilon_{i^*j^*}$ and reason that the remaining conditions of the lemma are satisfied. This comes from recalling how we proved Proposition \ref{prop: choosing removable vertices} in Section \ref{sec:scattered} (see also Figure~\ref{fig:scattered}).   
There, we applied Lemma \ref{lem:scattereddiamondtrees} to find a large $d_*$-scattered $K_r$-diamond tree $\cD_{sc}=(T_{sc},R_{sc},\Sigma_{sc})$.  
We had that $R_{sc}=X\cup Y$ was the set of removable vertices of $\cD_{sc}$ and $Y\subset R_{sc}$ was the set of leaves in $\cD_{sc}$. The conclusion of Proposition \ref{prop: choosing removable vertices} then followed readily as we could choose which leaves in $Y$ to include in a diamond subtree $\cD$ of $\cD_{sc}$.  
From this proof we see that we can partition $\Sigma_{sc}$ into $\Sigma_{sc}=: \Pi_{i^*j^*} \cup \Upsilon_{i^*j^*}$ where the $(r-1)$-cliques $\Pi_{i^*j^*}$ are interior cliques of the $K_r$-diamond subtree of $\cD_{sc}$ spanned by the removable vertices $Z_{i^*j^*}$. Furthermore, we can label $\Upsilon_{i^*j^*}$ with the vertices in $Y_{i^*j^*}$ so that  {\bfseries{($5_{i^*j^*}$)}} is satisfied. Indeed each vertex $v$ in $Y_{i^*j^*}$ corresponds to a leaf of the diamond tree $\cD_{sc}$ and so there is an interior clique $S_v\in \Sigma_{sc}$ such that any sub diamond-tree which contains the non-leaves $X$ of $\cD_{sc}$ can be extended by adding $v$ to the set of removable vertices and $S_v$ to the set of  interior cliques. As $\cD_{sc}$ is a well-defined $K_r$-diamond tree, we also have that condition  {\bfseries{($3_{i^*j^*}$)}} is satisfied and the size constraints of $\Pi_{i^*j^*}$ and $\Upsilon_{i^*j^*}$ in {\bfseries{($1_{i^*j^*}$)}} and {\bfseries{($2_{i^*j^*}$)}} are also immediate, noting that $|\Pi_{i^*j^*}|= |Z_{i^*j^*}|-1$ as the set of interior cliques of a diamond tree with removable vertices $Z_{i^*j^*}$.  
\end{proof}

\subsection{Completing the orchard}
\label{sec:second round}

We will now  use Lemma~\ref{lem: first round} to generate our orchard. This can be thought of as extending the parts of the diamond trees (the $Z_{ij}$ and $\Pi_{ij}$) which were fixed in Lemma~\ref{lem: first round}. The strategy is very similar to that of Lemma~\ref{lem:diamond trees with many transversal triangles} and Proposition \ref{prop:easy large order shrinkable}. Indeed we take random subsets of the pools of potential vertices in order to guarantee that the $K_r$-hypergraph generated by our final orchard is sufficiently dense. The key difference here is that, as opposed to fixing our orchard one diamond tree at a time, we appeal to Lemma~\ref{lem: first round} to fix part of all the diamond trees in  our orchard and then carry out the extensions on $(\ell,m)_r$-suborchards. That is, we apply the approach of Lemma~\ref{lem:diamond trees with many transversal triangles} on the whole suborchard as opposed to a singular $K_r$-diamond tree. After doing this process for all suborchards we end up with an orchard which generates a dense $K_r$-hypergraph. This is detailed in the following lemma.

\begin{lem} \label{lem:orchard with dense triangle hypergraph}
For any $3\leq r\in \NN$, $0<\alpha <1/2^{12r}$ and $0<\gamma<1$  there exists an $\eps>0$ such that the following holds for any $n$ vertex $(p,\beta)$-bijumbled graph $G$ with $\beta\leq \eps p^{r-1}n$, any vertex subset $U\subseteq V(G)$ with $|U|\geq n/2$ and any $k,m,t, \ell\in \mathbb{N}$ such that   \begin{equation} \label{eq: k t l m conds} k=t\ell, \qquad m\geq p^{r-1}n,  \qquad \ell m\geq p^{1-r}, \qquad \mbox{and} \qquad \alpha n\leq km\leq 2\alpha n. \end{equation}

Then there exists a $(k,m)_r$-orchard $\bc{O}$ in $G$ such that $V(\bc{O})\subset U$ and $\bc{O}$ can be partitioned into suborchards $\bc{Q}_1,\ldots, \bc{Q}_t$ such that  
each $\bc{Q}_i$ with $1 \leq i \leq t$ is a $(\ell,m)_r$-orchard 
and we have the following property. For any $i\in [t]$, any $\cD'\in \bc{O} $, any   suborchard $\bc{Q}'\subseteq \bc{Q}_i$ with $|\bc{Q}'| \geq \ell/4r$ and any set of disjoint suborchards $\bc{O}'_1,\ldots,\bc{O}'_{r-2}\subset \bc{O}$ with $|\bc{O}'_{i'}|\geq pk$ for $i'\in[r-3]$ and $|\bc{O}'_{r-2}|\geq \gamma k$, 
there exists a copy of $K_r$ traversing\footnote{Here as before, for a diamond tree $\cD$,  $R_\cD$ denotes the set of removable vertices of $\cD$ and for a suborchard $\bc{O}'\subseteq \bc{O}$, $R(\bc{O}')$ denotes the union of the set of removable vertices of diamond trees in $\bc{O}'$. That is, $R(\bc{O}'):=\cup_{\cD\in \bc{O}'}R_{\cD}$.} $R_{\cD'}$, $R(\bc{Q}')$ and $R(\bc{O}_{i'})$ for $i'\in[r-2]$.

\end{lem}
\begin{proof}

Fix $\eps>0$ small enough to apply Corollary \ref{cor:transversalcliques} with $\alpha_{\ref{cor:transversalcliques}}=\alpha':=\gamma\alpha$, small enough to apply Lemma~\ref{lem: first round} with $\alpha_{\ref{lem: first round}}=\alpha$ and small enough to force $n$ to be sufficiently large in what follows. We begin by applying Lemma~\ref{lem: first round} to get vertex sets $Z_{ij}, Y_{ij}$ and matchings of $K_{r-1}$ cliques $\Pi_{ij}$ and $\Upsilon_{ij}$ for $(i,j)\in[t]\times [\ell]$ satisfying \ref{cond: Z and Pi size} through \ref{cond: finding diamond trees} as listed in that lemma. 
Now we turn to finding the diamond trees $\cD_{ij}$ for $(i,j)\in[t]\times [\ell]$ which will form our orchard $\bc{O}$, so that the suborchard $\bc{Q}_{i}$ is defined to be $\bc{Q}_i:=\{\cD_{ij}:j\in[\ell]\}$ for each $i\in[t]$. 
We will appeal  in particular to condition  \ref{cond: finding diamond trees} of Lemma \ref{lem: first round} to find each $\cD_{ij}=(T_{ij},R_{ij},\Sigma_{ij})$. In more detail, for each $(i,j)\in[t] \times [\ell]$, we will find $\tilde{Y}_{ij} \subset Y_{ij}$ and apply \ref{cond: finding diamond trees} to find a diamond tree with removable set of vertices $R_{ij}:=Z_{ij} \cup \tilde{Y}_{ij}$.

Now for a set of indices $I'\subseteq [t]\times [\ell]$, we let $Z(I')=\cup_{(i,j)\in I'}Z_{ij}$. In order to guarantee the key property of $\bc{O}$  it suffices that for each $i\in[t]$ we have the following. For any 
 choice of $J\subseteq [\ell]$ with $|J|\geq \ell/4r$ and any choice of  $(i_0,j_0) \in [t]\times [\ell]$ and subsets $ I_{1},\ldots, I_{r-2}\subset [t]\times [\ell]$ with $|I_{i'}|\geq pk$ for $i'\in [r-2]$ and $|I_{r-2}|\geq \gamma k$, the following holds. There   exists a copy of $K_r$ traversing $\cup_{j\in J}\tilde{Y}_{ij}$,   $Z_{i_0j_0}$ and the sets  $Z(I_{i'})$ for $i'\in[r-2]$. This is what we prove in what follows as we select our sets $\tilde{Y}_{ij}$.

We work through the $i\in [t]$ in order. Let $W_0:= \bigcup_{(i,j)\in[t] \times[\ell]}(Z_{ij} \cup V(\Pi_{ij}))$ and initiate with $U_0=U\setminus W_0$. Suppose that we are at some step $i^*\in[t]$ and we have fixed $\cD_{ij}=(T_{ij},R_{ij},\Sigma_{ij})$ for all $i<i^*$. In this step, we will fix $\cD_{i^*j}$ for all $j\in[\ell]$.  
We define $W_{i^*}:=\left(\bigcup_{(i,j): i<i^*}V(\cD_{ij})\right) \cup W_0$. We further define  for each $J \subseteq  [\ell]$,
\[Y_J^{i^*}:=\big\{v\in U \setminus W_{i^*}: v\in Y_{i^*j} \mbox{ for some } j\in [J] \mbox{ and }S_v\in \Upsilon_{i^*j} \cap K_{r-1}\left(G[U\setminus W_{i^*}]\right)\big\}.\]
In words, $Y_J^{i^*}$ keeps track of the vertices $v$ which lie in one of the $Y_{i^*j}$ with $j\in J$ which we can still use, in that the vertex $v$ has not been used in a previous diamond tree and neither have the vertices of its associated copy of $K_{r-1}$, $S_v$. Note that $|W_{i^*}| \leq  4\alpha r n$ as a subset of vertices of a $(k,m)$-orchard with $km\leq 2 \alpha n$. Hence if $|J|\geq \ell/4r$, we have that 
\[|Y_J^{i^*}|\geq \frac{\ell}{4r}\left(\frac{\sqrt{\alpha}n}{\ell}\right)-4 \alpha r n \geq \left(\frac{\sqrt{\alpha}}{4r}-4\alpha r\right) n \geq  2\alpha n, \]
using conditions \ref{cond: Y size} and \ref{cond: external disjoint} of Lemma~\ref{lem: first round} for $i=i^*$ and our upper bound on $\alpha$. 

\vspace{2mm}

We now define a random subset $\tilde{Y}^{i^*}$ by taking each vertex $v\in Y_{[\ell]}^{i^*}$ into $\tilde{Y}^{i^*}$ independently with probability $q:=\frac{\ell m}{2\sqrt{\alpha} n}$, noting that $0<q\leq  1$ due to the fact that $\ell m\leq km \leq 2 \alpha n \leq 2\sqrt{\alpha} n$.  
For $j\in[\ell]$, we define $\tilde{Y}_{i^*j}:=\tilde{Y}^{i^*}\cap Y_{i^*j}$. It follows from \ref{cond: Y size} that $\EE[|\tilde{Y}_{i^*j}|]= q|Y_{i^*j}| \leq m/2$ for all $j\in [\ell]$ and an application of Theorem \ref{thm:chernoff} as well as a union bound gives that with high probability, $|\tilde{Y}_{i^*j}|\leq m$ for all $j\in [\ell]$.  Note that in order to show that the upper bound on the failure probability given by Theorem \ref{thm:chernoff} is strong enough to beat  a union bound of the  $\ell$ events, we use our lower bound on $m$ and Fact~\ref{fact:dense}.  Furthermore, we have that with high probability, for any 
 choice of $J\subseteq [\ell]$ with $|J|\geq \ell/4r$ and any choice of  $(i_0,j_0) \in [t]\times [\ell]$ and subsets $ I_{1},\ldots, I_{r-2}\subset [t]\times [\ell]$ with $|I_{i'}|\geq pk$ for $i'\in [r-2]$ and $|I_{r-2}|\geq \gamma k$, there   exists a copy of $K_r$ traversing $\cup_{j\in J}\tilde{Y}_{i^*j}$,   $Z_{i_0j_0}$ and the sets  $Z(I_{i'})$ for $i'\in[r-2]$.
Indeed, this follows from an application of Theorem \ref{thm:chernoff} very similar to the proof of Lemma \ref{lem:diamond trees with many transversal triangles}. 
For a fixed $J$, $(i_0,j_0)$ and  $I_{i'}$ for $i'\in[r-2]$ as above, we have that there is some subset $X$ of  $\alpha n$ vertices of $Y_J^{i^*}$ such that each vertex in $X$ has a copy of $K_{r-1}$ in its neighbourhood which traverses the sets $Z_{i_0j_0}$ and  $Z(I_{i'})$ for $i'\in[r-2]$.  
Indeed, $X$ can be found by repeated applications of Corollary \ref{cor:transversalcliques} \eqref{cor:r-2small}, deleting vertices from $Y_J^{i^*}$ and adding them to $X$ on each application. 
Therefore $\EE[|X\cap \tilde{Y}^{i^*}|]=q|X|=\sqrt{\alpha} \ell m/2$ and by Theorem \ref{thm:chernoff}, the probability that $X\cap \tilde{Y}^{i^*}=\emptyset$ is less than $e^{-\sqrt{\alpha}\ell m/4}$.  
Due to the fact that $\ell m\geq m \geq p^{r-1}n\geq n^{(r-2)/(2r-3)}$ because of  our lower bound on $m$ and Fact~\ref{fact:dense}, we have that this probability tends to $0$ as $n$ tends to infinity. Moreover as there are less than
$k \cdot (2^k)^{r-2} \cdot 2^{\ell}\leq 2^{r k}$ choices for such a  $(i_0,j_0)$, $I_{i'}$ for $i'\in[r-2]$ and $J$ we have that a union bound  gives the traversing copies of $K_r$ for all choices with high probability. Indeed, we have that 
\[rk\leq \frac{2\alpha r n}{m}\leq \frac{2\alpha r}{p^{r-1}}\leq 2\alpha r \ell m\leq  \frac{\sqrt{\alpha} \ell m}{8},\]
using our conditions on $km$, $m$, $ \ell m $ and our upper bound on $\alpha$ from the hypotheses. 

Therefore, we can fix an instance of $Y^{i^*}$ which satisfies the desired conclusions that  we have shown happen with high probability. For each $j\in[\ell]$, taking $\tilde{Y}_{i^*j}=\tilde{Y}^{i^*}\cap Y_{i^*j}$ and defining $\tilde{\Upsilon}_{i^*j}:=\{S_v\in \Upsilon_{i^*j}:v\in \tilde{Y}_{i^*j}\}$, we apply condition \ref{cond: finding diamond trees} for $i=i^*$ to get a family $\bc{Q}_{i^*}$ of $\cD_{i^*j}=(T_{i^*j},R_{i^*j}, \Sigma_{i^*j})$ for $j\in[\ell]$ such that for each $j$ we have that $R_{i^*j}=Z_{i^*j}\cup \tilde{Y}_{i^*j}$ and $\Sigma_{i^*j}=\Pi_{i^*j}\cup \tilde{\Upsilon}_{i^*j}$.  This completes the step for $i^*$ and we move to $i^*+1$ and repeat. Doing this for each $1 \leq i^* \leq t$ completes the proof.  
\end{proof}

\subsection{The existence of shrinkable orchards of smaller order} \label{sec:slightly smaller}

With Lemma \ref{lem:orchard with dense triangle hypergraph} in hand, we can now complete the proof of Proposition~\ref{prop:large shrinkable} as follows.

\begin{proof}[Proof of Proposition~\ref{prop:large shrinkable}]
Fix $\eps>0$ small enough to apply Lemmas~\ref{lem:density condition for shrinkable} and \ref{lem:orchard with dense triangle hypergraph} and Proposition~\ref{prop:easy large order shrinkable} all with the same  $\alpha$ and $\gamma$ and small enough that to gurantee that $p\geq Cn^{-1/(2r-3)}$ with $C:=(2/\alpha)^{1/r}$ (see Fact~\ref{fact:dense}). We also guarantee that $\eps>0$ is  small enough to ensure that $n$ is sufficiently large in what follows. 

Now  note that Proposition~\ref{prop:easy large order shrinkable} directly implies the existence of  the desired shrinkable orchard if $p^{1-r} \leq p^{r-1}n$ or if $ p^{r-1} n\leq p^{1-r} \leq m\leq n^{7/8}$ and so we can assume from now on that 
\begin{equation} \label{eq: m and p conds} p^{r-1}n\leq m \leq \min\{p^{1-r},n^{7/8}\} .
\end{equation}
We are therefore in a position (due to our lower bound on $m$) to apply Lemma~\ref{lem:orchard with dense triangle hypergraph} but we first need to fix $k,t,\ell\in\NN$ so that the conditions \eqref{eq: k t l m conds} are satisfied. We first fix $\ell\in \NN$ so that $p^{1-r}\leq \ell m\leq 2p^{1-r}$.  This is possible as $m\leq p^{1-r}$ and so there is a multiple of $m$ in the desired range. Next we fix $t\in\NN$ to be any integer such that $\alpha n\leq t\ell m\leq 2\alpha n$. Again, this is possible because $\ell m\leq 2p^{1-r}\leq \alpha n^{(r-1)/(2r-3)}\leq \alpha n$ using Fact \ref{fact:dense}. So there is indeed a multiple $t$ of $\ell m$ in the desired range. Finally, we fix $k=t\ell$ and so we have that all the conditions in \eqref{eq: k t l m conds} are satisfied with our choice of parameters. Before analysing the conclusion of Lemma~\ref{lem:orchard with dense triangle hypergraph}, let us point out a few further implications of our choices of parameters. Firstly, we have that 
\begin{equation}
\label{eq: p upper}
k^{-r^3\gamma}\geq \left(\frac{1}{k}\right)^{\frac{1}{2^r}}\geq \left( \frac{m}{n}\right)^{\frac{1}{2^r}}
\geq  \left( p^{r-1}\right)^{\frac{1}{2^r}} \geq p.\end{equation}
Moreover
\begin{equation} \label{eq: p ell compare}
\frac{\ell}{k} \leq \frac{2}{p^{r-1}mk}\leq  \frac{2}{\alpha p^{r-1}n}\leq p, 
\end{equation}
where we use the upper bound on $\ell m$ in the first inequality, the lower bound on $km$ in the second inequality and the fact that $p\geq C n^{-1/(2r-3)}\geq Cn^{-1/r}$ from Fact~\ref{fact:dense} in the last inequality. Putting \eqref{eq: p upper} and \eqref{eq: p ell compare} together then gives that 
\begin{equation} \label{eq: useful k ell bounds}
k^{1-r^3{\gamma}}\geq pk \geq \ell. 
    \end{equation}

Now we apply Lemma~\ref{lem:orchard with dense triangle hypergraph} and let $\bc{O}$ be the resulting $(k,m)_r$-orchard partitioned into $(\ell,m)_r$-suborchards $\bc{Q}_1,\ldots,\bc{Q}_t$. We will show that $\bc{O}$ is $\gamma$-shrinkable by appealing to Lemma~\ref{lem:density condition for shrinkable}. Firstly note that due to the upper bound of $m\leq n^{7/8}$ and the fact that $k=\Theta(n/m)$, by forcing $n$ to be sufficiently large, we can assume that $k$ is also sufficiently large to apply Lemma~\ref{lem:density condition for shrinkable}. We therefore just need to check the density condition \eqref{density cond} of the lemma. So fix some arbitrary $\cD\in \bc{O}$ and suborchard $\bc{P}\subset \bc{O}\setminus \{\cD\}$ such that $\bc{P}\geq k/(4r)$. By pigeonhole principle, there exists an $i\in [t]$ such that $|\bc{P}\cap \bc{Q}_i|\geq \ell/4r$. So fix such an $i$ and define $\bc{P}^*:=\bc{P}\cap \bc{Q}_i$, noting that we have that $|\bc{P}^*|\leq k^{1-r^3{\gamma}}$ due to \eqref{eq: useful k ell bounds}. Now we simply need to check that for any choice of suborchards $\bc{O}_1,\ldots,\bc{O}_{r-2}\subset \bc{P}\setminus \bc{P}^*$, with $|\bc{O}_{i'}|\geq k^{1-r^3{\gamma}}$ for $i'\in[r-3]$ and $|\bc{O}_{r-2}|\geq \gamma k$, there is a copy of $K_r$ in $G$ traversing $R_{\cD}$, $R(\bc{P}^*)$ and the sets $R(\bc{O}_{i'})$ for $i'\in[r-2]$.  This is verified by the conclusion of Lemma~\ref{lem:orchard with dense triangle hypergraph},  setting $\cD'=\cD$, $\bc{Q}'=\bc{P}^*$ and $\bc{O}'_{i'}=\bc{O}_{i'}$ for $i'\in[r-2]$, noting that the lower bounds on the sizes of the $\bc{O}'_{i'}$ are guaranteed by \eqref{eq: useful k ell bounds}. Hence $\bc{O}$ is indeed $\gamma$-shrinkable by Lemma~\ref{lem:density condition for shrinkable} and this concludes the proof. 
\end{proof}

\section{The final absorption} \label{sec:absorbingstructure}
The aim of this section is to prove Proposition \ref{prop:finalabsorption}  which we restate here for convenience.

\theoremstyle{plain}
\newtheorem*{prop:finalabsorption}{Proposition \ref{prop:finalabsorption}}
\begin{prop:finalabsorption}{(Restated)}
For any $3\leq r\in \NN$ and $0<\alpha,\eta<1/2^{3r}$ there exists an $\eps>0$ such that the following holds for any $n$ vertex $(p,\beta)$-bijumbled graph $G$ with $\beta\leq \eps p^{r-1}n$ 
and any vertex subset $W\subseteq V(G)$ with $|W|\geq n/2$.

 There exist vertex subsets $A,B\subset V$ such that $A\subset W$,  $|A|\leq \alpha n$, $|B|\leq \eta p^{2r-4} n$ and for  any $(k,m)_r$-orchard $ \bc{R}$  
whose vertices lie in $V(G)\setminus (A\cup B)$, we have that if $|A|+|V(\bc{R})|\in r \mathbb{N}$, $k\leq \alpha^2 n^{1/8}$ and $m\geq n^{7/8}$ then $G[A\cup V(\bc{R})]$ has a $K_r$-factor.  
\end{prop:finalabsorption}

In order to prove this, in Section \ref{sec:defabsorbingstructure} we first define an absorbing structure whose vertex set will play the role of $A$ in Proposition~\ref{prop:finalabsorption}. We then prove that it has the required absorbing property. Next, in Section \ref{sec:findingabsorbingstruc},  we prove that we can find the absorbing structure in a suitably pseudorandom graph  and show that  this implies Proposition \ref{prop:finalabsorption}.

\subsection{Defining an absorbing structure} \label{sec:defabsorbingstructure}
Recall from Section~\ref{sec:templates}  the definition of a template and the fact that templates of flexibility~$t$ with maximum degree~$40$ exist for all large enough~$t$ (Theorem~\ref{thm:montyexist}). We will use a template as an auxiliary graph to define an \emph{absorbing structure} which can contribute to a~$K_r$-factor in many ways.

\begin{dfn} \label{def:absorbingstruc}
Let $\pzc{T}$ be a template with flexibility $t$ on vertex sets $I$ and $J:=J_1 \cup J_2$ with $|I|=3t$ and $|J_1|=|J_2|=2t$. A $K_r$-\emph{absorbing structure} $\AA$ of \emph{order} $M$  with respect to $\pzc{T}$ in $G$ consists of a labelled matching of $(r-1)$-cliques in $G$, ~$\Xi(\AA):=\{S_i:i\in I\}\subset K_{r-1}(G)$ and a labelled  $(4t,M)_r$-orchard~$\bc{J}(\AA):=\{\cD_j:j\in J\}$ such that the following holds for each $i\in I$ and $j\in J$: 
\begin{itemize}
    \item $S_i\cap V(\cD_j)=\emptyset$;
    \item if $ij\in E(\pzc{T})$ then there is a vertex in the removable set of vertices $R_j$ of $\cD_j$ which forms a $K_r$ with $S_i$ in $G$. 
\end{itemize}
We say that $\AA$ has \emph{flexibility} $t$, which is inherited from the template by which $\AA$ is defined.  We refer to the vertices of the absorbing structure, denoted $V(\AA)$, which is all vertices which feature in cliques in $\Xi(\AA)$ or diamond trees in $\bc{J}(\AA)$.
\end{dfn}

\begin{figure}[h]
    \centering
  \includegraphics[scale=0.84]{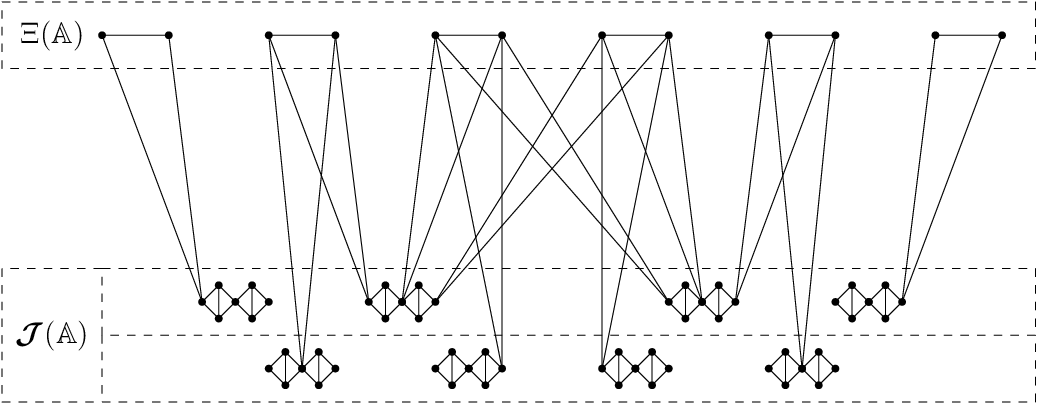}
    \caption{   \label{fig:absorbingstruc} A $K_3$-absorbing structure of order $3$ and flexibility $2$, whose defining template is the template $\pzc{T}$ displayed in Figure~\ref{fig:template}. } 
  \end{figure}

See Figure~\ref{fig:absorbingstruc} for an example of an absorbing structure. Note that a $K_r$-absorbing structure $\AA$ of flexibility $t$ and order $M$ has less than \begin{equation} \label{eq:absorbingstructuresize} 3t(r-1)+4t((2M-1)r+1)\leq 8rtM-rt+t\leq 8rtM\end{equation} vertices. The absorbing structure is defined in such a way that it inherits the robust property that the template has with respect to perfect matchings but has such a property with respect to $K_r$-factors. The following lemma demonstrates this and reduces Proposition \ref{prop:finalabsorption} to finding an appropriate absorbing structure in $G$. 

\begin{lem} \label{lem:reduction}
For any $3\leq r\in \NN$ and $0<\zeta,\eta<1/2^{2r}$ there exists an $\eps>0$ such that the following holds for any $n$ vertex $(p,\beta)$-bijumbled graph $G$ with $\beta\leq \eps p^{r-1}n$.
 Suppose that~$t,M\in \mathbb{N}$ such that $ tM \geq \zeta n$ and there exists an $K_r$-absorbing structure~$\AA$ of flexibility $t$ (with respect to some template $\pzc{T}$) and order $M$ in $G$.  Let $A:=V(\AA)$.

Then there exists some vertex subset $B\subset V(G)$, such that
$|B|\leq \eta p^{2r-4} n$ and for  any $(k,m)_r$-orchard $ \bc{R}$  
whose vertices lie in $V(G)\setminus (A\cup B)$, we have that if $|A|+|V(\bc{R})|\in r \mathbb{N}$, $k\leq \frac{t}{4r}$ and $m\geq M$ then $G[A\cup V(\bc{R})]$ has a  $K_r$-factor.

   \end{lem} 
\begin{proof}
Fix $\eps>0$ small enough to apply Lemma~\ref{lem:absorptionbetweenlayers} with $\zeta$ and $\eta$ as defined here 
and small enough to apply Corollary~\ref{cor:transversalcliques} with $\alpha_{\ref{cor:transversalcliques}}=\alpha:=\zeta/(r+1)$. 
Let $\bc{O}=\{\cD_j:j\in J_2\}$ be the suborchard of $\bc{J}(\AA)$ defined by those indices which lie in the flexible set $J_2$ of the template $\pzc{T}$ which defines $\AA$. Thus $\bc{O}$ is a $(2t,M)_r$-orchard. 
Therefore, applying Lemma \ref{lem:absorptionbetweenlayers}, we have that  there exists a set $B\subset V(G)$ with $|B|\leq \eta p^{2r-4} n$ and  for any $(k,m)_r$-orchard $\bc{R}$ as in the statement of the lemma, $\bc{O}$ absorbs $\bc{R}$. Indeed, note that in the notation of Lemma~\ref{lem:absorptionbetweenlayers}, we have that $k,m$ and $M$ are as defined here while $K=2t$. Hence the condition that $k\leq K/(8r)$ is precisely the same as our presumption that $k\leq t/(4r)$ whilst the condition that $kM\leq mK$ is guaranteed by the fact that $m\geq M$ and $k\leq K/(8r)$.  Unpacking the conclusion of Lemma \ref{lem:absorptionbetweenlayers}, we thus have that for any such $\bc{R}$  there exists some subfamily $\bc{P}_1\subseteq \bc{O}$ such that $|\bc{P}_1|=(r-1)k\leq t/2$ and $G[V(\bc{P}_1)\cup V(\bc{R})]$ has a $K_r$-factor. We will show that $G[A\setminus V(\bc{P}_1)]$ also has a $K_r$-factor which will complete the proof. 

Now note that for any $\bc{P}\subset \bc{O}$  such that $|\bc{P}|=t$, we have that $G[A\setminus V(\bc{P})]$ has a $K_r$-factor. Indeed let $\overline{J}:=\{j\in J_2: \cD_j\in \bc{P}\}$ be the indices of diamond trees that feature in the suborchard $\bc{P}$. By the definition of the template $\pzc{T}$, Definition \ref{def:template}, we know that there is a perfect matching $F\subset E(\pzc{T})$ in $\pzc{T}[V(\pzc{T})\setminus \overline{J}]$.
Now for $ij\in F$, we can take a $K_r$-factor on $S_i\cup V(\cD_j)$ in $G$ guaranteed by the fact that $S_i$ forms a copy of $K_r$ with a removable vertex of $\cD_j$ (Definition \ref{def:absorbingstruc}) and the key property of the removable vertices of a $K_r$-diamond tree (Observation \ref{obs:diamondtree}).  
As $F$ is a perfect matching in $\pzc{T}[V(\pzc{T})\setminus \overline{J}]$, we see that by taking these $K_r$-factors for each $ij\in F$, we obtain a $K_r$-factor in $G[A\setminus V(\bc{P})]$ as required.

If we had that $|\bc{P}_1|=t$, this would complete the proof. However we have that $\bc{P}_1$ is actually much smaller than this. Indeed $|\bc{P}_1|\leq t/2$. We will proceed by finding some $\bc{P}_2\subset \bc{O}\setminus \bc{P}_1$ such that $G[V(\bc{P}_2)]$ has a $K_r$-factor and $|\bc{P}_1|+|\bc{P}_2|=t$.
We  build $\bc{P}_2$ by the following greedy process. We initiate with $\bc{O}'=\bc{O}\setminus \bc{P}_1$ and $\bc{P}_2=\emptyset$. Then at each time step, as long as $|\bc{P}_1|+|\bc{P}_2|+r\leq t$ we partition $\bc{O}'$ into $r$ parts  $\bc{O}'=\bc{O}_1 \cup \cdots \cup\bc{O}_r$ such that the sizes of the parts are as equal as possible. We let $R_x$ be the union of the removable vertices of diamond trees in the orchard $\bc{O}_x$ for $x\in[r]$. We have that each $R_x$ has size at least $tM/(r+1) \geq \zeta n/(r+1)=\alpha n $. Therefore,  by Corollary \ref{cor:transversalcliques} \eqref{cor:r-2small},  
there is a copy of $K_r$ traversing the  $R_x$, $x\in[r]$. This gives some $r$-tuple of diamond trees $\cD_1,\ldots, \cD_r$ such that $\cD_x \in \bc{O}_x$ for all $x\in[r]$ and there is a $K_r$-factor in $G[V(\cD_1)\cup \cdots \cup V(\cD_r)]$, given by taking the copy of $K_r$ that traverses their sets of removable vertices and applying Observation \ref{obs:diamondtree}. 
We add $\cD_1,\ldots,\cD_r$ to $\bc{P}_2$ and delete these diamond trees from the orchard $\bc{O}'$ which completes this time step. Clearly at all points in this process there is a $K_r$-factor in $G[V(\bc{P}_2)]$ and we claim that this process terminates when $|\bc{P}_1|+|\bc{P}_2|$ is exactly equal to $t$.  
Indeed if this is not the case, as we increase the size of $|\bc{P}_2|$ by exactly $r$ in each step, we have that $|\bc{P}_1|+|\bc{P}_2|=t-s$ for some $s\in [r-1]$ at the end of the process. Let $\bc{P}_3 \subset \bc{O} \setminus (\bc{P}_1\cup \bc{P}_2)$ be a set of $s$ $K_r$-diamond trees. Now as $V_1:=V(\bc{R})\cup V(\bc{P}_1) \cup V(\bc{P}_2)$ hosts a $K_r$-factor, we have that $r||V_1|$. Likewise, we know from above that $V_2=A\setminus (\cup_{i=1}^3V(\bc{P}_i))$ hosts a $K_r$-factor and so $r||V_2|$. Due to the fact that $r$ divides the size of $A\cup V(\bc{R})$ and $A\cup V(\bc{R})=V_1 \cup V_2 \cup V(\bc{P}_3)$, we can defer that $r||V(\bc{P}_3)|$. This is a contradiction as $\bc{P}_3$ is a set of $s$ $K_r$-diamond trees for some $1\leq s \leq r-1$ and the number of vertices in any diamond tree is $1 \mod r$. Therefore we can find a $\bc{P}_2$ as claimed.

Finally, taking $\bc{P}:=\bc{P}_1 \cup \bc{P}_2$, we are then done by taking our $K_r$-factor in $G[A\cup V(\bc{R})]$ to be the union of the $K_r$-factor in $G[V(\bc{R}) \cup V(\bc{P}_1)]$, the $K_r$-factor in $G[V(\bc{P}_2)]$ and the $K_r$-factor in $G[A \setminus V(\bc{P})]$. 
\end{proof}

\subsection{Finding an absorbing structure}
\label{sec:findingabsorbingstruc}

Lemma \ref{lem:reduction} reduces Proposition \ref{prop:finalabsorption} to finding an appropriate absorbing structure in $G$. In this section we prove that this is possible by proving the following proposition.

\begin{prop} \label{prop:findingabsorbingstruc}
For any $3\leq r\in \NN$ and $0<\alpha<1/2^{3r}$ there exists an $\eps>0$ such that the following holds for any $n$ vertex $(p,\beta)$-bijumbled graph $G$ with $\beta\leq \eps p^{r-1}n$ 
and any vertex subset $W\subseteq V(G)$ with $|W|\geq n/2$. 
There exists an $K_r$-absorbing structure $\AA$ in $G$ of flexibility $t=\alpha n^{1/8}$ and order $M=n^{7/8}$ such that $V(\AA)\subseteq W$.
\end{prop}

With Lemma \ref{lem:reduction} and Proposition \ref{prop:findingabsorbingstruc}, the proof of Proposition \ref{prop:finalabsorption} follows readily as we now show.

\begin{proof}[Proof of Proposition \ref{prop:finalabsorption}]
Fix $\zeta:=\alpha/(8r)$ and $\eps>0$ small enough to apply Lemma \ref{lem:reduction} with $\zeta$ and $\eta$ as defined here and 
small enough to apply Proposition~\ref{prop:findingabsorbingstruc} with $\alpha_{\ref{prop:findingabsorbingstruc}}=\zeta$.  We can therefore apply Proposition \ref{prop:findingabsorbingstruc} to get an absorbing structure $\AA$ in $G$ with flexibility $t=\zeta n^{1/8}$ and order $M=n^{7/8}$. We have that $A=V(\AA)\subset W$ has size $|A|\leq 8rtM=\alpha n$  (see \eqref{eq:absorbingstructuresize}). The conclusion then follows from Lemma \ref{lem:reduction} noting that $k\leq \alpha^2n^{1/8}$ implies that $k\leq t/(4r)=\zeta n^{1/8}/(4r)=\alpha n^{1/8}/(48r^2)$ due to our upper bound on $\alpha$.   
\end{proof}

Now in order to prove the existence of an absorbing structure as in \Cref{prop:findingabsorbingstruc}, we will first fix some template $\pzc{T}$ which will define $\AA$. Next, we will set aside a large matching $\Pi\subset K_{r-1}(G)$ of $(r-1)$-cliques. These will be candidates for the   matching of $(r-1)$-cliques $\Xi(\AA)$ in our absorbing structure but we start with a much bigger set $\Pi$ of size $\Omega(n^{2/3})$. Moreover, to each $(r-1)$-clique~$S\in \Pi$ we will associate a set of vertices $X_S\subset G$ such that $X_S\subset N^G(S)$, $|X_S|=\Omega(n^{1/3})$ and, crucially, the sets $\{X_S:S\in \Pi\}$ are \emph{disjoint}. We will find this matching of $(r-1)$-cliques $\Pi$ with a simple greedy procedure, appealing to Corollary~\ref{cor:transversalcliques}~\eqref{cor:popular K_r-1} to find each $S\in \Pi$ (and its corresponding neighbourhood set $X_S$), one by one. 

After finding $\Pi$, we then turn to constructing the $(4t,M)$-orchard $\bc{J}(\AA)$ for the absorbing structure $\AA$. Again, this will be done greedily, fixing the diamond trees $\cD\in \bc{J}(\AA)$ one at a time. Let us consider fixing some diamond tree $\cD_j\in \bc{J}(\AA)$. Note that as we fix $\cD_j$, we immediately get restrictions on which $S\in \Pi$ remain as candidates to play the r\^ole of certain~$S_i\in \Xi(\AA)$. Indeed, if the removable vertices of $\cD_j$ are disjoint from $N^G(S)$ and $ij$ is an edge in the template~$\pzc{T}$ defining $\AA$, then there is no way $S$ can play the role of $S_i$ in $ \Xi(\AA)$. Therefore as we fix our diamond trees, we will aim to have that their sets of removable vertices intersect as many of the~$X_S$ (and hence neighbourhoods $N^G(S)$) for $S\in \Pi$, as possible. 

In order to do this, we will use the following lemma, which shows that we can find diamond trees whose removable vertices intersect many prescribed sets (in our case this will be the sets $X_S$). The proof of this lemma is a simple application of~\cref{prop: choosing removable vertices}.  

\begin{lem} \label{lem:diamond tree that intersects many sets}
For any $3\leq r\in \NN$ and  $0<\alpha<1/2^{2r}$,  there exists an $\eps>0$ such that the following holds for any $n$ vertex $(p,\beta)$-bijumbled graph $G$ with $\beta\leq \eps p^{r-1}n$.

Suppose $\frac{\alpha^2}{2} n^{2/3}\leq \ell\leq \alpha n^{2/3}$  and we have disjoint vertex subsets $W,U_1,\ldots,U_\ell$ such that $|W|\geq n/4$ and $|U_i|\geq  n^{1/3}$ for all $i\in[\ell]$. Then there exists a diamond tree $\cD=(T,R,\Sigma)$ in $G$ such that the following conditions hold:

\begin{enumerate}[label=(\roman*)]
    \item \label{cond:i}$\Sigma\subset K_{r-1}(G[W])$ is a matching of $(r-1)$-cliques in $W$;
    \item \label{cond:ii} $R\subset \cup_{i=1}^{\ell}U_i$ and $R$ intersects $\ell'$ of the sets $U_i$ for some $\ell'\geq \ell/(4r)$;
    \item \label{cond:iii} The order of $\cD$ is at most $n^{2/3}$;
    \item \label{cond:iv} For all but at most $n^{1/2}$ of the indices $i\in[\ell]$, we have that $|V(\cD)\cap U_i|\leq n^{1/6}$.
    \end{enumerate}

\end{lem}

\begin{proof}

We begin by fixing $\gamma:=\ell/n^{2/3}$ so that $\frac{\alpha^2}{2}\leq\gamma\leq \alpha$ and we fix $\eps$ small enough to apply Propsition~\ref{prop: choosing removable vertices} with $\alpha_{\ref{prop: choosing removable vertices}}=\alpha':=\gamma/(4r)$ and small enough to guarantee that $p\geq Cn^{-1/(2r-3)}$ with $C=4/\alpha'$ (see Fact~\ref{fact:dense}). Now shrink each set $U_i$ so that it has exactly $ n^{1/3}$ vertices and define $U:=\cup_{i=1}^{\ell} U_i$. Furthermore fix $d_*:=\alpha'^2p^{r-1}n$ and apply Proposition~\ref{prop: choosing removable vertices} with $U$, $W$ and $z=\alpha'n$. So we get disjoint subsets $X,Y\subset U$ as in the outcome of Proposition~\ref{prop: choosing removable vertices}. 

Now firstly note that as $|X|+|Y|=z=\alpha'n$, $|U|= \ell n^{1/3}=\gamma n=4rz$ and each of the $U_i$ have equal size, we must have that $X\cup Y$ intersects at least $\ell/4r$ of the sets $U_i$. We will choose our $\cD=(T,R,\Sigma)$ so that $R$ intersects all the sets $U_i$ that $X\cup Y$ intersects, thus guaranteeing condition \ref{cond:ii}. Indeed, if we let $Y'\subset Y$ be the \emph{minimal} subset of $Y$ such that there exists no $i\in [\ell]$ with  $Y\cap U_i\neq \emptyset$ and $Y'\cap U_i=\emptyset$, Proposition~\ref{prop: choosing removable vertices} gives the existence of a diamond tree $\cD=(T,R,\Sigma)$ so that $R=X\cup Y'$ and $\Sigma\subset K_{r-1}(G[W])$  and so conditions \ref{cond:i} and \ref{cond:ii} are satisfied. 

In order to establish condition \ref{cond:iii}, note that $|Y'|\leq \ell \leq \frac{n^{2/3}}{2}$ and  if $|X|>1$ then 
\[|X|\leq \frac{2z}{d_*}\leq \frac{2}{\alpha'p^{r-1}}\leq \frac{2n^{(r-1)/(2r-3)}}{\alpha'C^{r-1}}\leq \frac{n^{2/3}}{2},\]
due to our definition of $C$ and the fact that $\frac{r-1}{2r-3}\leq \frac{2}{3}$ for all $r\geq 3$.
 Finally,  
condition \ref{cond:iv} is a simple consequence of \ref{cond:iii}. Indeed, if \ref{cond:iv} were not true, then as the $U_i$ are pairwise disjoint, we would have that $\cD$ has order greater than $n^{1/2}\cdot n^{1/6}\geq n^{2/3}$, a contradiction. 
\end{proof}

Let us return to sketching the proof of \Cref{prop:findingabsorbingstruc}, considering now that we can use Lemma~\ref{lem:diamond tree that intersects many sets} to find diamond trees $\cD\in \bc{J}(\AA)$. As discussed above, the key property of diamond trees generated by \Cref{lem:diamond tree that intersects many sets} is \ref{cond:ii}, allowing us to find diamond trees that intersect many of sets $\{X_S:S\in \Pi\}$ which we begin the proof with. The property \ref{cond:iv} will also be useful as it shows that in the process of building $\bc{J}(\AA)$ one by one, we do not destroy many of the sets~$X_S$ and most of them remain large and can be used by other $\cD\in \bc{J}(\AA)$. 

One potentially troublesome consequence of Lemma~\ref{lem:diamond tree that intersects many sets} is that the diamond trees it finds are far too small \ref{cond:iii}. Indeed the diamond trees in our orchard $\bc{J}(\AA)$ are supposed to be of order~$M=n^{7/8}$. It turns out that this is not such a big hurdle as we can find a large diamond tree disjoint from all the $X_S$ and connect it to the diamond tree $\cC$ output by~\Cref{lem:diamond tree that intersects many sets}. In more detail, we can apply Proposition~\ref{prop: choosing removable vertices} to create a large (linear) pool $Y$ of vertices that can be removable vertices of some diamond tree which will be disjoint from all the vertices in the sets~$X_S$. We also consider the large (linear) pool of vertices $Z$ that lie in  some $X_S\setminus V(\cC)$ with~$S\in \Pi$ such that the removable vertices of $\cC$ intersect $X_S$. It is not hard to show (see for example Corollary~\ref{cor:transversalcliques}~\eqref{cor:popular K_r-1}) that there is a copy of $K^-_{r+1}$ with one degree $r-1$ vertex in $Y$ and the other in $X_{S^*}\subset Z$ for some $S^*\in \Pi$. By also taking $S^*$ into $\cD$ and choosing an appropriate~$Y'\subset Y$ to apply the key property of Proposition~\ref{prop: choosing removable vertices}, we can obtain a diamond tree $\cD$ of the correct size that contains the diamond tree $\cC$ output by Lemma~\ref{lem:diamond tree that intersects many sets}.

More troublesome is the fact that the condition~\ref{cond:ii} which gives that the removable vertices of~$\cC$ intersect many of the desired sets $X_S$ is, in fact, not strong enough. Indeed, consider some fixed~$i\in I$ for which we want to find a copy $S_i$ of $K_{r-1}$  to lie in $\Xi(\AA)$.   If $j,j'\in N^{\pzc{T}}(i)$ and the sets~$\{X_S:S\in \Pi,R_{\cD_{j}}\cap X_S\neq \emptyset\}$ and $\{X_S:S\in \Pi, R_{\cD_{j'}}\cap X_S\neq \emptyset\}$ (here, as usual, we use $R_\cD$ to denote the removable vertices of $\cD$) are disjoint, then already there are no candidates for $S_i$ in $\Pi$. To fix this, we actually need that when we choose a diamond tree $\cD\in \bc{J}(\AA)$, we want~$R_\cD$ to intersect \emph{almost all} of the sets $\{X_S:S\in \Pi\}$. We achieve this by iterating Lemma~\ref{lem:diamond tree that intersects many sets}, creating constantly many disjoint diamond trees $\cC$ that together hit almost all of the $X_S$ with their removable vertices. We then connect all of these diamond trees $\cC$ with a large diamond tree disjoint from the sets $X_S$ to obtain the desired diamond tree $\cD\in \bc{J}(\AA)$. This connecting process is similar (although slightly more involved) to the connecting process outlined in the previous paragraph. We now give the full details for the proof of Proposition~\ref{prop:findingabsorbingstruc}, concluding this section and chapter.

\begin{proof}[Proof of Proposition~\ref{prop:findingabsorbingstruc}] 
We begin by fixing $\eps>0$ small enough to apply Corollary~\ref{cor:transversalcliques} with $\alpha_{\ref{cor:transversalcliques}}=\alpha':=\alpha^2/16r$ and to apply Proposition~\ref{prop: choosing removable vertices} and Lemma~\ref{lem:diamond tree that intersects many sets} each with $\alpha_{\ref{prop: choosing removable vertices}}=\alpha_{\ref{lem:diamond tree that intersects many sets}}=\alpha$. We also make sure that $\eps>0$ is small enough to force $n$ to be sufficiently large
 in what follows and small enough to guarantee that $p\geq C'n^{-1/(2r-3)}$ for $C':=2/\alpha'^2$ using Fact~\ref{fact:dense}.  We further fix some template~$\pzc{T}$ with vertex sets~$I$ and~$J=J_1 \cup J_2$ of flexibility~$t$ and maximum degree 40 which we know exists for~$n$ (and hence~$t$) sufficiently large due to Montgomery~\cite{M14a} (see~\Cref{thm:montyexist}). 
 
 We will find an absorbing structure with respect to $\pzc{T}$ and so must prove the existence of a matching  $\Xi(\AA)=\{S_{i}:i \in I\}\subset K_{r-1}(G[W])$  of $3t$ copies of $K_{r-1}$, and a $(4t,M)$-orchard  $\bc{J}=\bc{J}(\AA)=\{\cD_j:j\in J\}$ such that the conditions of Definition \ref{def:absorbingstruc} are satisfied.   We will do this in three stages. 
In Claim \ref{clm:candidateedges}, we fix some large matching $\Pi\subset K_{r-1}(G[W])$ of $(r-1)$-cliques which will be candidates for the $(r-1)$-cliques which will feature in $\Xi(\AA)$. We will guarantee that the cliques in $\Pi$ are contained in many copies of $K_r$ which will help as we proceed to build our absorbing structure. In Claim \ref{clm:fixingdiamondtrees}, we will fix the $K_r$-diamond trees which will form our orchard $\bc{J}$ for our $K_r$-absorbing structure. We will carefully control how these diamond trees intersect the cliques in our candidate set $\Pi$ and their  neighbourhoods. Finally, we will show that we can find a suitable $\Xi(\AA)\subset  \Pi$ so that we obtain the desired absorbing structure.

\begin{clm} \label{clm:candidateedges}
There exists a matching $ \Pi=\{S_1,\ldots,S_{\ell}\}\subset K_{r-1}(G[W])$of $\ell:=\alpha n^{2/3}$ copies of $K_{r-1}$  and sets $X_h\subset W \setminus V(\Pi)$ 
 for each $h \in[\ell]$, such that  the $X_h$ are pairwise disjoint, each has size $|X_h|=  2 n^{1/3}$ and for all $h\in[\ell]$ we have that $X_h\subset N^G_W(S_h)$. 
\end{clm}
\begin{claimproof}
We can do this by way of a simple greedy process choosing such an $(r-1)$-clique $S_h$ and set $X_h$ in order for $h=1,\ldots,\ell$. When choosing $S_h$ and $X_h$, we look at the set of vertices $V_h\subset W$ which have not been used in previous choices of $S_{h'}$ or $X_{h'}$. We  have that \[|V_h|\geq |W|-|\cup_{h'<h}(X_{h'}\cup S_{h'})|\geq \frac{n}{2}-(\ell-1)(r-1+2n^{1/3})\geq \left(\frac{1}{2}-2\alpha\right)n\geq  \frac{n}{4}\] and an application of Corollary \ref{cor:transversalcliques} \eqref{cor:popular K_r-1}  
with $W_0=W_1=W_2=V_h$
gives the desired $S_h$ and $X_h$ in $V_h$  using that  \[\alpha'^2p^{r-1}n\geq \alpha'^2 C' n^{1-(r-1)/(2r-3)}\geq 2n^{1/3}\]
  due to Fact \ref{fact:dense}. 
\end{claimproof}

\vspace{2mm}
Next we turn to fixing our $(4t,M)$-orchard $\bc{J}$.

\begin{clm} \label{clm:fixingdiamondtrees}
Let $S_h$ and $X_h$ for $h=1,\ldots, \ell $ be as in Claim \ref{clm:candidateedges}. Then there exists  a $(4t,M)$-orchard $\bc{J}=\{\cD_1,\ldots,\cD_{4t}\}$   such that $V(\bc{J})\subset W$ and the following properties hold for each $\cD_j=(T_j,R_j,\Sigma_j)$ with $j\in[4t]$:
\begin{enumerate}

    \item \label{cond:1} The set of removable vertices $R_j$ intersects at least $(1-\alpha)\ell$ of the sets $X_h$ with $h\in[\ell]$;

    \item \label{cond:2} $V(\cD_j)$ intersects at most $ C:=\frac{\log\left(\frac{2}{\alpha}\right)}{\log\left(\frac{4r}{4r-1}\right)}$ of the $S_h$ with $h\in[\ell]$;

\end{enumerate}
\end{clm}
Before proving the claim, let us see how it implies the proposition.
Indeed, taking the $(4t,M)_r$ orchard  $\bc{J}$ from Claim \ref{clm:fixingdiamondtrees} as $\bc{J}(\AA)$, we  just need to choose a matching of $(r-1)$-cliques $\Xi(\AA)=\{S_{i}:i\in I\}$  so that $S_{i}\cap V(\bc{J})=\emptyset$ for all $i\in I$ and whenever $ij\in E(\pzc{T})$, we have that there is a vertex in $R_j$ which forms a copy of $K_r$ with $S_{i}$. We  do this greedily, showing that for each $i=1,2,\ldots ,3t$ in order, there is a suitable choice for $S_{i}$ in $\Pi$.  We initiate by fixing $L\subseteq  [\ell]$ to be the indices $h\in[\ell]$ such that $S_h\cap V(\bc{J})=\emptyset$.  By condition \eqref{cond:2} in Claim \ref{clm:fixingdiamondtrees}, for large $n$ we have that 
\[ |L|\geq \ell-4Ct\geq (1-\alpha)\ell,\]
at the beginning of this process, recalling that $\ell=\alpha n^{2/3}$ and $t=\alpha n^{1/8}$. Now for $i=1,\ldots,3t$, we find an index $h=h(i)\in L$ such that $S_h$ forms a copy of $K_r$ with a vertex in $R_j$ for all $j$ such that $ij\in E(\pzc{T})$. We fix $S_{i}=S_h$ and delete $h$ from $L$. If this process succeeds in finding a suitable $h=h(i)$ for each $i\in I$ then the resulting $\Xi(\AA)=\{S_i:i\in I\}$ along with $\bc{J}$ form the desired $K_r$-absorbing structure.

It remains to check that we are successful at each step. So consider step $i^*\in[3t]$.  We have that $|L|\geq (1-\alpha)\ell-(i^*-1)\geq (1-2\alpha)\ell$ at the beginning of the step. Now for each $j\in J$ which is a neighbour of $i^*$ in the template $\pzc{T}$, we have by Claim \ref{clm:fixingdiamondtrees} \eqref{cond:1} that there are at most $\alpha \ell$ indices $h\in[\ell]$ such that no vertex of $R_j$  forms a $K_r$ with $S_h$ in $G$. Indeed for almost all choices of $h$, we have that $R_j\cap X_h\neq \emptyset$ and $X_h\subset N^G(S_h)$. Given that $\pzc{T}$ has maximum degree 40, we have that this gives at most $40\alpha \ell$ indices $h\in L$  that would not be a good choice for $h(i^*)$. Therefore there are at least $(1-42\alpha)\ell$ indices $h\in L$ which can be chosen as $h(i^*)$ and we simply choose one arbitrarily. 

This shows that the algorithm is successful in generating the desired absorbing structure and so it only remains to prove Claim~\ref{clm:fixingdiamondtrees}, which we do now.

\vspace{2mm}

\begin{altclaimproof}
We will find the diamond trees $\cD_j$, $j=1,\ldots,4t$ one by one so that they are vertex disjoint and satisfy the two conditions in the statement of the claim as well as the further following condition:
\begin{enumerate}\addtocounter{enumi}{2}
   \item \label{cond:3}$V(\cD_j)$ intersects all but at most $ Cn^{1/2}$ of the $X_h$ with $h\in[\ell]$ in more than $2Cn^{1/6}$ vertices.
 \end{enumerate}

We will initiate the process with $\Lambda:=[\ell]$ and  $U_h=X_h$ for all $h\in [\ell]$. These sets $U_h$ will keep track of vertices in $X_h$ that we are still allowed to use, that is, those vertices which have not been used in previously chosen diamond trees. Furthermore,  the set $\Lambda\subset [\ell]$ will keep track of all indices which are \emph{alive}. When we choose a $\cD_j$ for some $j\in [4t]$, we \emph{kill} (and remove from $\Lambda$) all the indices $h \in [\ell]$ such that $V(\cD_j)$ intersects $X_h$ in more than $2Cn^{1/6}$ vertices. We also kill any index $h$ such that $V(\cD_j)$ intersects $S_h$. Due to our conditions \eqref{cond:2} and \eqref{cond:3}, we have that throughout the process, 
\[|\Lambda|\geq \ell- 4t(C+Cn^{1/2})\geq \left(1-\frac{\alpha}{2}\right)\ell, \]
for $n$ large recalling that $\ell=\alpha n^{2/3}$ and $t=\alpha n^{1/8}$.  Moreover, due to condition \eqref{cond:3}, at any point in the process, for all
alive indices $h$ in $\Lambda$, the size of $U_h\subset X_h$ is at least 
\[|U_h|\geq |X_h|- \sum_j |V(\cD_j) \cap X_h| \geq 2 n^{1/3}-8tCn^{1/6}\geq n^{1/3},\]
for $n$ sufficiently large. We remark that it is crucial in the previous two calculations that $t=n^{1/8}$ and so when choosing our diamond trees, we do not kill too many indices or make too many of the sets $X_h$ too small to be used by subsequent diamond trees. In fact any $t$ polynomially smaller than $n^{1/6}$ would suffice for this.

So let us suppose that we are at step $j^*\in[4t]$ where we look to find $\cD_{j^*}$ and we have some fixed set $\Lambda$ of alive indices and subsets $U_h\subset X_h$ for $h\in \Lambda$. We run a sub-algorithm that finds $\cD_{j^*}$ in two phases. 
We begin by setting $\Gamma=\Lambda$ and $\bc{C}=\emptyset$. The first phase of the sub-algorithm works by finding at most $C$  small order diamond trees whose removable vertices intersect many of the $U_h$ for $h\in\Lambda$. The family  $\bc{C}$ will collect these small order  diamond trees and the set $\Gamma$ will keep track of the indices $h$ in $\Lambda$ for  which we have not yet intersected $U_h$. In the second phase of the algorithm, we will form $\cD_{j^*}$ by joining together the  diamond trees in $\bc{C}$ so that they form one diamond tree. By guaranteeing that our diamond trees in $\bc{C}$ have removable vertices that intersect most of the sets $U_h$, we will guarantee condition \eqref{cond:1} of the claim.  Before starting, we also initiate by setting $W'\subset W$ to be 
\[W'=W\setminus \left(\cup_{h\in[\ell]}(S_h\cup X_h) \cup_{j<j^*}V(\cD_j)\right).\]
In words, $W'$ is the subset of vertices of $W$ that has not been used in any of the structures that we have found so far. Finally we initiate a counter by setting $s=1$. 

\vspace{2mm}

At step $s$, we apply Lemma \ref{lem:diamond tree that intersects many sets} on the sets $W'$ and $\{U_h:h\in \Gamma\}$. We thus find a $K_r$-diamond tree $\cC_s=(T,R,\Sigma)$ which we add to $\bc{C}$,
which has the following properties guaranteed by Lemma \ref{lem:diamond tree that intersects many sets}:
\begin{enumerate}[label=(\roman*)]
\item \label{cond:sigma} $\Sigma\subset K_{r-1}(G[W'])$ and we delete $V(\Sigma)$ from $W'$;
\item \label{cond:high intesects} $R\subset\cup_{h\in\Gamma} U_h$ and defining $\Gamma_s\subset \Gamma$ to be $\Gamma_s:=\{h':R\cap U_{h'}\neq \emptyset\}$, we have that $|\Gamma_s|\geq |\Gamma|/4r$. We delete $\Gamma_s$ from $\Gamma$;

\item \label{cond:not too big} The order of $\cC_s$ is at most $n^{2/3}$;
\item  \label{cond:not too many big intersects} There is a set $\Phi_s\subset \Gamma_s \subset  \Gamma\subset [\ell]$ of at most $n^{1/2}$ indices, such that for all $h\in [\ell]\setminus \Phi_s$ we have that $|V(\cC_s)\cap U_h|\leq n^{1/6}$. 
\end{enumerate}
Finding such a $\cC_s$ concludes this step $s$. If $|\Gamma|<\alpha \ell/2$, we terminate this phase and move on to the next phase. If $|\Gamma|\geq \alpha \ell/2$, we move to step $s+1$. 

We must check that the conditions for Lemma~\ref{lem:diamond tree that intersects many sets} are satisfied throughout this phase in order to find the required diamond trees $\cC_s$ at each step.  Indeed this follows because 
\[\frac{\alpha^2}{2}n^{2/3}=\frac{\alpha}{2}\ell\leq |\Gamma|\leq \ell =\alpha n^{2/3}\]
throughout 
and we have that $|U_h|\geq n^{1/3}$ for all $h\in \Gamma$ as $\Gamma\subset \Lambda$ is a subset of alive indices. Finally we have $|W'|\geq n/4$ throughout this  process. Indeed, note that due to condition \ref{cond:high intesects} and the fact that we only continue until $|\Gamma|\leq \alpha \ell/2$, we have that the process runs for a maximum of $C$ steps, recalling the definition of $C$ from condition \eqref{cond:2} of the claim. That is, we have that $|\bc{C}|\leq C$ throughout and  so 
\begin{align}
 \nonumber   |W'|&\geq|W|- \sum_{h\in[\ell]}(|S_h|+ |X_h|) -\sum_{j<j^*}|V(\cD_j)|-\sum_{\cC\in \bc{C}}V(\cC) \\ \nonumber
    &\geq \frac{n}{2}-\ell\cdot 3n^{1/3}-8trM-Cn^{2/3} \\ \label{eq:W' lower bound}
    &\geq \left(\frac{1}{2}-(4+8r)\alpha\right)n \geq \frac{n}{4},
\end{align}
due to our upper bound on $\alpha$, for $n$ sufficiently large. This verifies that we find $\cC_s$ at every step $s$ of this process and so we finish this phase with $|\Gamma|<\alpha \ell/2$ and some family $\bc{C}=\{\cC_1,\ldots,\cC_c\}$ of  $c\leq C$ vertex disjoint $K_r$-diamond trees.

\vspace{2mm}

Now we describe how we generate  $\cD_{j^*}$ which will have all the diamond trees $\cC_s\in \bc{C}$ as sub-diamond trees. We refer the reader to Figure~\ref{fig:Dj0} to keep on track of the many components that contribute to our $\cD_{j^*}$.  One thing to note is that the sum of the orders of the diamond trees in $\bc{C}$ is far too small for us to just build $\cD_{j^*}$ from the diamond trees in $\bc{C}$. Indeed the sum of the orders is $O(n^{2/3})$ and we want $\cD_{j^*}$ to have order  $M=n^{7/8}$. Therefore we will have to find the majority of the  $K_r$-diamond tree $\cD_{j^*}$ elsewhere. In order to prepare for this, we first split $W'$ arbitrarily into $U_0,W_0$ and $Z_0$ of roughly equal size  and note that due to our lower bound \eqref{eq:W' lower bound} on $|W'|$, we have that each of these sets has size at least $n/16$. Next we fix $d_*:=\alpha^2p^{r-1}n$ and $z=\alpha^2 n$ and apply Proposition~\ref{prop: choosing removable vertices} with respect to the sets $U_0$ and $W_0$ to get disjoint sets $X,Y\subset U_0$ as detailed there. Note that $|X|\leq 2n^{2/3}$. Indeed if $|X|>1$, then we have that $|X| \leq 2z/d_*\leq 2p^{1-r}\leq 2n^{2/3}$ due to Fact~\ref{fact:dense}.

\begin{figure}
    \centering
  \includegraphics[scale=0.84]{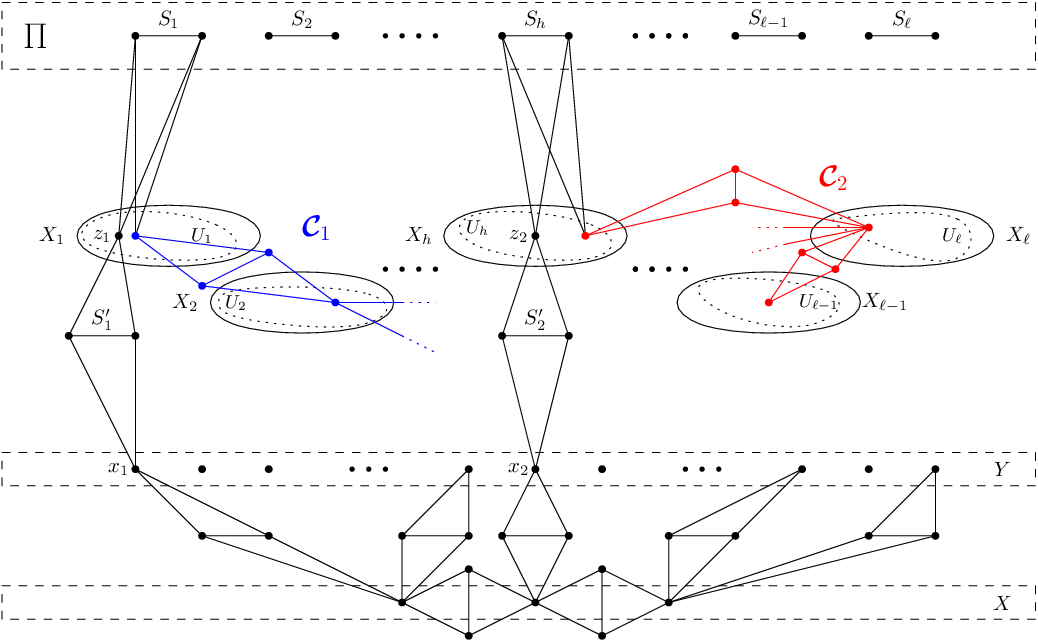}
    \caption{   \label{fig:Dj0} An example of $\cD_{j^*}$ and its components. In this case, we have $c=2$, $h_1=1$ and $h_2=h$.} 
  \end{figure}

Now for $1\leq s\leq c$, define $Z_s:=\cup_{h\in\Gamma_s\setminus \Phi_s}(U_h\setminus V(\cC_s))$.  In words, $Z_s$ is the union of the sets $U_h$ which $\cC_s$ intersects, after removing  the sets $U_{h'}$ which $\cC_s$ intersects in too many vertices and then removing the  vertices of $\cC_s$. Now we have that for each $s\in [c]$, 
\[|Z_s|\geq (|\Gamma_s|-|\Phi_s|)(n^{1/3}-n^{1/6})\geq \frac{\alpha \ell n^{1/3}}{8r}-2n^{5/6}\geq \frac{\alpha^2n}{16r}\geq \alpha'n,\]
for $n$ large, as the $U_h$ are pairwise disjoint.  Note also that as the $\Gamma_s$ are pairwise disjoint, the $Z_s$ are also pairwise disjoint for $s\in[c]$. Now for $1\leq s\leq c$, apply Corollary \ref{cor:transversalcliques} \eqref{cor:popular K_r-1}  to find an $(r-1)$-clique $S'_s\subset K_{r-1}(G[Z_0])$ such that there is a vertex $z_s\in Z_s\cap N^G(S'_s)$ and a vertex $x_s\in (X\cup Y)\cap N^G(S'_s)$. We delete the $S'_s$ from $Z_0$ and move to the next index $s+1$ or finish if $s=c$.

Now choose some $Y'\subset Y$ such that $x_s\in X\cup Y'$ for all $s\in [c]$ and \[|Y'|+|X|+\sum_{s\in[c]}(|R_{\cC_{s}}|+1)=M.\] This is easily done as~$|X|+|Y|=\alpha^2 n$ is linear and $|X|,|R_{\cC_s}|\leq 2n^{2/3}$ for all $s\in [c]$ which is much smaller than $M=n^{7/8}$. By Proposition \ref{prop: choosing removable vertices}, there is a $K_r$-diamond tree $\tilde{\cD}=(\tilde{T},\tilde{R},\tilde{\Sigma})$ with $\tilde{R}=X\cup Y'$ and $\tilde{\Sigma}\subset K_{r-1}(G[W_0])$ a matching of $(r-1)$-cliques in $W_0\subset W'$. Our diamond tree $\cD_{j^*}$ is then obtained by connecting $\tilde{\cD}$ and all the $\cC_s\in \bc{C}$. 
In more detail, for each $s\in [c]$, there exists some $h_s\in\Gamma_s$ such that $z_s\in U_{h_s}\subset X_{h_s}$. We define
\[R_{j^*}:=\tilde{R}\cup_{s\in[c]}(R_{\cC_s} \cup \{z_s\}) \qquad \mbox{and} \qquad \Sigma_{j^*}:=\tilde{\Sigma}\cup_{s\in[c]}(\Sigma_{\cC_s}\cup \{S'_s\} \cup \{S_{h_s}\}  ),\] where $\Sigma_{\cC_s}$ is the set of interior $(r-1)$-cliques of $\cC_s$, we have that $S'_{s}\in K_{r-1}(G[Z_0])$ is the $(r-1)$-clique which forms a clique with both $z_s$ and $x_s$ defined above and $S_{h_s}$ is the $(r-1)$-clique corresponding to the set $X_{h_s}$ (which contains $z_s$) in Claim~\ref{clm:candidateedges}. 
We claim that there exists a diamond tree $\cD_{j^*}$ of order $M$ which has $R_{j^*}$ as a set of removable vertices and $\Sigma_{j^*}$ as a set of interior $(r-1)$-cliques. Indeed we can form the defining auxiliary tree $T_{j^*}$ by starting with the forest of the disjoint union of $\tilde{T}$ and the $T_{\cC_{s}}$ for $s\in [c]$, where $T_{\cC_{s}}$ denotes the defining tree for the $K_r$-diamond tree $\cC_s$. 
For each $s\in[c]$, we then add a path of length two (with two edges)  between some vertex in $V(T_{\cC_s})$ and $V(\tilde{T})$. The edges of this path correspond exactly to the internal $(r-1)$-cliques $S_{h_s}$ and $S'_{s}$ and thus the vertices of this path correspond to $x_s$, $z_s$ and some vertex in $R_{\cC_s}\cap U_{h_s}$ for each $s\in[c]$. 

\vspace{2mm}

This thus defines $\cD_{j^*}$ and so we update all the $U_h$ to be $U_h\setminus V(\cD_{j^*})$ for $h\in[\ell]$ and kill any indices $h\in\Lambda$ such that either $V(\cD_{j^*})$ intersects $S_h$ or $|X_h\cap V(\cD_{j^*})|\geq 2Cn^{1/6}$. We now need to check that the conditions \eqref{cond:1}, \eqref{cond:2} and \eqref{cond:3} hold for $\cD_{j^*}$. To see \eqref{cond:1}, note that $R_{j^*}$ contains all the $R_{\cC_{s}}$ for $s\in [c]$ and so intersects $X_h$ for all $h\in \cup_{s\in[c]}\Gamma_s$. Moreover taking $\Gamma$ as defined at the end of finding the $\cC_s$, we have that $|\Gamma\cup_{s\in[c]}\Gamma_s|\geq (1-\alpha/2)\ell$ and $|\Gamma|\leq \alpha \ell/2$ and so this confirms \eqref{cond:1}. To see \eqref{cond:2}, note that the only times we used vertices of the $S_h$ with $h\in[\ell]$ to construct $\cD_{j^*}$ was when we added the $S_{h_s}$ for $s\in[c]$ to the set of interior cliques. Thus we intersected exactly $c\leq C$ of these with $V(\cD_{j^*})$. Finally we have that \eqref{cond:3} for $\cD_{j^*}$ is implied by the conditions \ref{cond:not too many big intersects} when we found the $\cC_s$. Indeed, we have that $R_{j^*}\cap (\cup_{h\in[\ell]}X_h)=\cup_{s\in [c]}(R_{\cC_s}\cup\{z_s\})$ and so for any  index $h$ that does not lie in $\cup_{s\in [c]}\Lambda_{s}$ (which has size at most $Cn^{1/2}$),  we have that \[|V(\cD_{j^*})\cap X_h|\leq \sum_{s\in [c]}|(V(\cC_{s})\cup\{z_s\})\cap X_h|+\leq C(n^{1/6}+1)\leq 2Cn^{1/6}.\] This concludes the finding of $\cD_{j^*}$ and doing this for all $j^*\in[4t]$ gives the desired claim and hence the proposition. 
\end{altclaimproof}

\end{proof}

\section{Concluding Remarks}\label{sec:conclude}
In this paper, we showed that a condition of $\beta=o( p^{k}n)$ in an $n$ vertex $(p,\beta)$-bijumbled graph guarantees a $K_{k+1}$-factor. We conjecture that the same condition in fact guarantees any subgraph with maximum degree $k$. 

\begin{conj} \label{conj:universal}
For any $2\leq k\in\NN$ and $c>0$ there exists an $\eps>0$ such that any $n$-vertex $(p,\beta)$-bijumbled graph with $\delta(G)\geq cpn$ and  $\beta \leq \eps p^{k}n$ is~\emph{$k$-universal}. That is, given any graph $F$ on at most $n$ vertices, with maximum degree at most $k$, $G$ contains a copy of $F$.
\end{conj}

Note that Corollary~\ref{cor:2-factor} settles Conjecture~\ref{conj:universal} for the case $k=2$.  For $k\geq 3$, the best known result comes from the sparse blow-up lemma of Allen, B\"ottcher, H\`an, Kohayakawa and Person \cite{ABHKP16} which gives a condition of $\beta=o(p^{(3k+1)/2}n)$ guaranteeing $k$-universality in a $(p,\beta)$-bijumbled  graph. 

\vspace{2mm}

The conjecture echoes the notion that a $K_{k+1}$-factor is the `\emph{hardest}' maximum degree $k$ graph to find. This idea has manifested in various other settings. For example, we know from the theorem of Hajnal and Szemer\'edi (Theorem~\ref{thm:HajnalSzem}) that any $n$ vertex graph $G$ with $\delta(G)\geq (k/(k+1))n$ contains a $K_{k+1}$-factor and that this is tight. Bollob\'as and Eldridge~\cite{bollobas1978packings}, and independently Catlin~\cite{Catlin}, conjectured that the same minimum degree condition actually guarantees $k$-universality. This has been proven for $k=2,3$~\cite{aigner1993embedding,alon19962,csaba2003proof} but remains open in general. In the case of random graphs, Johansson, Kahn and Vu~\cite{JohanssonKahnVu_FactorsInRandomGraphs} proved that the threshold for the appearance of a $K_{k+1}$-factor is of the order of \[p_k^*(n):=n^{-2/(k+1)}(\log n)^{2/(k^2+k)}.\] A recent breakthrough result of Frankston, Kahn, Narayanan and Park~\cite{frankston2019thresholds} implies that for any $n$ vertex graph $F$ with maximum degree $k$, the threshold for the appearance of $F$ in $G(n,p)$ is at most $p_k^*(n)$. Note that this is \emph{not} implying that $G(n,p)$ is $k$-universal whp when $p=\omega(p^*_k(n))$ as we can only guarantee that some fixed $F$ appears whp. However, the stronger version that $p^*_k(n)$ is the threshold for $k$-universality is believed to be true but only verified for $k=2$~\cite{FKL16}.

\vspace{2mm}

One thing that sets aside the pseudorandom setting in stark contrast to the other settings discussed above is that it might be possible to replace a $K_{k+1}$-factor as the benchmark for the `hardest' graph to find in the host graph, by a single copy of $K_{k+1}$. Indeed, various authors~\cite{CFZ12,fox2013chromatic,KS06,sudakov2005generalization} have stipulated that $n$ vertex $K_{k+1}$-free $(p,\beta)$-bijumbled graphs exist with $\beta=\Theta(p^{k}n)$.  Such graphs would witness the tightness of both Theorem~\ref{thm:main} and  Conjecture~\ref{conj:universal}   for all values of  $k\geq 2$ (taking $r=k+1$ in the setting of  Theorem~\ref{thm:main}).  Focusing on optimally pseudorandom graphs (that is, fixing $\beta=\Theta(\sqrt{pn})$ in $(p,\beta)$-bijumbled graphs), we expect to be able to  find $K_{k+1}$-free optimally pseudorandom graphs with $p=\Omega(n^{-1/(2k-1)})$. These are only known to exist when $k=2$. Indeed, we discussed the triangle-free construction of Alon in the introduction and other  constructions~\cite{conlon2017sequence,Kopparty11} have also been given which are (near-)optimal. For $k\geq 3$ however, this remains a key challenge in the understanding of pseudorandom graphs with the best known general construction coming from a recent  improvement of Bishnoi, Ihringer and Pepe~\cite{bishnoi2020construction} who give $K_{k+1}$-free optimally pseudorandom graphs of density $p=\Theta(n^{-1/k})$. Further interest in finding denser such graphs comes from a recent remarkable  connection discovered by Mubayi and Verstra\"ete~\cite{mubayi2019note} that shows that if, as we expect, the  $K_{k+1}$-free optimally pseudorandom graphs with density $p=\Omega(n^{-1/(2k-1)})$  do exist, then it is possible to improve the lower bound on the  off-diagonal Ramsey numbers to match the upper bound and thus determine the asymptotics of this extremal function. In detail, they show that if these pseudorandom graphs exist, then the off-diagonal Ramsey number is $R(k+1,t)=t^{k+o(1)}$ as $t$ tends to infinity. In fact, even a construction with $p=\omega(n^{-1/{(k+1)}})$ would improve on the current best known lower bound on off-diagonal Ramsey numbers due to Bohman and Keevash~\cite{bohman2010early}.

\vspace{2mm}

We conclude by noting that Theorem~\ref{thm:maintriangle} is, in some sense, the first result of its kind, giving a tight condition on pseudorandomness to guarantee the existence of a spanning structure. Indeed, the case of Hamilton cycles remains an intriguing open problem.  Krivelevich and Sudakov~\cite{krivelevich2003sparse} conjectured that a condition of $\lambda=o(d)$ is sufficient in $(n,d,\lambda)$-graphs and proved the currently best known bound of 
\[\lambda=o\left(\frac{(\log\log n)^2d}{\log n(\log \log \log n)}\right).\]
For hypergraphs of higher uniformity, one can easily generalise the notion of bijumbledness in Definition~\ref{def:bijumbled} but the picture becomes considerably more complex. Indeed it turns out that the only subgraphs that one can guarantee by imposing conditions on bijumbledness are \emph{linear} subgraphs, those in which pairs of  hyperedges intersect in at most one vertex. Building on previous work~\cite{conlon2012weak,kohayakawa2010weak,lenz2016hamilton,lenz2016perfect} mainly concerned with dense hypergraphs (the so-called quasirandom regime), Hi\d{\^{e}}p H\`an, Jie Han and the author~\cite{han2020factors,han2020factorsjournal} recently gave the best-known conditions on pseudorandomness  that guarantee  different linear subgraphs of hypergraphs. These include all fixed sized linear subgraphs as well as $F$-factors for linear $F$ (including perfect matchings) and loose Hamilton cycles. The tightness of these results is unclear as no good constructions are known for $F$-free pseudorandom hypergraphs. In general, the appearance of subgraphs in sparse pseudorandom (hyper-)graphs remains a fascinating area which is far from being understood.

\bibliography{Biblio}

\end{document}